\documentclass[french,english]{paper}
\usepackage[T1]{fontenc}
\usepackage[latin9]{inputenc}
\usepackage{geometry}
\usepackage{babel}
\makeatletter
\addto\extrasfrench{%
   \providecommand{\fg}{\ifdim\lastskip>\z@\unskip\fi~\frqq}%
}

\makeatother
\usepackage{amsmath}
\usepackage{amssymb}
\usepackage{graphicx}
\usepackage{setspace}
\usepackage{esint}
\usepackage[unicode=true,pdfusetitle,
 bookmarks=true,bookmarksnumbered=false,bookmarksopen=false,
 breaklinks=false,pdfborder={0 0 1},backref=false,colorlinks=false]
 {hyperref}

\makeatletter

\newcommand{\noun}[1]{\textsc{#1}}
\providecommand{\tabularnewline}{\\}

\numberwithin{figure}{section}
\numberwithin{equation}{section}
\numberwithin{table}{section}
\newcommand{\lyxaddress}[1]{
\par {\raggedright #1
\vspace{1.4em}
\noindent\par}
}

\usepackage[title,titletoc]{appendix}
\pdfoutput=1

\makeatother

\begin{document}

\title{{\Large{Finite size effects in the correlation structure of stochastic
neural networks: analysis of different connectivity matrices and failure
of the mean-field theory}}}

\author{Diego Fasoli{*}, Olivier Faugeras{*}}

\maketitle

\lyxaddress{{*}NeuroMathComp Laboratory, INRIA Sophia-Antipolis, France. Email:
firstname.name@inria.fr}

\section*{Abstract}

\noindent We quantify the finite size effects in a stochastic network
made up of rate neurons, for several kinds of recurrent connectivity
matrices.

\noindent This analysis is performed by means of a perturbative expansion
of the neural equations, where the perturbative parameters are the
intensities of the sources of randomness in the system.

\noindent In detail, these parameters are the variances of the background
or input noise, of the initial conditions and of the distribution
of the synaptic weights.

\noindent The technique developed in this article can be used to study
systems which are invariant under the exchange of the neural indices
and it allows us to quantify the correlation structure of the network,
in terms of pairwise and higher order correlations between the neurons.

\noindent We also determine the relation between the correlation and
the external input of the network, showing that strong signals coming
from the environment reduce significantly the amount of correlation
between the neurons.

\noindent Moreover we prove that in general the phenomenon of propagation
of chaos does not occur, even in the thermodynamic limit, due to the
correlation structure of the $3$ sources of randomness considered
in the model.

\noindent Furthermore, we show that the propagation of chaos does
not depend only on the number of neurons in the network, but also
and mainly on the number of incoming connections per neuron.

\noindent To conclude, we prove that for special values of the parameters
of the system the neurons become perfectly correlated, a phenomenon
that we have called \textit{stochastic synchronization}.

\noindent These discoveries clearly prevent the use of the mean-field
theory in the description of the neural network.

\section{\label{sec:Introduction} Introduction}

\noindent According to the theory of complexity \cite{Weaver48}\cite{bar1997dynamic}\cite{h2006information},
in order to have emergent behaviors in systems made up of many interacting
units, what really matters are not the properties of the single units
themselves, but rather the way they interact with each other.

\noindent The most famous example of this phenomenon is represented
by flocks of birds, since it is possible to recreate in a computer
simulation their ability to form stable and complicated patterns and
to rejoin when the group is splitted, throught the implementation
of very simple rules of interaction.

\noindent In fact, according to the famous artificial life program
\textit{Boids} \cite{Reynolds:1987:FHS:37401.37406}, it is possible
to reproduce this emergent behavior assuming that every bird has to
fly in the same direction of the neighbours, with the same speed and
avoiding obstacles or to bump into other birds.

\noindent This example clearly shows that the complexity of the system
is a consequence of the interaction between birds, and not of the
model used to described a single bird.

\noindent Therefore, in the context of the brain, where the elementary
units are represented by neurons, the priority is not to study extremely
biologically realistic models of single neurons.

\noindent The really important problem is to describe in an accurate
way the interaction between them, or in other terms their synaptic
connectivity matrix.

\noindent For this reason we believe that it should be more relevant
to use simplified neural models (like the so called \textit{rate model}
\cite{Wilson72excitatoryand}\cite{WilsonCowan}\cite{Amari_1972}\cite{PhD_amari1977}\cite{journals/siamads/TouboulHF12})
with complex connectivity matrices, than to use more biologically
plausible neural models (like the \textit{Hodgkin-Huxley model} \cite{citeulike:851137})
with simple connections.

\noindent So in this article we focus mainly on the differences in
the behavior of the network induced by different topologies of the
synaptic connectivity.

\noindent Following this current of thought, a great effort has been
devoted to finding the pattern of the synaptic connections of the
human brain \cite{10.1371/journal.pcbi.0010042}\cite{sporns2011the}\cite{HumanConnectomeProject}.

\noindent Therefore an important question we have to try to answer
is: how does the behavior of the brain change when we modify its connectivity
matrix?

\noindent This is the key to understand how the brain processes the
information it receives from the environment, and hopefully also the
necessary ingredient to explain how its higher cognitive functions
naturally emerge from the interaction of the neurons.

\noindent It is therefore fundamental to develop a theory that is
able to determine the behavior of a neural network once its external
inputs and connectivity matrix are known, for many different kinds
of input and connection structures.

\bigskip{}

\noindent The behavior of a neural network that this theory must be
able to describe is twofold.

\noindent In fact the theory has to provide the time evolution of
the membrane potentials and firing rates of the neurons but, from
a probabilistic point of view, it must also be able to determine their
statistics.

\noindent It is well known that the neurons are not reliable units
and that the variability of their behavior is described by adding
a source of stochasticity in their equations \cite{rolls2010noisy}\cite{citeulike:732425}\cite{White:TINS:2000}.

\noindent This transforms the system into a non-deterministic network,
and therefore one has to compute the probability density of the neurons
and their pairwise correlation structure.

\noindent Actually the latter is receiving an increasing attention
from the scientific community since it may be involved in the extraordinary
information processing capabilities of the brain.

\noindent It is also known as \textit{functional connectivity}, a
term coined in order to distinguish the correlation structure from
the wiring pattern of all the synaptic connections, that is known
instead as \textit{structural or anatomical connectivity}.

\noindent The problem of finding the relation between these two kinds
of connectivity is currently intensively investigated \cite{citeulike:61}\cite{Ponten_Daffertshofer_Hillebrand_Stam_2010}\cite{10.1371/journal.pcbi.1002059}\cite{koch2002investigation}\cite{Eickhoff:9544}\cite{journals/neuroimage/CabralHKD12},
and the theory we are looking for should provide such a link.

\bigskip{}

\noindent In Section \ref{sec:Description of the model} we develop
a perturbative approch that let us determine the behavior of a neural
network with finite size, made up of a generic number of neurons described
by rate equations.

\noindent The perturbative parameters are the variances of the $3$
sources of randomness considered in the model: the background or input
noise, the initial conditions and the synaptic weights.

\noindent This allows us to study the network for many different types
of structural connectivity matrices, with a special emphasis on the
evaluation of their corresponding functional connectivity matrices.

\noindent Therefore in Section \ref{sec:Correlation structure of the network}
we provide an explicit formula that relates the two connectivity matrices,
for many different kinds of structural connectivities, and we also
explain how to use this model to compute higher order correlations,
like those between triplets and quadruplets of neurons.

\noindent In particular, in Section \ref{sec:Calculation of the fundamental matrix}
we consider the special cases of block circulant matrices with circulant
blocks, and of symmetric matrices.

\noindent The goodness of the perturbative approach is proved in Section
\ref{sec:Numerical comparison} through many numerical simulations.

\noindent Moreover, using these formulae, in Section \ref{sec:Correlation as a function of the input}
we determine the relation between the correlation structure and the
external input of the network, showing that strong signals coming
from the environment reduce significantly the amount of correlation
between the neurons.

\noindent Moreover in Section \ref{sec:Failure of the mean field theory}
we prove that the phenomenon known as \textit{propagation of chaos}
\cite{journals/siamads/TouboulHF12}\cite{baladron:inserm-00732288}\cite{samuelides:inria-00529560}\cite{Touboul11}
in general does not occur, even in the thermodynamic limit, due to
the correlation structure of the $3$ sources of randomness.

\noindent Furthermore, we show that propagation of chaos does not
depend only on the number of neurons in the network, but also and
mainly on the number of incoming connections per neuron.

\noindent This model predicts also arbitrarily high values of correlation
between the neurons for special values of the parameters of the system,
a phenomenon that we have called \textit{stochastic synchronization}. 

\noindent The direct consequence of these results is therefore the
impossibility to apply in general the mean-field theory in order to
describe the activity of the network.

\bigskip{}

\noindent It is also interesting to observe that since this approach
works for a generic and finite number of neurons $N$, it is able
to quantify the finite size effects of the network.

\noindent Therefore simply increasing the value of the parameter $N$,
in principle we can evaluate the differences that occur in the behavior
of a neural network when we switch from a microscopic scale ($N\sim10^{0}$,
i.e. single neurons) to mesoscopic ($N\sim10^{1}\div10^{5}$, i.e.
neural masses and cortical columns) and then to macroscopic scales
($N\sim10^{6}\div10^{11}$, i.e. extended brain areas, like the visual
cortex).

\noindent This in principle would allow us to show if there are actually
emerging properties of the system that are triggered by its size.

\noindent Such properties can be identified for example in a difference
in the information processing capabilities of the system that emerges
when we increase $N$.

\noindent These capabilities can be evaluated with our model in the
linear approximation regime, since in that case the system is described
by a multivariate normal distribution and therefore it allows us to
compute analytically all the information quantities of the system.

\noindent However we will not show this analysis here since it is
beyond the purpose of the current article and also because the study
of macroscopic areas of the brain requires the refinement of the model
using neural fields equations \cite{PhD_amari1977}\cite{Coombes:2005:WBP:1081979.1081982}\cite{veltz:hal-00784425}
in order to describe the spatial extension of the areas and also the
delays in the propagation of the electric signals.

\section{\noindent \label{sec:Description of the model}Description of the
model}

\begin{flushleft}
We suppose that the neural network is described by the following rate
model:
\par\end{flushleft}

\begin{onehalfspace}
\begin{center}
{\small{
\begin{equation}
dV_{i}\left(t\right)=\left[-\frac{1}{\tau}V_{i}\left(t\right)+\sum_{j=0}^{N-1}J_{ij}\left(t\right)S\left(V_{j}\left(t\right)\right)+I_{i}\left(t\right)\right]dt+\sigma_{1}dB_{i}\left(t\right)\label{eq:rate-model-exact-equations-2}
\end{equation}
}}
\par\end{center}{\small \par}
\end{onehalfspace}

\noindent with $i=0,1,...,N-1$.

\noindent Here $V_{i}\left(t\right)$ is the mebrane potential of
the $i$-th neuron, $I_{i}\left(t\right)$ is its external input current
and $\tau$ is a time constant that determines the speed of convergence
of the membrane potential to its rest state $V_{i}\left(t\right)=0$
in the case of disconnected neurons.

\noindent Moreover $N$ is the number of neurons in the network, $J_{ij}\left(t\right)$
is the synaptic weight from the the $j$-th neuron to the $i$-th
neuron, and $S\left(\cdot\right)$ is a sigmoid function, which converts
the membrane potential of a neuron into the rate of the spikes it
produces, according to the law:

\noindent \begin{flushleft}
\[
S\left(V\right)=\frac{T_{MAX}}{1+e^{-\lambda\left(V-V_{T}\right)}}
\]

\par\end{flushleft}

\noindent \begin{flushleft}
where $T_{MAX}$ is the maximum amplitude of the function, $\lambda$
is a parameter that determines its slope, and $V_{T}$ is the horizontal
shift along the $V$ axis.
\par\end{flushleft}

\bigskip{}

\noindent \begin{flushleft}
$\sigma_{1}$ is the noise intensity, that for simplicity is supposed
to be the same for all the neurons and constant in time.
\par\end{flushleft}

\noindent The functions $B_{i}\left(t\right)$ are Brownian motions,
which can be equivalently interpreted as a background noise for the
membrane potentials $V_{i}\left(t\right)$ or as the stochastic component
of the external input $I_{i}\left(t\right)$.

\noindent In general they are correlated according to a covariance
matrix $\Sigma_{1}$, whose components are:

\begin{onehalfspace}
\begin{center}
{\small{
\begin{align}
\left[\Sigma_{1}\right]_{ij}= & Cov\left(\frac{dB_{i}\left(t\right)}{dt},\frac{dB_{j}\left(s\right)}{ds}\right)=C_{ij}^{1}\delta\left(t-s\right)\nonumber \\
\label{eq:Brownian-covariance}\\
C_{ij}^{1}= & \begin{cases}
1 & \begin{array}{ccc}
\mathrm{if} &  & i=j\end{array}\\
\\
C_{1} & \begin{array}{ccc}
\mathrm{if} &  & i\neq j\end{array}
\end{cases}\nonumber 
\end{align}
}}
\par\end{center}{\small \par}
\end{onehalfspace}

\noindent where $\delta\left(\cdot\right)$ is the Dirac delta function,
while $C_{1}$ represents the correlation between two different Brownian
motions (here the derivative of the Brownian motion is meant in the
weak sense of distributions and is interpreted as white noise). In
order to be a true covariance matrix, $\Sigma_{1}$ must be positive-semidefinite.
Since it is symmetric, then it is positive-semidefinite if and only
if its eigenvalues are non-negative. But $\Sigma_{1}$ is a circulant
matrix, therefore its eigenvalues are $e_{0}=1+C_{1}\left(N-1\right)$
and $e_{i}=1-C_{1}$, for $i=1,2,...,N-1$. Therefore $\Sigma_{1}$
is positive-semidefinite if and only if $\frac{1}{1-N}\leq C_{1}\leq1$.
We could increase the complexity of this correlation structure, since
there is no technical difficulty in doing that, but we keep it simple
for the sake of clarity.

\bigskip{}

\noindent We also suppose that the initial conditions are distributed
according to the following multivariate normal probability density:

\begin{onehalfspace}
\begin{center}
{\small{
\begin{equation}
\overrightarrow{V}\left(0\right)\sim\mathcal{N}\left(\overrightarrow{\mu},\Sigma_{2}\right)\label{eq:initial-conditions}
\end{equation}
}}
\par\end{center}{\small \par}
\end{onehalfspace}

\noindent where for simplicity:

\begin{onehalfspace}
\begin{center}
{\small{
\begin{align}
\mu_{i}= & \begin{array}{ccc}
\mu, &  & i=0,1,...,N-1\end{array}\label{eq:initial-conditions-mean}\\
\nonumber \\
\Sigma_{2}= & \sigma_{2}^{2}\left[\begin{array}{cccc}
1 & C_{2} & \cdots & C_{2}\\
C_{2} & 1 & \cdots & C_{2}\\
\vdots & \vdots & \ddots & \vdots\\
C_{2} & C_{2} & \cdots & 1
\end{array}\right]\label{eq:initial-conditions-covariance}
\end{align}
}}
\par\end{center}{\small \par}
\end{onehalfspace}

\noindent Here $\sigma_{2}$ represents the initial standard deviation
of each neuron, while $C_{2}$ is the initial correlation between
pairs of neurons. As before, the matrix $\Sigma_{2}$ must be positive-semidefinite,
and this is true if and only if $\frac{1}{1-N}\leq C_{2}\leq1$. Again,
we could increase the complexity of this correlation structure, if
desired.

\bigskip{}

\noindent For the synaptic connectivity matrix $J\left(t\right)$,
we suppose that its entries have a deterministic temporal evolution,
but they are distributed randomly over many repetitions of the network.
With this model we can match our results with those appeared in \cite{10.3389/neuro.10.001.2009}.
So, in detail, we suppoose that the synaptic connectivity matrix $J\left(t\right)$
has random entries distributed according to the law:

\begin{onehalfspace}
\begin{center}
{\small{
\begin{equation}
J\left(t\right)\sim\mathcal{MN}\left(\overline{J}+\sigma_{4}Z\left(t\right),\Omega_{3},\Sigma_{3}\right)\label{eq:synaptic-weights-1}
\end{equation}
}}
\par\end{center}{\small \par}
\end{onehalfspace}

\noindent This is the so called \textit{matrix normal distribution}
\cite{MatrixNormalDistribution}, namely the generalization of the
multivariate normal distribution to the case of matrix-valued random
variables. Here $\overline{J}$, $Z\left(t\right)$, $\Omega_{3}$
and $\Sigma_{3}$ are $N\times N$ deterministic matrices. In particular,
$\overline{J}+\sigma_{4}Z\left(t\right)$ represents the mean of $J\left(t\right)$,
while $\Omega_{3}$ and $\Sigma_{3}$ are its covariance matrices.
We suppose that $\overline{J}$ has only two different kinds of entries,
namely $0$ (absence of connection) and $\Lambda$, where $\Lambda$
is a free non-zero parameter. We also suppose that $Z\left(t\right)$
has general entries (with the obvious exception that $Z_{ij}\left(t\right)=0$
if there is no connection from the $j$-th neuron to the $i$-th neuron,
namely if $J_{ij}\left(t\right)=0$ $\forall t$), therefore it is
a source of inhomogeneity and time-variability for the connectivity
matrix. We use for simplicity specific structures of the covariance
matrices $\Omega_{3}$ and $\Sigma_{3}$. Supposing that all the non-zero
entries of $J\left(t\right)$ have the same standard deviation $\sigma_{3}$,
it is possible to rewrite the matrix $J\left(t\right)$ in the following
equivalent way:

\begin{onehalfspace}
\begin{center}
{\small{
\begin{align}
J\left(t\right)= & \overline{J}+\sigma_{3}W+\sigma_{4}Z\left(t\right)\label{eq:synaptic-weights-2}\\
\nonumber \\
W\sim & \mathcal{MN}\left(0,\widetilde{\Omega}_{3},\widetilde{\Sigma}_{3}\right)\label{eq:W-distribution}
\end{align}
}}
\par\end{center}{\small \par}
\end{onehalfspace}

\noindent where $\widetilde{\Omega}_{3}$ and $\widetilde{\Sigma}_{3}$
are normalized covariance matrices. Their explicit structure is not
important, and the only thing that we need to know is that they are
chosen in order to have:

\begin{onehalfspace}
\begin{center}
{\small{
\begin{equation}
Cov\left(W_{ij},W_{kl}\right)=\begin{cases}
0 & \begin{array}{ccc}
\mathrm{if} &  & \left(g\left(i,j\right)=0\right)\vee\left(g\left(k,l\right)=0\right)\end{array}\\
\\
1 & \begin{array}{ccc}
\mathrm{if} &  & \left(i=k\right)\wedge\left(j=l\right)\wedge\left(g\left(i,j\right)=1\right)\end{array}\\
\\
C_{3} & \mathrm{otherwise}
\end{cases}\label{eq:synaptic-weights-covariance}
\end{equation}
}}
\par\end{center}{\small \par}
\end{onehalfspace}

\noindent where $g\left(x,y\right)=0$ if there is no synaptic connection
from the the $y$-th neuron to the $x$-th neuron (namely if $J_{xy}\left(t\right)=0$
$\forall t$), and $1$ otherwise, while $C_{3}$ is the correlation
between two different and non-zero synaptic weights. We observe that
the range of the possible values of $C_{3}$ in general depends on
the topology of the connectivity matrix, and $W_{ij}=0$ if there
is no connection from the $j$-th neuron to the $i$-th neuron (this
is a consequence of the formulae \ref{eq:synaptic-weights-2} and
\ref{eq:synaptic-weights-covariance}). Again, as for the Brownian
motions and the initial conditions, we could increase the complexity
of this correlation structure, if desired.

\noindent We suppose that every neuron has \textit{the same number
of incoming connections}, that we call $M$. We observe that our assumptions
imply that \textit{the network is invariant under exchange of the
neuronal indices}, which is the main hypothesis of this article. When
$M$ increases, each neuron receives a larger and larger input from
the remainder of the network, therefore in order to fix this divergence
we normalize the synaptic weight in the following way:

\begin{onehalfspace}
\begin{center}
\textit{\small{
\[
J\left(t\right)\rightarrow\frac{J\left(t\right)}{M}
\]
}}
\par\end{center}{\small \par}
\end{onehalfspace}

\noindent This normalization is intended to be used only when $M\neq0$,
because otherwise we obtain $J_{ij}=\frac{0}{0}$. For $M=0$ the
neurons have no incoming connections, therefore we have simply to
set $J_{ij}=0$.

\bigskip{}

\noindent To conclude, we also suppose that the external input current
is deterministic (if we interpret $B_{i}\left(t\right)$ as the noise
of the membrane potential) and given by:

\begin{onehalfspace}
\begin{center}
{\small{
\begin{equation}
\overrightarrow{I}\left(t\right)=\overrightarrow{\overline{I}}+\sigma_{5}\overrightarrow{H}\left(t\right)\label{eq:external-input}
\end{equation}
}}
\par\end{center}{\small \par}
\end{onehalfspace}

\noindent where the vector $\overrightarrow{\overline{I}}$ is time-independent
and such that $\overline{I}_{i}=\overline{I}$, for $i=0,1,...,N-1$.
The vector $\overrightarrow{H}\left(t\right)$ has in general different
and time-variable entries, so it is a source of inhomogeneity and
time-variability.

\bigskip{}

\noindent Now we define the following 2nd order perturbative expansion
of the membrane potential:

\begin{onehalfspace}
\begin{center}
{\small{
\begin{equation}
V_{i}\left(t\right)\approx\mu+\sum_{m=1}^{5}\sigma_{m}Y_{m}^{i}\left(t\right)+{\displaystyle \sum\limits _{\substack{m,n=1\\
m\leq n
}
}^{5}}\sigma_{m}\sigma_{n}Y_{m,n}^{i}\left(t\right)\label{eq:membrane-potential-perturbative-expansion}
\end{equation}
}}
\par\end{center}{\small \par}
\end{onehalfspace}

\noindent which will be used to obtain an approximate analytic solution
of the system \ref{eq:rate-model-exact-equations-2}.

\subsection{\noindent \label{sub:The system of equations}The system of equations}

\noindent Now we put the perturbative expansion \ref{eq:membrane-potential-perturbative-expansion}
and the expressions \ref{eq:synaptic-weights-2} and \ref{eq:external-input}
for, respectively, the synaptic weights and the external input current,
inside the system \ref{eq:rate-model-exact-equations-2}. If all the
parameters $\sigma_{m}$ are small enough, we can expand the sigmoid
function in a Taylor series around $\mu$ (see \ref{eq:initial-conditions-mean}).
In order to be rigorous, we have to determine the radius of convergence
of the Taylor expansion of $S\left(V\right)$ for every value of $V$
and to check if it is big enough compared to $\sigma_{m}$, because
otherwise our technique cannot be applied. In fact, the various $\sigma_{m}$
determine the order of magnitude of the fluctuations of $V$ around
$\mu$, therefore it is important to check that $V$ is inside the
interval of convergence of the Taylor expansion of $S\left(V\right)$.
To our knowledge, this calculation has been performed only for $V=0$,
so in Appendix \ref{sec:Radius of convergence of the sigmoid and arctangent functions}
we show the general analysis, obtaining that in general the radius
of convergence decreases with the slope parameter $\lambda$ of the
sigmoid function. So, supposing that $\lambda$ is small enough, if
we call:

\begin{onehalfspace}
\begin{center}
\textit{\small{
\[
\zeta_{j}=\sum_{m=1}^{5}\sigma_{m}Y_{m}^{j}\left(t\right)+{\displaystyle \sum\limits _{\substack{m,n=1\\
m\leq n
}
}^{5}}\sigma_{m}\sigma_{n}Y_{m,n}^{j}\left(t\right)
\]
}}
\par\end{center}{\small \par}
\end{onehalfspace}

\noindent the Taylor expansion of the sigmoid function is:

\begin{onehalfspace}
\begin{center}
{\small{
\begin{align*}
S\left(\mu+\zeta_{j}\right)\approx & S\left(\mu\right)+S'\left(\mu\right)\zeta_{j}+\frac{1}{2}S''\left(\mu\right)\zeta_{j}^{2}\hphantom{\,\,\,\,\,\,\,\,\,\,\,\,\,\,\,\,\,\,\,\,\,\,\,\,\,\,\,\,\,\,\,\,\,\,\,\,\,\,\,\,\,\,\,\,\,\,\,\,\,\,\,\,\,\,\,\,\,\,\,\,\,\,\,\,\,\,\,\,\,\,\,\,\,\,\,\,\,\,\,\,\,\,\,\,\,\,\,\,\,\,\,\,\,}
\end{align*}
}}
\par\end{center}{\small \par}

\begin{center}
{\small{
\begin{align*}
\approx & S\left(\mu\right)+S'\left(\mu\right)\sum_{m=1}^{5}\sigma_{m}Y_{m}^{j}\left(t\right)\\
\\
 & +{\displaystyle \sum\limits _{\substack{m,n=1\\
m<n
}
}^{5}}\sigma_{m}\sigma_{n}\left[S'\left(\mu\right)Y_{m,n}^{j}\left(t\right)+S''\left(\mu\right)Y_{m}^{j}\left(t\right)Y_{n}^{j}\left(t\right)\right]\\
\\
 & +\sum_{m=1}^{5}\sigma_{m}^{2}\left[S'\left(\mu\right)Y_{m,m}^{j}\left(t\right)+\frac{1}{2}S''\left(\mu\right)\left(Y_{m}^{j}\left(t\right)\right)^{2}\right]
\end{align*}
}}
\par\end{center}{\small \par}
\end{onehalfspace}

\noindent having neglected the terms with order higher than $2$.
Now we substitute this expansion of the sigmoid function inside the
neural equation system and we equate the terms with the same $\sigma$
coefficients, obtaining (here we report only the equations that we
will actually use to compute the correlation structure in Section
\ref{sec:Correlation structure of the network}):

\begin{onehalfspace}
\begin{center}
{\small{
\begin{align}
\mu= & \tau\left[\Lambda S\left(\mu\right)+\overline{I}\right]\label{eq:perturbative-equation-1}\\
\nonumber \\
dY_{1}^{i}\left(t\right)= & \left[-\frac{1}{\tau}Y_{1}^{i}\left(t\right)+S'\left(\mu\right)\sum_{j=0}^{N-1}\overline{J}_{ij}Y_{1}^{j}\left(t\right)\right]dt+dB_{i}\left(t\right)\label{eq:perturbative-equation-2}\\
\nonumber \\
dY_{2}^{i}\left(t\right)= & \left[-\frac{1}{\tau}Y_{2}^{i}\left(t\right)+S'\left(\mu\right)\sum_{j=0}^{N-1}\overline{J}_{ij}Y_{2}^{j}\left(t\right)\right]dt\label{eq:perturbative-equation-3}\\
\nonumber \\
dY_{3}^{i}\left(t\right)= & \left[-\frac{1}{\tau}Y_{3}^{i}\left(t\right)+S'\left(\mu\right)\sum_{j=0}^{N-1}\overline{J}_{ij}Y_{3}^{j}\left(t\right)+S\left(\mu\right)\sum_{j=0}^{N-1}W_{ij}\right]dt\label{eq:perturbative-equation-4}\\
\nonumber \\
dY_{4}^{i}\left(t\right)= & \left[-\frac{1}{\tau}Y_{4}^{i}\left(t\right)+S'\left(\mu\right)\sum_{j=0}^{N-1}\overline{J}_{ij}Y_{4}^{j}\left(t\right)+S\left(\mu\right)\sum_{j=0}^{N-1}Z_{ij}\left(t\right)\right]dt\label{eq:perturbative-equation-5}\\
\nonumber \\
dY_{5}^{i}\left(t\right)= & \left[-\frac{1}{\tau}Y_{5}^{i}\left(t\right)+S'\left(\mu\right)\sum_{j=0}^{N-1}\overline{J}_{ij}Y_{5}^{j}\left(t\right)+H_{i}\left(t\right)\right]dt\hphantom{\,\,\,\,\,\,\,\,\,\,\,\,\,\,\,\,\,\,\,\,\,\,\,\,\,\,\,\,\,\,\,\,\,\,\,\,\,\,\,\,\,\,\,\,\,\,\,\,\,\,\,\,\,\,\,\,\,\,\,\,\,\,\,\,\,\,\,\,\,\,\,\,\,\,\,\,\,\,\,\,\,\,\,\,\,\,\,\,\,\,\,\,\,\,\,\,\,\,\,\,\,\,\,\,\,\,\,\,\,\,\,\,\,\,\,\,\,\,\,\,\,\,\,\,\,\,\,\,\,\,\,\,\,\,\,\,}\label{eq:perturbative-equation-6}
\end{align}
}}
\par\end{center}{\small \par}

\begin{center}
{\small{
\begin{align}
\vdots\nonumber \\
dY_{1,4}^{i}\left(t\right)= & \left[-\frac{1}{\tau}Y_{1,4}^{i}\left(t\right)+S'\left(\mu\right)\sum_{j=0}^{N-1}\overline{J}_{ij}Y_{1,4}^{j}\left(t\right)+S'\left(\mu\right)\sum_{j=0}^{N-1}Z_{ij}\left(t\right)Y_{1}^{j}\left(t\right)+S''\left(\mu\right)\sum_{j=0}^{N-1}\overline{J}_{ij}Y_{1}^{j}\left(t\right)Y_{4}^{j}\left(t\right)\right]dt\label{eq:perturbative-equation-7}\\
\nonumber \\
dY_{1,5}^{i}\left(t\right)= & \left[-\frac{1}{\tau}Y_{1,5}^{i}\left(t\right)+S'\left(\mu\right)\sum_{j=0}^{N-1}\overline{J}_{ij}Y_{1,5}^{j}\left(t\right)+S''\left(\mu\right)\sum_{j=0}^{N-1}\overline{J}_{ij}Y_{1}^{j}\left(t\right)Y_{5}^{j}\left(t\right)\right]dt\label{eq:perturbative-equation-8}\\
\vdots\nonumber \\
dY_{2,4}^{i}\left(t\right)= & \left[-\frac{1}{\tau}Y_{2,4}^{i}\left(t\right)+S'\left(\mu\right)\sum_{j=0}^{N-1}\overline{J}_{ij}Y_{2,4}^{j}\left(t\right)+S'\left(\mu\right)\sum_{j=0}^{N-1}Z_{ij}\left(t\right)Y_{2}^{j}\left(t\right)+S''\left(\mu\right)\sum_{j=0}^{N-1}\overline{J}_{ij}Y_{2}^{j}\left(t\right)Y_{4}^{j}\left(t\right)\right]dt\label{eq:perturbative-equation-9}\\
\nonumber \\
dY_{2,5}^{i}\left(t\right)= & \left[-\frac{1}{\tau}Y_{2,5}^{i}\left(t\right)+S'\left(\mu\right)\sum_{j=0}^{N-1}\overline{J}_{ij}Y_{2,5}^{j}\left(t\right)+S''\left(\mu\right)\sum_{j=0}^{N-1}\overline{J}_{ij}Y_{2}^{j}\left(t\right)Y_{5}^{j}\left(t\right)\right]dt\label{eq:perturbative-equation-10}\\
\vdots\nonumber \\
dY_{3,4}^{i}\left(t\right)= & \left[-\frac{1}{\tau}Y_{3,4}^{i}\left(t\right)+S'\left(\mu\right)\sum_{j=0}^{N-1}\overline{J}_{ij}Y_{3,4}^{j}\left(t\right)+S'\left(\mu\right)\sum_{j=0}^{N-1}Z_{ij}\left(t\right)Y_{3}^{j}\left(t\right)\right.\nonumber \\
 & \left.+S'\left(\mu\right)\sum_{j=0}^{N-1}W_{ij}Y_{4}^{j}\left(t\right)+S''\left(\mu\right)\sum_{j=0}^{N-1}\overline{J}_{ij}Y_{3}^{j}\left(t\right)Y_{4}^{j}\left(t\right)\right]dt\label{eq:perturbative-equation-11}\\
\nonumber \\
dY_{3,5}^{i}\left(t\right)= & \left[-\frac{1}{\tau}Y_{3,5}^{i}\left(t\right)+S'\left(\mu\right)\sum_{j=0}^{N-1}\overline{J}_{ij}Y_{3,5}^{j}\left(t\right)+S'\left(\mu\right)\sum_{j=0}^{N-1}W_{ij}Y_{5}^{j}\left(t\right)+S''\left(\mu\right)\sum_{j=0}^{N-1}\overline{J}_{ij}Y_{3}^{j}\left(t\right)Y_{5}^{j}\left(t\right)\right]dt\label{eq:perturbative-equation-12}\\
\vdots\nonumber 
\end{align}
}}
\par\end{center}{\small \par}
\end{onehalfspace}

\noindent Equation \ref{eq:perturbative-equation-1} is algebraic
and non-linear, therefore must be solved numerically. \ref{eq:perturbative-equation-2}
is the only stochastic differential equation of the set and can be
solved analytically, since it is linear with constant coefficients.
Equations \ref{eq:perturbative-equation-3} - \ref{eq:perturbative-equation-6}
are ordinary, and can be solved in the same way as \ref{eq:perturbative-equation-2}.
To conclude, equations \ref{eq:perturbative-equation-7} - \ref{eq:perturbative-equation-12}
determine the functions $Y_{m,n}^{i}\left(t\right)$, and depend on
the terms $Y_{m}^{i}\left(t\right)$, which have been calculated at
the previous step. Being linear and with constant coefficients, they
can be integrated analytically as a function of the already known
functions $Y_{m}^{i}\left(t\right)$.

\subsection{\noindent \label{sub:The initial conditions}The initial conditions}

\noindent The perturbative expansion \ref{eq:membrane-potential-perturbative-expansion}
at $t=0$ gives:

\begin{onehalfspace}
\begin{center}
{\small{
\[
V_{i}\left(0\right)\approx\mu+\sum_{m=1}^{5}\sigma_{m}Y_{m}^{i}\left(0\right)+{\displaystyle \sum\limits _{\substack{m,n=1\\
m\leq n
}
}^{5}}\sigma_{m}\sigma_{n}Y_{m,n}^{i}\left(0\right)
\]
}}
\par\end{center}{\small \par}
\end{onehalfspace}

\noindent Moreover, according to \ref{eq:initial-conditions}, we
have $V_{i}\left(0\right)\sim\mathcal{N}\left(\mu,\sigma_{2}^{2}\right)=\mu+\sigma_{2}\mathcal{N}\left(0,1\right)$,
so from the comparison it must be:

\begin{onehalfspace}
\begin{center}
\begin{align}
 & Y_{2}^{i}\left(0\right)\sim\mathcal{N}\left(0,1\right)\label{eq:initial-conditions-decomposed-1}\\
\nonumber \\
 & \begin{array}{ccc}
Y_{m}^{i}\left(0\right)=0, &  & m=1,3,4,5\end{array}\label{eq:initial-conditions-decomposed-2}\\
\nonumber \\
 & \begin{array}{ccc}
Y_{m,n}^{i}\left(0\right)=0, &  & \forall\left(m,n\right):\, m\leq n\end{array}\label{eq:initial-conditions-decomposed-3}
\end{align}

\par\end{center}
\end{onehalfspace}

\noindent So we have $V_{i}\left(0\right)=\mu+\sigma_{2}Y_{2}^{i}\left(0\right)$
and therefore:

\begin{onehalfspace}
\begin{center}
{\small{
\[
Cov\left(V_{i}\left(0\right),V_{j}\left(0\right)\right)=\sigma_{2}^{2}Cov\left(Y_{2}^{i}\left(0\right),Y_{2}^{j}\left(0\right)\right)
\]
}}
\par\end{center}{\small \par}
\end{onehalfspace}

\noindent But from \ref{eq:initial-conditions-covariance} we also
know that:

\begin{onehalfspace}
\begin{center}
{\small{
\[
Cov\left(V_{i}\left(0\right),V_{j}\left(0\right)\right)=\begin{cases}
\sigma_{2}^{2} & \begin{array}{ccc}
\mathrm{if} &  & i=j\end{array}\\
\\
\sigma_{2}^{2}C_{2} & \begin{array}{ccc}
\mathrm{if} &  & i\neq j\end{array}
\end{cases}
\]
}}
\par\end{center}{\small \par}
\end{onehalfspace}

\noindent so from the comparison it must be that:

\begin{onehalfspace}
\begin{center}
{\small{
\begin{equation}
Cov\left(Y_{2}^{i}\left(0\right),Y_{2}^{j}\left(0\right)\right)=\begin{cases}
1 & \begin{array}{ccc}
\mathrm{if} &  & i=j\end{array}\\
\\
C_{2} & \begin{array}{ccc}
\mathrm{if} &  & i\neq j\end{array}
\end{cases}\label{eq:initial-conditions-covariance-decomposed}
\end{equation}
}}
\par\end{center}{\small \par}
\end{onehalfspace}

\subsection{\noindent \label{sub:Solutions of the equations}Solutions of the
equations}

\noindent As we said at the end of Section \ref{sub:The system of equations},
the algebraic equation \ref{eq:perturbative-equation-1} is non-linear,
therefore it cannot be solved exactly. However, the differential equations
satisfied by all the functions $Y_{m}^{i}\left(t\right)$ and $Y_{m,n}^{i}\left(t\right)$
are linear with constant coefficients, therefore they can be solved
analytically. In particular, the equations \ref{eq:perturbative-equation-2}
- \ref{eq:perturbative-equation-6} can be solved directly. Instead
the remaining equations are functions of the previous $Y_{m}^{i}\left(t\right)$,
that we have already calculated. For example, according to \ref{eq:perturbative-equation-7},
$Y_{1,4}^{i}\left(t\right)$ can be determined analytically as a function
of $Y_{1}^{i}\left(t\right)$ and $Y_{4}^{i}\left(t\right)$, which
are already known from the equations \ref{eq:perturbative-equation-2}
and \ref{eq:perturbative-equation-5}. Now we introduce the \textit{fundamental
matrix} $\Phi\left(t\right)$ such that:

\begin{onehalfspace}
\begin{center}
{\small{
\begin{align}
\Phi\left(t\right)= & e^{At}\nonumber \\
\label{eq:Phi-and-A-matrices}\\
A_{ij}= & \begin{cases}
-\frac{1}{\tau} & \begin{array}{ccc}
\mathrm{if} &  & i=j\end{array}\\
\\
\overline{J}_{ij}S'\left(\mu\right) & \begin{array}{ccc}
\mathrm{if} &  & i\neq j\end{array}
\end{cases}\nonumber 
\end{align}
}}
\par\end{center}{\small \par}
\end{onehalfspace}

\noindent where:

\begin{onehalfspace}
\begin{center}
\textit{\small{
\begin{equation}
\mathcal{J}=\overline{J}S'\left(\mu\right)\label{eq:effective-connectivity-matrix}
\end{equation}
}}
\par\end{center}{\small \par}
\end{onehalfspace}

\noindent is the \textit{effective connectivity} matrix of the network.
Therefore the solutions of all the functions $Y_{m}^{i}\left(t\right)$
can be obtained straightforwardly as follows:

{\small{
\begin{align}
 & Y_{1}^{i}\left(t\right)=\sum_{j=0}^{N-1}\int_{0}^{t}\left[\Phi\left(t-s\right)\right]_{ij}dB_{j}\left(s\right)\hphantom{\,\,\,\,\,\,\,\,\,\,\,\,\,\,\,\,\,\,\,\,\,\,\,\,\,\,\,\,\,\,\,\,\,\,\,\,\,\,\,\,\,\,\,\,\,\,\,\,\,\,\,\,\,\,\,\,\,\,\,\,\,\,\,\,\,\,\,\,\,\,\,\,\,\,\,\,\,\,\,\,\,\,\,\,\,\,\,\,\,\,\,\,\,\,\,\,\,\,\,\,\,\,\,\,\,\,\,\,\,\,\,\,\,\,\,\,\,\,\,\,\,\,\,\,\,\,\,\,\,\,\,\,\,\,\,\,\,\,\,\,\,\,\,\,\,\,\,\,\,\,\,\,\,\,\,\,\,\,\,\,\,\,\,\,\,\,\,\,\,\,\,\,\,\,\,\,\,\,\,\,\,\,\,\,\,\,\,\,\,\,\,\,\,\,\,\,\,}\label{eq:solution-perturbative-equation-1}
\end{align}
}}{\small \par}

{\footnotesize{
\begin{align}
Y_{2}^{i}\left(t\right)= & \sum_{j=0}^{N-1}\Phi_{ij}\left(t\right)Y_{2}^{j}\left(0\right)\label{eq:solution-perturbative-equation-2}\\
\nonumber \\
Y_{3}^{i}\left(t\right)= & S\left(\mu\right)\sum_{j,k=0}^{N-1}W_{jk}\int_{0}^{t}\left[\Phi\left(t-s\right)\right]_{ij}ds\label{eq:solution-perturbative-equation-3}\\
\nonumber \\
Y_{4}^{i}\left(t\right)= & S\left(\mu\right)\sum_{j,k=0}^{N-1}\int_{0}^{t}\left[\Phi\left(t-s\right)\right]_{ij}Z_{jk}\left(s\right)ds\label{eq:solution-perturbative-equation-4}\\
\nonumber \\
Y_{5}^{i}\left(t\right)= & \sum_{j=0}^{N-1}\int_{0}^{t}\left[\Phi\left(t-s\right)\right]_{ij}H_{j}\left(s\right)ds\label{eq:solution-perturbative-equation-5}\\
\vdots\nonumber \\
Y_{1,4}^{i}\left(t\right)= & S'\left(\mu\right)\sum_{j,k,l=0}^{N-1}\int_{0}^{t}\left[\Phi\left(t-s\right)\right]_{ij}\left\{ \int_{0}^{s}\left[\Phi\left(s-u\right)\right]_{kl}dB_{l}\left(u\right)\right\} Z_{jk}\left(s\right)ds\nonumber \\
 & +S\left(\mu\right)S''\left(\mu\right)\sum_{j,k,l,m,n=0}^{N-1}\overline{J}_{jk}\int_{0}^{t}\left[\Phi\left(t-s\right)\right]_{ij}\left\{ \int_{0}^{s}\left[\Phi\left(s-u\right)\right]_{kl}dB_{l}\left(u\right)\right\} \left\{ \int_{0}^{s}\left[\Phi\left(s-u\right)\right]_{km}Z_{mn}\left(u\right)du\right\} ds\label{eq:solution-perturbative-equation-6}\\
\nonumber \\
Y_{1,5}^{i}\left(t\right)= & S''\left(\mu\right)\sum_{j,k,l,m=0}^{N-1}\overline{J}_{jk}\int_{0}^{t}\left[\Phi\left(t-s\right)\right]_{ij}\left\{ \int_{0}^{s}\left[\Phi\left(s-u\right)\right]_{kl}dB_{l}\left(u\right)\right\} \left\{ \int_{0}^{s}\left[\Phi\left(s-u\right)\right]_{km}H_{m}\left(u\right)du\right\} ds\label{eq:solution-perturbative-equation-7}\\
\vdots\nonumber \\
Y_{2,4}^{i}\left(t\right)= & S'\left(\mu\right)\sum_{j,k,l=0}^{N-1}Y_{2}^{l}\left(0\right)\int_{0}^{t}\left[\Phi\left(t-s\right)\right]_{ij}\Phi_{kl}\left(s\right)Z_{jk}\left(s\right)ds\nonumber \\
 & +S\left(\mu\right)S''\left(\mu\right)\sum_{j,k,l,m,n=0}^{N-1}\overline{J}_{jk}Y_{2}^{l}\left(0\right)\int_{0}^{t}\left[\Phi\left(t-s\right)\right]_{ij}\Phi_{kl}\left(s\right)\left\{ \int_{0}^{s}\left[\Phi\left(s-u\right)\right]_{km}Z_{mn}\left(u\right)du\right\} ds\label{eq:solution-perturbative-equation-8}\\
\nonumber \\
Y_{2,5}^{i}\left(t\right)= & S''\left(\mu\right)\sum_{j,k,l,m=0}^{N-1}\overline{J}_{jk}Y_{2}^{l}\left(0\right)\int_{0}^{t}\left[\Phi\left(t-s\right)\right]_{ij}\Phi_{kl}\left(s\right)\left\{ \int_{0}^{s}\left[\Phi\left(s-u\right)\right]_{km}H_{m}\left(u\right)du\right\} ds\label{eq:solution-perturbative-equation-9}\\
\vdots\nonumber \\
Y_{3,4}^{i}\left(t\right)= & S\left(\mu\right)S'\left(\mu\right)\sum_{j,k,l,m=0}^{N-1}W_{lm}\int_{0}^{t}\left[\Phi\left(t-s\right)\right]_{ij}\left\{ \int_{0}^{s}\left[\Phi\left(s-u\right)\right]_{kl}du\right\} Z_{jk}\left(s\right)ds\nonumber \\
 & +S\left(\mu\right)S'\left(\mu\right)\sum_{j,k,l,m=0}^{N-1}W_{jk}\int_{0}^{t}\left[\Phi\left(t-s\right)\right]_{ij}\left\{ \int_{0}^{s}\left[\Phi\left(s-u\right)\right]_{kl}Z_{lm}\left(u\right)du\right\} ds\nonumber \\
 & +S^{2}\left(\mu\right)S''\left(\mu\right)\sum_{j,k,l,m,n,p=0}^{N-1}\overline{J}_{jk}W_{lm}\int_{0}^{t}\left[\Phi\left(t-s\right)\right]_{ij}\left\{ \int_{0}^{s}\left[\Phi\left(s-u\right)\right]_{kl}du\right\} \left\{ \int_{0}^{s}\left[\Phi\left(s-u\right)\right]_{kn}Z_{np}\left(u\right)du\right\} ds\label{eq:solution-perturbative-equation-10}
\end{align}
}}{\footnotesize \par}

{\footnotesize{
\begin{align}
Y_{3,5}^{i}\left(t\right)= & S'\left(\mu\right)\sum_{j,k,l=0}^{N-1}W_{jk}\int_{0}^{t}\left[\Phi\left(t-s\right)\right]_{ij}\left\{ \int_{0}^{s}\left[\Phi\left(s-u\right)\right]_{kl}H_{l}\left(u\right)du\right\} ds\nonumber \\
 & +S\left(\mu\right)S''\left(\mu\right)\sum_{j,k,l,m,n=0}^{N-1}\overline{J}_{jk}W_{lm}\int_{0}^{t}\left[\Phi\left(t-s\right)\right]_{ij}\left\{ \int_{0}^{s}\left[\Phi\left(s-u\right)\right]_{kl}du\right\} \left\{ \int_{0}^{s}\left[\Phi\left(s-u\right)\right]_{kn}H_{n}\left(u\right)du\right\} ds\label{eq:solution-perturbative-equation-11}\\
\vdots\nonumber 
\end{align}
}}{\footnotesize \par}

\noindent To conclude, we have performed a perturbative expansion
around a stationary state $\mu$ because in this way the equations
\ref{eq:perturbative-equation-2} - \ref{eq:perturbative-equation-12}
have constant coefficients and therefore they can be solved exactly
using the fundamental matrix \ref{eq:Phi-and-A-matrices}. Had we
performed the perturbative expansion around a non-stationary state,
we would have obtained a system of differential equations with time-varying
coefficients, whose general solution is not known. In this case, the
best thing that we can try is to write the solution in terms of the
\textit{Magnus expansion} \cite{CPA:CPA3160070404}, but this introduces
another approximation to the real solution of the neural network.

\bigskip{}

\noindent In this article we have also supposed that the system is
invariant under exchange of the neural indices: for this reason we
have used the same stationary solution $\mu$, the same (unperturbed)
input current $\overline{I}$ and the same number of incoming connections
for all the neurons in the network. This invariance is required in
order to ensure that the effective connectivity matrix $\mathcal{J}$
given by \ref{eq:effective-connectivity-matrix} has the same structure
as the real and unperturbed connectivity matrix $\overline{J}$. In
this way the fundamental matrix $\Phi\left(t\right)$ can be calculated
using the properties of $\overline{J}$, as explained in Section \ref{sec:Calculation of the fundamental matrix}.
If the system is not invariant under exchange of the neural indices,
$\mathcal{J}$ does not inherit the structure of $\overline{J}$,
therefore the technique introduced in this article cannot be used
anymore (see also the discussion at the end of Section \ref{sub:Symmetric matrices}).
To conclude, it is important to observe that even if we have chosen
structures of $\Sigma_{1}$, $\Sigma_{2}$, $\Sigma_{3}$ and $\Omega_{3}$
that are invariant under exchange of the neural indices, their invariance
is not required here: we have used it only to simplify the final formulae
that we will obtain in Section \ref{sec:Correlation structure of the network}.
Therefore in principle inhomogeneous structures can be used for these
covariance matrices.

\section{\noindent \label{sec:Correlation structure of the network}Correlation
structure of the network}

\noindent In this section we want to calculate the correlation structure
of the membrane potentials, according to the perturbative expansion
\ref{eq:membrane-potential-perturbative-expansion}. Since the covariance
function is bilinear, we have to compute it for all the possible combinations
of the pairs $\left(Y_{m}^{i}\left(t\right),Y_{n}^{j}\left(t\right)\right)$,
$\left(Y_{m}^{i}\left(t\right),Y_{n,p}^{j}\left(t\right)\right)$
and $\left(Y_{m,n}^{i}\left(t\right),Y_{p,q}^{j}\left(t\right)\right)$.
However we do not have to consider the terms of order $4$, like $\sigma_{1}^{2}\sigma_{2}^{2}Cov\left(Y_{1,1}^{i}\left(t\right),Y_{2,2}^{j}\left(t\right)\right)$,
because they are incomplete. In effect, in the perturbative expansion
of $V_{i}\left(t\right)$, we did not consider the terms of order
$3$, like $\sigma_{1}^{2}\sigma_{2}Y_{1,1,2}^{i}\left(t\right)$,
that generate contributions of order $4$ in the formula of the covariance.
So the terms of order $4$ cannot be considered in the expansion of
the covariance, therefore the final formula is of order $3$.

\bigskip{}

\noindent For simplicity we suppose that the Brownian motions, the
initial conditions and the uncertainty of the synaptic weights are
$3$ \textit{independent} random processes (and indeed there is a
priori no obvious reason to think that they are correlated), so all
the cross terms like $\sigma_{1}\sigma_{2}Cov\left(Y_{1}^{i}\left(t\right),Y_{2}^{j}\left(t\right)\right)$,
$\sigma_{1}\sigma_{2}Cov\left(Y_{2}^{i}\left(t\right),Y_{1}^{j}\left(t\right)\right)$,
$\sigma_{1}^{2}\sigma_{3}Cov\left(Y_{1,1}^{i}\left(t\right),Y_{3}^{j}\left(t\right)\right)$,
... are equal to zero (however, if desired, we could assume non-zero
correlations between these $3$ sources of randomness, since there
is no technical difficulty in the calculations, only the problem to
compute many non-zero cross terms). Let us show it with an example:

\begin{onehalfspace}
\begin{center}
{\small{
\begin{align*}
Cov\left(Y_{1}^{i}\left(t\right),Y_{2}^{j}\left(t\right)\right)= & Cov\left(\int_{0}^{t}\sum_{k=0}^{N-1}\left[\Phi\left(t-s\right)\right]_{ik}dB_{k}\left(s\right),\sum_{l=0}^{N-1}\Phi_{jl}\left(t\right)Y_{2}^{l}\left(0\right)\right)\\
\\
= & \sum_{k,l=0}^{N-1}\Phi_{jl}\left(t\right)Cov\left(\int_{0}^{t}\left[\Phi\left(t-s\right)\right]_{ik}dB_{k}\left(s\right),Y_{2}^{l}\left(0\right)\right)\\
\\
= & 0
\end{align*}
}}
\par\end{center}{\small \par}
\end{onehalfspace}

\noindent since $B_{k}\left(s\right)$ and $Y_{2}^{l}\left(0\right)$
are independent by assumption. Moreover, due to the Isserlis' theorem
\cite{1918}, we obtain also that all the terms in the covariance
proportional to $\sigma_{m}^{2}\sigma_{n}$ with $m,n=1,2,3$ are
equal to zero, like $\sigma_{1}^{2}\sigma_{2}$ and $\sigma_{3}^{3}$.
The same thing happens to all the terms proportional to $\sigma_{m}^{2}\sigma_{n}$,
with $m=4,5$ and $n=1,2,3$. This is due to the fact that, according
to the Isserlis' theorem again, the mean of the product of any odd
number of zero-mean normal processes is equal to zero. We show it
with an example:

\begin{onehalfspace}
\begin{center}
{\small{
\begin{align*}
 & Cov\left(Y_{2}^{i}\left(t\right),Y_{2,2}^{j}\left(t\right)\right)\\
\\
 & =Cov\left(\sum_{k=0}^{N-1}\Phi_{ik}\left(t\right)Y_{2}^{k}\left(0\right),\frac{1}{2}S''\left(\mu\right)\sum_{l,m,n,p=0}^{N-1}\overline{J}_{lm}Y_{2}^{n}\left(0\right)Y_{2}^{p}\left(0\right)\int_{0}^{t}\left[\Phi\left(t-s\right)\right]_{jl}\Phi_{mn}\left(s\right)\Phi_{mp}\left(s\right)ds\right)\\
\\
 & =\frac{1}{2}S''\left(\mu\right)\sum_{k,l,m,n,p=0}^{N-1}\Phi_{ik}\left(t\right)\overline{J}_{lm}\left\{ \int_{0}^{t}\left[\Phi\left(t-s\right)\right]_{jl}\Phi_{mn}\left(s\right)\Phi_{mp}\left(s\right)ds\right\} Cov\left(Y_{2}^{k}\left(0\right),Y_{2}^{n}\left(0\right)Y_{2}^{p}\left(0\right)\right)\\
\\
 & =0
\end{align*}
}}
\par\end{center}{\small \par}
\end{onehalfspace}

\noindent because:

\begin{onehalfspace}
\begin{center}
{\small{
\[
Cov\left(Y_{2}^{k}\left(0\right),Y_{2}^{n}\left(0\right)Y_{2}^{p}\left(0\right)\right)=\mathbb{E}\left[Y_{2}^{k}\left(0\right)Y_{2}^{n}\left(0\right)Y_{2}^{p}\left(0\right)\right]-\mathbb{E}\left[Y_{2}^{k}\left(0\right)\right]\mathbb{E}\left[Y_{2}^{n}\left(0\right)Y_{2}^{p}\left(0\right)\right]=0
\]
}}
\par\end{center}{\small \par}
\end{onehalfspace}

\noindent since $\mathbb{E}\left[Y_{2}^{k}\left(0\right)Y_{2}^{n}\left(0\right)Y_{2}^{p}\left(0\right)\right]=0$
by the Isserlis' theorem and $\mathbb{E}\left[Y_{2}^{k}\left(0\right)\right]=0$,
because $Y_{2}^{k}\left(0\right)\sim\mathcal{N}\left(0,1\right)$
from \ref{eq:initial-conditions-decomposed-1}.

\noindent \bigskip{}

\noindent In the final formula of the covariance, also the terms proportional
to $\sigma_{m}\sigma_{n}$ and $\sigma_{m}^{2}\sigma_{n}$ with $m,n=4,5$
are zero, because the functions $Y_{m}^{i}\left(t\right)$ and $Y_{m,n}^{i}\left(t\right)$
are deterministic for $m,n=4,5$. In fact, for example, from the formulae
\ref{eq:solution-perturbative-equation-4} and \ref{eq:solution-perturbative-equation-5}
we can easily see that the functions $Y_{4}^{i}\left(t\right)$ and
$Y_{5}^{i}\left(t\right)$ depend only on deterministic functions
($\Phi\left(t\right)$, $Z_{jk}\left(t\right)$ and $H_{j}\left(t\right)$),
deterministic parameters ($\tau$ and all the parameters of $S\left(\cdot\right)$)
and deterministic initial conditions ($Y_{4}^{i}\left(0\right)=Y_{5}^{i}\left(0\right)=0$,
from \ref{eq:initial-conditions-decomposed-2}), and therefore they
are deterministic as well. Also the terms proportional to $\sigma_{m}\sigma_{n}\sigma_{p}$
for $m=4,5$ and $n\neq p$ are zero, due to the independence of the
sources of randomness or to the fact that $Y_{m}^{i}\left(t\right)$
is deterministic for $m=4,5$. In the same way the terms obtained
from the covariance of $Y_{m}^{i}\left(t\right)$ for $m=4,5$ with
$Y_{n,n}^{i}\left(t\right)$ for $n=1,2,3$ are zero due to the fact
that the first function is deterministic.

\noindent To conclude, the only non-zero terms in the final formula
of the covariance are those proportional to $\sigma_{m}^{2}$ for
$m=1,2,3$, and those obtained from the covariance of $Y_{m,n}^{i}\left(t\right)$
with $Y_{m}^{i}\left(t\right)$, for $m=1,2,3$ and $n=4,5$. So the
final formula for the covariance is:

\begin{onehalfspace}
\begin{center}
{\footnotesize{
\begin{align}
 & Cov\left(V_{i}\left(t\right),V_{j}\left(t\right)\right)\nonumber \\
\nonumber \\
 & =\sigma_{1}^{2}Cov\left(Y_{1}^{i}\left(t\right),Y_{1}^{j}\left(t\right)\right)+\sigma_{2}^{2}Cov\left(Y_{2}^{i}\left(t\right),Y_{2}^{j}\left(t\right)\right)+\sigma_{3}^{2}Cov\left(Y_{3}^{i}\left(t\right),Y_{3}^{j}\left(t\right)\right)\nonumber \\
\nonumber \\
 & +\sigma_{4}\left\{ \sigma_{1}^{2}\left[Cov\left(Y_{1}^{i}\left(t\right),Y_{1,4}^{j}\left(t\right)\right)+Cov\left(Y_{1,4}^{i}\left(t\right),Y_{1}^{j}\left(t\right)\right)\right]+\sigma_{2}^{2}\left[Cov\left(Y_{2}^{i}\left(t\right),Y_{2,4}^{j}\left(t\right)\right)+Cov\left(Y_{2,4}^{i}\left(t\right),Y_{2}^{j}\left(t\right)\right)\right]\right.\nonumber \\
\nonumber \\
 & \left.+\sigma_{3}^{2}\left[Cov\left(Y_{3}^{i}\left(t\right),Y_{3,4}^{j}\left(t\right)\right)+Cov\left(Y_{3,4}^{i}\left(t\right),Y_{3}^{j}\left(t\right)\right)\right]\right\} \nonumber \\
\nonumber \\
 & +\sigma_{5}\left\{ \sigma_{1}^{2}\left[Cov\left(Y_{1}^{i}\left(t\right),Y_{1,5}^{j}\left(t\right)\right)+Cov\left(Y_{1,5}^{i}\left(t\right),Y_{1}^{j}\left(t\right)\right)\right]+\sigma_{2}^{2}\left[Cov\left(Y_{2}^{i}\left(t\right),Y_{2,5}^{j}\left(t\right)\right)+Cov\left(Y_{2,5}^{i}\left(t\right),Y_{2}^{j}\left(t\right)\right)\right]\right.\nonumber \\
\nonumber \\
 & \left.+\sigma_{3}^{2}\left[Cov\left(Y_{3}^{i}\left(t\right),Y_{3,5}^{j}\left(t\right)\right)+Cov\left(Y_{3,5}^{i}\left(t\right),Y_{3}^{j}\left(t\right)\right)\right]\right\} \label{eq:covariance}
\end{align}
}}
\par\end{center}{\footnotesize \par}
\end{onehalfspace}

\noindent Even if the third order terms can be calculated exactly
using the Isserlis' theorem (and even if in principle we can extend
this perturbative expansion to any higher order), due to their complexity
in this article we consider only the second order terms, that is equivalent
to say that we truncate the perturbative expansion \ref{eq:membrane-potential-perturbative-expansion}
of the membrane potential at the first order. After some algebra we
obtain:

\begin{onehalfspace}
\begin{center}
{\small{
\begin{align}
Cov\left(Y_{1}^{i}\left(t\right),Y_{1}^{j}\left(t\right)\right)= & \sum_{k=0}^{N-1}\int_{0}^{t}\left[\Phi\left(t-s\right)\right]_{ik}\left[\Phi\left(t-s\right)\right]_{jk}ds\nonumber \\
 & +C_{1}{\displaystyle \sum\limits _{\substack{k,l=0\\
k\neq l
}
}^{N-1}}\int_{0}^{t}\left[\Phi\left(t-s\right)\right]_{ik}\left[\Phi\left(t-s\right)\right]_{jl}ds\label{eq:covariance-part-1}\\
\nonumber \\
Cov\left(Y_{2}^{i}\left(t\right),Y_{2}^{j}\left(t\right)\right)= & \sum_{k=0}^{N-1}\Phi_{ik}\left(t\right)\Phi_{jk}\left(t\right)+C_{2}{\displaystyle \sum\limits _{\substack{k,l=0\\
k\neq l
}
}^{N-1}}\Phi_{ik}\left(t\right)\Phi_{jl}\left(t\right)\label{eq:covariance-part-2}\\
\nonumber \\
Cov\left(Y_{3}^{i}\left(t\right),Y_{3}^{j}\left(t\right)\right)= & \frac{S^{2}\left(\mu\right)}{M}\sum_{k=0}^{N-1}\left\{ \int_{0}^{t}\left[\Phi\left(t-s\right)\right]_{ik}ds\right\} \left\{ \int_{0}^{t}\left[\Phi\left(t-s\right)\right]_{jk}ds\right\} \nonumber \\
 & +C_{3}S^{2}\left(\mu\right)\left\{ \sum_{k,l=0}^{N-1}\left\{ \int_{0}^{t}\left[\Phi\left(t-s\right)\right]_{ik}ds\right\} \left\{ \int_{0}^{t}\left[\Phi\left(t-s\right)\right]_{jl}ds\right\} \right.\nonumber \\
 & \left.-\frac{1}{M}\sum_{k=0}^{N-1}\left\{ \int_{0}^{t}\left[\Phi\left(t-s\right)\right]_{ik}ds\right\} \left\{ \int_{0}^{t}\left[\Phi\left(t-s\right)\right]_{jk}ds\right\} \vphantom{\sum_{k,l=0}^{N-1}}\right\} \label{eq:covariance-part-3}
\end{align}
}}
\par\end{center}{\small \par}
\end{onehalfspace}

\noindent So now the covariance $Cov\left(V_{i}\left(t\right),V_{j}\left(t\right)\right)$
is known for all the possible pairs $\left(i,j\right)$, with $i,j=0,1,...,N-1$,
therefore we can determine the correlation structure of the network
using the formula for the Pearson's correlation coefficient:

\begin{onehalfspace}
\begin{center}
{\small{
\begin{equation}
Corr_{2}\left(V_{i}\left(t\right),V_{j}\left(t\right)\right)=\frac{Cov\left(V_{i}\left(t\right),V_{j}\left(t\right)\right)}{\sqrt{Var\left(V_{i}\left(t\right)\right)Var\left(V_{j}\left(t\right)\right)}}\label{eq:correlation}
\end{equation}
}}
\par\end{center}{\small \par}
\end{onehalfspace}

\noindent where:

\begin{onehalfspace}
\begin{center}
\begin{equation}
Var\left(V_{i}\left(t\right)\right)=Cov\left(V_{i}\left(t\right),V_{i}\left(t\right)\right)\label{eq:variance}
\end{equation}

\par\end{center}
\end{onehalfspace}

\noindent is the variance of the stochastic process $V_{i}\left(t\right)$.
The subscript ``$2$'' means that this is a correlation between
a pair of neurons.

\noindent In order to determine the higher order correlations between
triplets, quadruplets, quintuplets etc of neurons, we have to extend
the Pearson's formula in the following way. The natural generalization
of the covariance for $n$ functions is:

\begin{onehalfspace}
\begin{center}
{\small{
\begin{equation}
\kappa_{n}\left(V_{i_{0}}\left(t\right),V_{i_{1}}\left(t\right),...,V_{i_{n-1}}\left(t\right)\right)=\mathbb{E}\left[\prod_{j=0}^{n-1}\left(V_{i_{j}}\left(t\right)-\overline{V}_{i_{j}}\left(t\right)\right)\right]\label{eq:joint-cumulant}
\end{equation}
}}
\par\end{center}{\small \par}
\end{onehalfspace}

\noindent This is known as the \textit{joint cumulant} of the functions
$V_{i_{0}}\left(t\right),V_{i_{1}}\left(t\right),...,V_{i_{n-1}}\left(t\right)$.
Unfortunately this is not enough, because as with the Pearson's correlation
coefficient, we want to normalize the joint cumulant in order to find
a function that is in the range $\left[-1,1\right]$. To this purpose,
we can observe that:

\begin{onehalfspace}
\begin{center}
{\small{
\[
\left|\mathbb{E}\left[\prod_{j=0}^{n-1}\left(V_{i_{j}}\left(t\right)-\overline{V}_{i_{j}}\left(t\right)\right)\right]\right|\leq\mathbb{E}\left[\left|\prod_{j=0}^{n-1}\left(V_{i_{j}}\left(t\right)-\overline{V}_{i_{j}}\left(t\right)\right)\right|\right]\leq\left\{ \prod_{j=0}^{n-1}\mathbb{E}\left[\left|V_{i_{j}}\left(t\right)-\overline{V}_{i_{j}}\left(t\right)\right|^{n}\right]\right\} ^{\frac{1}{n}}
\]
}}
\par\end{center}{\small \par}
\end{onehalfspace}

\noindent having used the fact that $\left|a+b\right|\leq\left|a\right|+\left|b\right|$
at the first step and a special case of the H�lder's inequality at
the second. Therefore we have:

\begin{onehalfspace}
\begin{center}
{\small{
\begin{equation}
\left|\frac{\mathbb{E}\left[{\displaystyle \prod_{j=0}^{n-1}}\left(V_{i_{j}}\left(t\right)-\overline{V}_{i_{j}}\left(t\right)\right)\right]}{\sqrt[n]{{\displaystyle \prod_{j=0}^{n-1}}\mathbb{E}\left[\left|V_{i_{j}}\left(t\right)-\overline{V}_{i_{j}}\left(t\right)\right|^{n}\right]}}\right|\leq1\label{eq:normalization-fulfilled}
\end{equation}
}}
\par\end{center}{\small \par}
\end{onehalfspace}

\noindent This means that the function:

\begin{onehalfspace}
\begin{center}
{\small{
\begin{equation}
Corr_{n}\left(V_{i_{0}}\left(t\right),V_{i_{1}}\left(t\right),...,V_{i_{n-1}}\left(t\right)\right)\overset{def}{=}\frac{\mathbb{E}\left[{\displaystyle \prod_{j=0}^{n-1}}\left(V_{i_{j}}\left(t\right)-\overline{V}_{i_{j}}\left(t\right)\right)\right]}{\sqrt[n]{{\displaystyle \prod_{j=0}^{n-1}}\mathbb{E}\left[\left|V_{i_{j}}\left(t\right)-\overline{V}_{i_{j}}\left(t\right)\right|^{n}\right]}}\label{eq:higher-order-correlation-1}
\end{equation}
}}
\par\end{center}{\small \par}
\end{onehalfspace}

\noindent is in the range $\left[-1,1\right]$, therefore it is a
good formula to express higher order correlations. We can see that
for $n=2$ it gives the Pearson's formula, as it should be. Now, all
these means $\mathbb{E}$ can be computed using the Isserlis' theorem
as we did for the covariance, so in principle we can determine also
the higher order correlation structure of the neural network. However,
in practice, this gives rise to combinatorial problems with different
levels of complexity when $V_{i_{j}}\left(t\right)$ does not have
the same behavior for different values of $i_{j}$, namely if the
deterministic matrix $\overline{J}_{ij}+\sigma_{4}Z_{ij}\left(t\right)$
and the input vector $\overrightarrow{I}\left(t\right)$ do not have
strong symmetries. Therefore, for simplicity, in the Appendix \ref{sec:Higher order correlations for a fully connected neural network}
we show only the fully connected case with the same synaptic weights
and the same input current for all the neurons.

\section{\noindent \label{sec:Calculation of the fundamental matrix}Calculation
of the fundamental matrix}

\noindent As we can see from the formulae \ref{eq:covariance-part-1},
\ref{eq:covariance-part-2} and \ref{eq:covariance-part-3}, the correlation
structure is a function of the matrices $\Phi\left(t\right)$ and
$\Phi\left(t\right)\Phi^{T}\left(t\right)$. Therefore we need to
compute them for different kinds of connectivity matrices $\overline{J}$.
In general this is not an easy task, but however in some special cases
they can be obtained as discussed in Sections \ref{sub:Block circulant matrices with circulant blocks}
and \ref{sub:Symmetric matrices}.

\subsection{\noindent \label{sub:Block circulant matrices with circulant blocks}Block
circulant matrices with circulant blocks}

\noindent Given two positive integers $R$ and $S$, with $1\leq R,S\leq N$,
we suppose that $\overline{J}$ is an $N\times N$ block circulant
matrix (with $N=RS$) of the form:

\begin{onehalfspace}
\begin{center}
{\small{
\begin{equation}
\overline{J}=\frac{\Lambda}{M}\left[\begin{array}{cccc}
b^{\left(0\right)} & b^{\left(1\right)} & \cdots & b^{\left(R-1\right)}\\
b^{\left(R-1\right)} & b^{\left(0\right)} & \cdots & b^{\left(R-2\right)}\\
\vdots & \vdots & \ddots & \vdots\\
b^{\left(1\right)} & b^{\left(2\right)} & \cdots & b^{\left(0\right)}
\end{array}\right]\label{eq:block-circulant-matrix}
\end{equation}
}}
\par\end{center}{\small \par}
\end{onehalfspace}

\noindent where $b^{\left(0\right)},b^{\left(1\right)},...,b^{\left(R-1\right)}$
are $S\times S$ circulant matrices:

\begin{onehalfspace}
\begin{center}
{\small{
\begin{equation}
b^{\left(i\right)}=\left[\begin{array}{cccc}
b_{0}^{\left(i\right)} & b_{1}^{\left(i\right)} & \cdots & b_{S-1}^{\left(i\right)}\\
b_{S-1}^{\left(i\right)} & b_{0}^{\left(i\right)} & \cdots & b_{S-2}^{\left(i\right)}\\
\vdots & \vdots & \ddots & \vdots\\
b_{1}^{\left(i\right)} & b_{2}^{\left(i\right)} & \cdots & b_{0}^{\left(i\right)}
\end{array}\right]\label{eq:circulant-block}
\end{equation}
}}
\par\end{center}{\small \par}
\end{onehalfspace}

\begin{figure}
\begin{centering}
\includegraphics[scale=0.3]{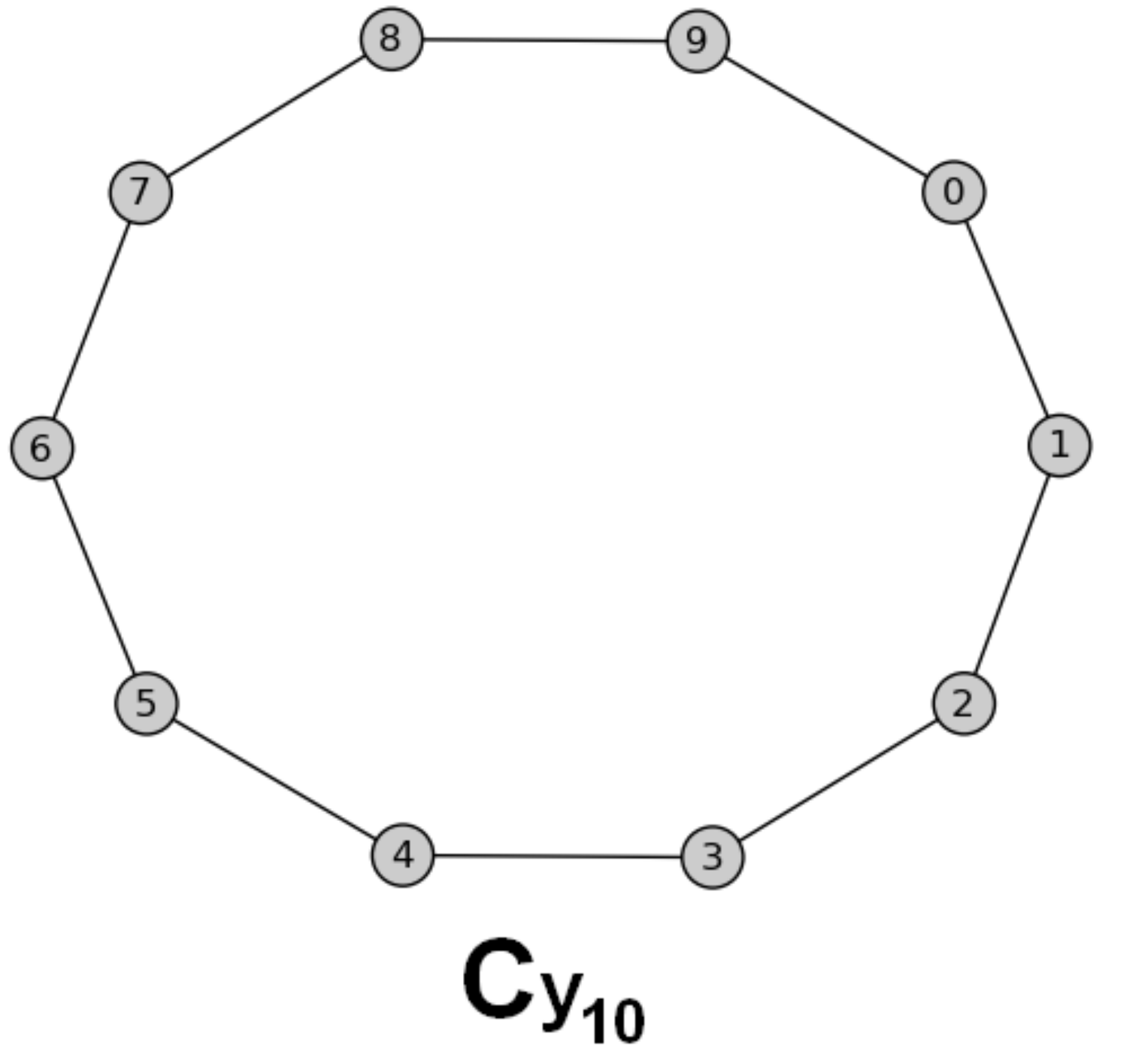}\includegraphics[scale=0.29]{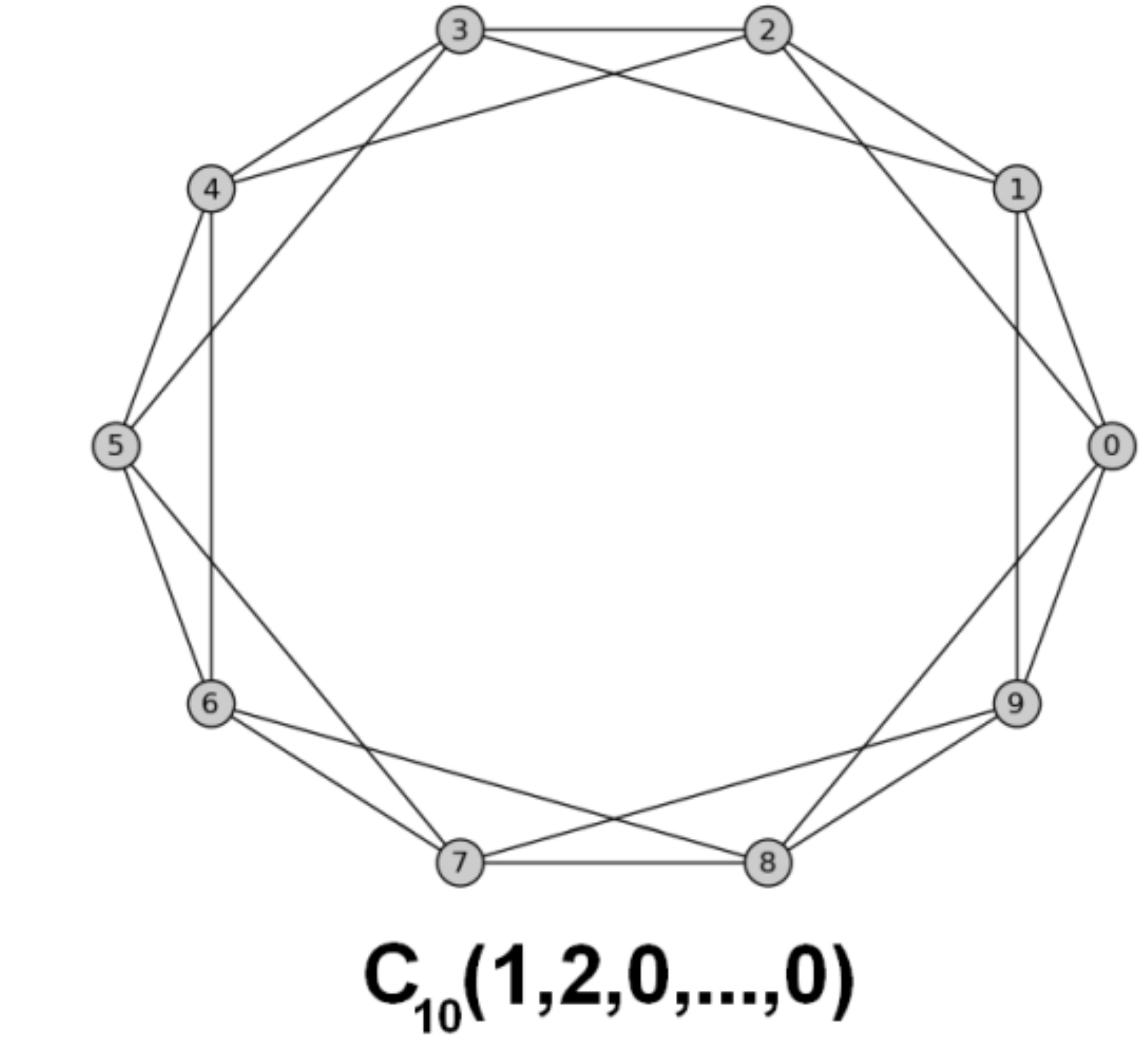}
\par\end{centering}

\begin{centering}
\includegraphics[scale=0.3]{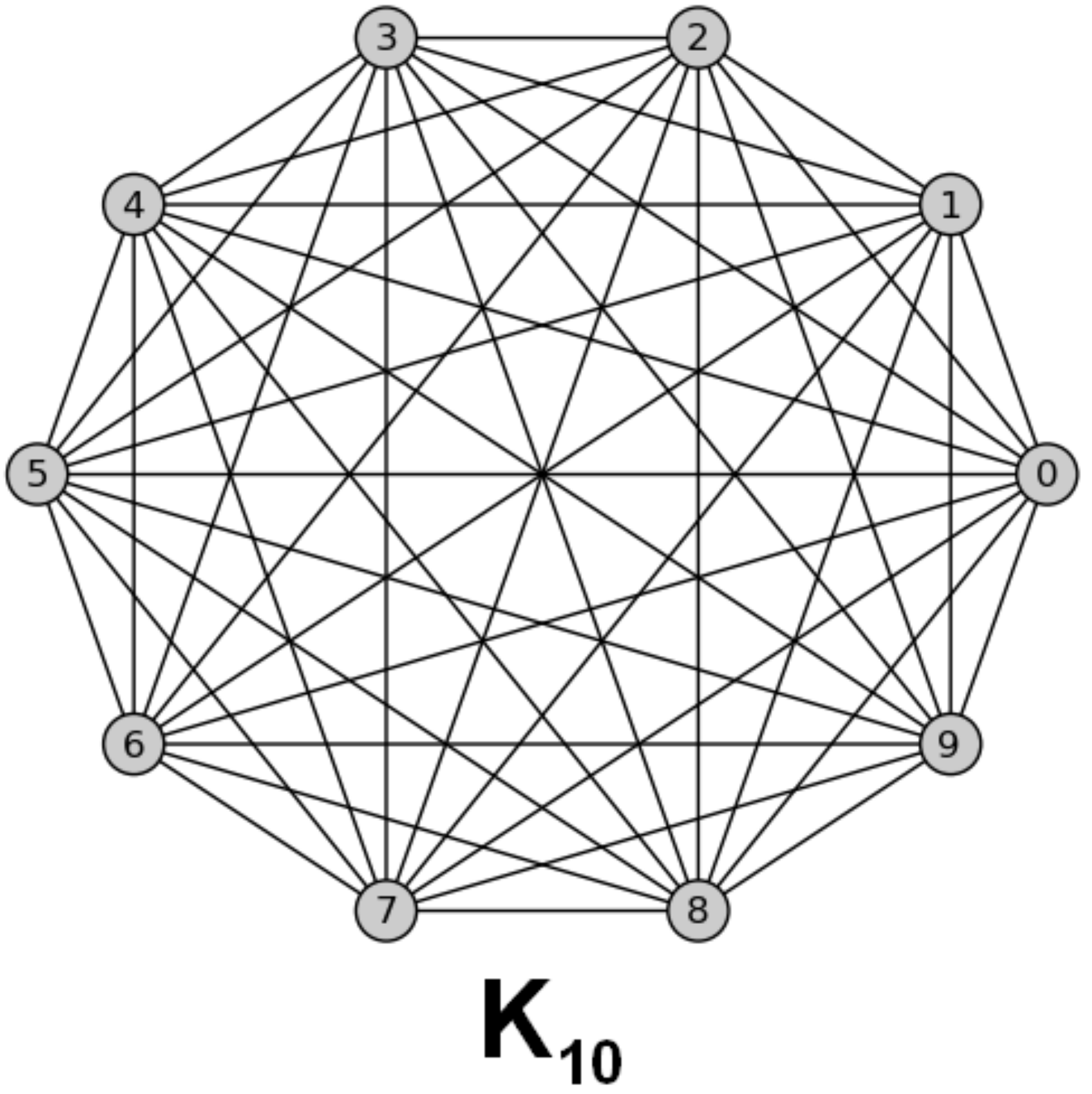}
\par\end{centering}

\caption[{\small{Circulant graphs}}]{{\small{}}{\footnotesize{\label{fig:circulant-graphs}Three examples
of circulant graphs: $Cy_{N}=C_{N}\left(1,0,0,...,0\right)$ (top-left),
also known as }}\textit{\footnotesize{cycle graph}}{\footnotesize{,
$C_{N}\left(1,2,0,...,0\right)$ (top-right) and $K_{N}=C_{N}\left(1,2,...,\left\lfloor \frac{N}{2}\right\rfloor \right)$
(bottom), also known as }}\textit{\footnotesize{complete graph}}{\footnotesize{
or }}\textit{\footnotesize{fully connected network}}{\footnotesize{.}}}
\end{figure}

\noindent All the entries $b_{j}^{\left(i\right)}$, for $i=0,1,...,R-1$
and $j=0,1,...,S-1$, can only be equal to $0$ or $1$, with only
the exception of $b_{0}^{\left(0\right)}$ that must always be equal
to $0$ in order to avoid the self-connections. $R$ can be interpreted
as the number of neural populations, and $S$ as the number of neurons
per population. Due to this particular structure of the connectivity
matrix, all the neurons have the same number of incoming synaptic
connections $M$, as required. This analysis includes the special
case when the matrix $\overline{J}$ is circulant (obtained for $R=1$
or $S=1$). In the context of \textit{Graph Theory}, a network whose
adjacency matrix is circulant is called \textit{circulant graph} (see
Figure \ref{fig:circulant-graphs}) and is usually represented by
the notation $C_{N}\left(1,2,...,q\right)$.

\noindent 

\noindent Moreover we have to recall that even if in Graph Theory
the connections are often represented through undirected unweighted
graphs, which means that the connectivity matrix is symmetric, in
this section we do not assume in general that $\overline{J}$ is symmetric.

\noindent Now we want to calculate the matrices $\Phi\left(t\right)$
and $\Phi\left(t\right)\Phi^{T}\left(t\right)$ in terms of the eigenquantities
of $\overline{J}$. The eigenvalues of $\overline{J}$ are the collection
of the eigenvalues of the following matrices:

\begin{onehalfspace}
\begin{center}
{\small{
\begin{equation}
\widetilde{b}^{\left(i\right)}=\sum_{j=0}^{R-1}e^{\frac{2\pi}{R}ij\iota}b^{\left(j\right)}\label{eq:discrete-Fourier-transform}
\end{equation}
}}
\par\end{center}{\small \par}
\end{onehalfspace}

\noindent where $\iota=\sqrt{-1}$. Since the matrices $\widetilde{b}^{\left(i\right)}$
are circulant, we can compute their eigenvalues $e_{j}^{\left(i\right)}$
as follows:

\begin{onehalfspace}
\begin{center}
{\small{
\begin{equation}
e_{j}^{\left(i\right)}=\sum_{k=0}^{S-1}e^{\frac{2\pi}{S}jk\iota}\left[\widetilde{b}^{\left(i\right)}\right]_{0k}=\sum_{k=0}^{S-1}\sum_{l=0}^{R-1}e^{2\pi\left(\frac{jk}{S}+\frac{il}{R}\right)\iota}b_{k}^{\left(l\right)}\label{eq:block-circulant-matrix-eigenvalues}
\end{equation}
}}
\par\end{center}{\small \par}
\end{onehalfspace}

\noindent Instead the matrix of the eigenvectors of $\overline{J}$
is:

\begin{onehalfspace}
\begin{center}
{\small{
\begin{align}
Q= & F_{R}\otimes F_{S}\nonumber \\
\label{eq:block-circulant-matrix-eigenvectors}\\
\left[F_{K}\right]_{ij}= & \begin{array}{ccccc}
\frac{1}{\sqrt{K}}e^{\frac{2\pi}{K}ij\iota}, &  & K=R,S, &  & i,j=0,1,...,K-1\end{array}\nonumber 
\end{align}
}}
\par\end{center}{\small \par}
\end{onehalfspace}

\noindent where $\otimes$ is the Kronecker product. Now, for $k=0,1,...,N-1$,
we call $a_{k}$ the eigenvalues of $A=-\frac{1}{\tau}Id_{N}+\overline{J}S'\left(\mu\right)$
(where $Id_{N}$ is the $N\times N$ identity matrix) and $e_{k}$
the eigenvalues of $\overline{J}$ (namely the collection of all the
$e_{j}^{\left(i\right)}$, with $k=iS+j$), while we call $\overrightarrow{v}_{k}$
and $\overrightarrow{w}_{k}$ their respective eigenvectors. Therefore
we have $a_{k}=-\frac{1}{\tau}+e_{k}S'\left(\mu\right)$ and $\overrightarrow{v}_{k}=\overrightarrow{w}_{k}$.
Moreover, using also the fact that the matrix $e^{At}$ can be diagonalized
and is real, we can write:

\begin{onehalfspace}
\begin{center}
{\small{
\begin{align*}
\Phi\left(t\right)= & e^{At}=QD\left(t\right)Q^{*}\\
\\
\Phi\left(t\right)\Phi^{T}\left(t\right)= & e^{At}\left(\left[e^{A\left(t\right)}\right]^{T}\right)^{*}=QD\left(t\right)Q^{*}QD^{*}\left(t\right)Q^{*}=QD\left(t\right)D^{*}\left(t\right)Q^{*}
\end{align*}
}}
\par\end{center}{\small \par}
\end{onehalfspace}

\noindent where $*$ is the element-by-element complex conjugation,
and $D\left(t\right)=diag\left(e^{a_{0}t},e^{a_{1}t}...,e^{a_{N-1}t}\right)$.
Here we have used the fact that $D\left(t\right)$ and $Q$ are symmetric
matrices and also the identity:

\begin{onehalfspace}
\begin{center}
{\small{
\[
Q^{*}Q=\left(F_{R}^{*}\otimes F_{S}^{*}\right)\left(F_{R}\otimes F_{S}\right)=\left(F_{R}^{*}F_{R}\right)\otimes\left(F_{S}^{*}F_{S}\right)=Id_{RS}=Id_{N}
\]
}}
\par\end{center}{\small \par}
\end{onehalfspace}

\noindent due to the mixed-product property of the Kronecker product
and to the elementary identity $F_{K}^{*}F_{K}=Id_{K}$. Now, since:

\begin{onehalfspace}
\begin{center}
{\small{
\[
\left[F_{R}\otimes F_{S}\right]_{ij}=\left[F_{R}\right]_{mn}\left[F_{S}\right]_{pq}=\frac{1}{\sqrt{N}}e^{2\pi\left(\frac{mn}{R}+\frac{pq}{S}\right)\iota}
\]
}}
\par\end{center}{\small \par}

\begin{center}
{\small{
\[
\begin{array}{ccccccc}
m=\left\lfloor \frac{i}{S}\right\rfloor , &  & n=\left\lfloor \frac{j}{S}\right\rfloor , &  & p=i-mS, &  & q=j-nS\end{array}
\]
}}
\par\end{center}{\small \par}
\end{onehalfspace}

\noindent we conclude that:

\begin{onehalfspace}
\begin{center}
{\small{
\begin{align}
\Phi_{ij}\left(t\right)= & \frac{1}{N}{\displaystyle \sum_{k=0}^{N-1}}e^{\left[-\frac{1}{\tau}+e_{k}S'\left(\mu\right)\right]t}f_{ijk}\nonumber \\
\label{eq:Phi-matrix-block-circulant-case}\\
\left[\Phi\left(t\right)\Phi^{T}\left(t\right)\right]_{ij}= & \frac{1}{N}{\displaystyle \sum_{k=0}^{N-1}}e^{2\left[-\frac{1}{\tau}+\Re\left(e_{k}\right)S'\left(\mu\right)\right]t}f_{ijk}\nonumber 
\end{align}
}}
\par\end{center}{\small \par}
\end{onehalfspace}

\noindent where $\Re\left(e_{k}\right)$ represents the real part
of $e_{k}$, while:

\begin{onehalfspace}
\noindent \begin{center}
{\small{
\[
f_{ijk}=\left[F_{R}\otimes F_{S}\right]_{ik}\left[F_{R}\otimes F_{S}\right]_{kj}^{*}=e^{2\pi\left\{ \frac{1}{R}\left\lfloor \frac{k}{S}\right\rfloor \left(\left\lfloor \frac{i}{S}\right\rfloor -\left\lfloor \frac{j}{S}\right\rfloor \right)+\frac{k}{S}\left(i-j\right)\right\} \iota}
\]
}}
\par\end{center}{\small \par}
\end{onehalfspace}

\noindent These formulae seem to give complex-valued functions, but
due to the particular structure of the eigenvalues $e_{k}$ and of
the function $f_{ijk}$, their imaginary parts are equal to zero (see
Appendix \ref{sec:Proof that formula 4.6 gives real functions}).
Therefore the covariance is a real function, as it should be.

\noindent Now we show an explicit example of this technique, namely
the case when the blocks of the matrix $\overline{J}$ have the following
symmetric circulant band structure:

\begin{onehalfspace}
\noindent \begin{center}
{\scriptsize{
\begin{equation}
b^{\left(i\right)}=\left[\begin{array}{cccccccccc}
1-\delta_{i0} & 1 & \cdots & 1 & 0 & \cdots & 0 & 1 & \cdots & 1\\
1 & 1-\delta_{i0} & \ddots &  & \ddots & \ddots &  & \ddots & \ddots & \vdots\\
\vdots & \ddots & \ddots & \ddots &  & \ddots & \ddots & 0 & \ddots & 1\\
1 &  & \ddots & \ddots & \ddots &  & \ddots & \ddots &  & 0\\
0 & \ddots &  & \ddots & \ddots & \ddots &  & \ddots & \ddots & \vdots\\
\vdots &  & \ddots &  & \ddots & \ddots & \ddots &  & \ddots & 0\\
0 &  &  & \ddots &  & \ddots & \ddots & \ddots &  & 1\\
1 & \ddots & 0 &  & \ddots &  & \ddots & \ddots & \ddots & \vdots\\
\vdots & \ddots & \ddots &  &  & \ddots &  & \ddots & 1-\delta_{i0} & 1\\
1 & \cdots & 1 & 0 & \cdots & 0 & 1 & \cdots & 1 & 1-\delta_{i0}
\end{array}\right]\label{eq:symmetric-circulant-band-matrix}
\end{equation}
}}
\par\end{center}{\scriptsize \par}
\end{onehalfspace}

\noindent where, supposing for simplicity that $S\geq3$, the first
row of $b^{\left(i\right)}$ (excluding the term $\left[b^{\left(i\right)}\right]_{00}$,
which is $0$ for $i=0$ and $1$ for $i>0$) can be written explicitly
as:

\begin{onehalfspace}
\begin{center}
{\small{
\[
\begin{cases}
\left[b^{\left(i\right)}\right]_{0j}= & \begin{array}{ccc}
1, &  & \left(1\leq j\leq\nu_{i}\right)\vee\left(\rho_{i}\leq j\leq S-1\right)\end{array}\\
\\
\left[b^{\left(i\right)}\right]_{0j}= & \begin{array}{ccc}
0, &  & \nu_{i}<j<\rho_{i}\end{array}
\end{cases}
\]
}}
\par\end{center}{\small \par}

\begin{center}
{\small{
\begin{align*}
\rho_{i}= & S-\nu_{i}+H\left(\nu_{i}-\left\lfloor \frac{S}{2}\right\rfloor +\left(-1\right)^{S}\right)\\
\\
H\left(x\right)= & \begin{cases}
0, & \begin{array}{cc}
 & x\leq0\end{array}\\
\\
1, & \begin{array}{cc}
 & x>0\end{array}
\end{cases}
\end{align*}
}}
\par\end{center}{\small \par}
\end{onehalfspace}

\noindent with $1\leq\nu_{i}\leq\left\lfloor \frac{S}{2}\right\rfloor $.
Here we have to suppose that $S\geq3$ because otherwise it is not
possible to distinguish the diagonal band from the corner elements.
Now, the bandwidth of $b^{\left(i\right)}$ is $2\nu_{i}+1$, so this
defines the integer parameters $\nu_{i}$. Moreover, $2\nu_{0}-H\left(\nu_{0}-\left\lfloor \frac{S}{2}\right\rfloor +\left(-1\right)^{S}\right)$
represents the number of connections that every neuron in a given
population receives from the neurons in the same population. Instead
$2\nu_{i}+1-H\left(\nu_{i}-\left\lfloor \frac{S}{2}\right\rfloor +\left(-1\right)^{S}\right)$,
for $i=1,2,...,R-1$, is the number of connections that every neuron
in the the $k$-th population receives from the neurons in the $\left(i+k\right)$-th
mod $R$ population, for $k=0,1,...,R-1$. So the total number of
incoming connections per neuron is $M=R-1+\sum_{i=0}^{R-1}\left[2\nu_{i}-H\left(\nu_{i}-\left\lfloor \frac{S}{2}\right\rfloor +\left(-1\right)^{S}\right)\right]$.
It is important to observe that even if all the matrices $b^{\left(i\right)}$
are symmetric, the matrix $\overline{J}$ in general is not, since
the number of connections in every block is different (the case of
symmetric connectivity matrices is studied in Section \ref{sub:Symmetric matrices}).
Now, using formula \ref{eq:block-circulant-matrix-eigenvalues}, we
obtain that:

\begin{onehalfspace}
\begin{center}
{\small{
\begin{align}
e_{mS+n}= & \begin{cases}
\frac{\Lambda}{M}\left[R-1+{\displaystyle \sum_{k=0}^{R-1}}f\left(n,\nu_{k},S\right)\right], & \begin{array}{cc}
 & m=0,\,\forall n\end{array}\\
\\
\frac{\Lambda}{M}\left[-1+{\displaystyle \sum_{k=0}^{R-1}}e^{\frac{2\pi}{R}mk\iota}f\left(n,\nu_{k},S\right)\right], & \begin{array}{cc}
 & m\neq0,\,\forall n\end{array}
\end{cases}\nonumber \\
\label{eq:symmetric-circulant-band-matrix-eigenvalues}\\
f\left(n,\nu_{k},S\right)= & \begin{cases}
2\nu_{k}-H\left(\nu_{k}-\left\lfloor \frac{S}{2}\right\rfloor +\left(-1\right)^{S}\right), & \begin{array}{cc}
 & n=0,\:\forall\nu_{k}\end{array}\\
\\
-1, & \begin{array}{cc}
 & n\neq0,\:\nu_{k}=\left\lfloor \frac{S}{2}\right\rfloor \end{array}\\
\\
\frac{sin\left(\frac{\pi n\left(2\nu_{k}+1\right)}{S}\right)}{sin\left(\frac{\pi n}{S}\right)}-1, & \begin{array}{cc}
 & n\neq0,\:\nu_{k}<\left\lfloor \frac{S}{2}\right\rfloor \end{array}
\end{cases}\nonumber 
\end{align}
}}
\par\end{center}{\small \par}
\end{onehalfspace}

\noindent with $m=0,1,...,R-1$ and $n=0,1,...,S-1$.

\noindent Many different special cases can be studied. The simplest
one is obtained for $\nu_{0}=\nu_{1}=...=\nu_{R-1}\overset{def}{=}\nu$,
and in this case formula \ref{eq:symmetric-circulant-band-matrix-eigenvalues}
gives:

\begin{onehalfspace}
\begin{center}
{\small{
\begin{equation}
e_{mS+n}=\begin{cases}
\frac{\Lambda}{M}\left[R-1+Rf\left(n,\nu,S\right)\right], & \begin{array}{cc}
 & m=0,\,\forall n\end{array}\\
\\
-\frac{\Lambda}{M}, & \begin{array}{cc}
 & m\neq0,\,\forall n\end{array}
\end{cases}\label{eq:symmetric-circulant-band-matrix-eigenvalues-simplified}
\end{equation}
}}
\par\end{center}{\small \par}
\end{onehalfspace}

\noindent with $M=R-1+R\left[2\nu-H\left(\nu-\left\lfloor \frac{S}{2}\right\rfloor +\left(-1\right)^{S}\right)\right]$.
Therefore in this case the eigenvalues are real, as it must be, since
with this special choice of the parameters the matrix $\overline{J}$
is symmetric. For $R=1$ and $\nu<\left\lfloor \frac{N}{2}\right\rfloor $
we have $M=2\nu$ and formula \ref{eq:symmetric-circulant-band-matrix-eigenvalues-simplified}
gives the eigenvalues of the circulant network:

\begin{onehalfspace}
\begin{center}
{\small{
\begin{equation}
e_{n}=\begin{cases}
\Lambda, & \begin{array}{cc}
 & n=0\end{array}\\
\\
\frac{\Lambda}{2\nu}\left[\frac{sin\left(\frac{\pi n\left(2\nu+1\right)}{N}\right)}{sin\left(\frac{\pi n}{N}\right)}-1\right], & \begin{array}{cc}
 & n\neq0\end{array}
\end{cases}\label{eq:circulant-network-eigenvalues}
\end{equation}
}}
\par\end{center}{\small \par}
\end{onehalfspace}

\noindent Instead for $\nu=\left\lfloor \frac{S}{2}\right\rfloor $
and $\forall R,\, S$ we have $M=N-1$ and formula \ref{eq:symmetric-circulant-band-matrix-eigenvalues-simplified}
gives the eigenvalues of the fully connected network:

\begin{onehalfspace}
\begin{center}
{\small{
\begin{equation}
e_{n}=\begin{cases}
\Lambda, & \begin{array}{cc}
 & n=0\end{array}\\
\\
-\frac{\Lambda}{N-1}, & \begin{array}{cc}
 & n\neq0\end{array}
\end{cases}\label{eq:complete-network-eigenvalues}
\end{equation}
}}
\par\end{center}{\small \par}
\end{onehalfspace}

\subsection{\noindent \label{sub:Symmetric matrices}Symmetric matrices}

\noindent Another case where the matrices $\Phi\left(t\right)$ and
$\Phi\left(t\right)\Phi^{T}\left(t\right)$ can be computed easily
is when we have a general symmetric matrix $\overline{J}$. Since
its entries are real, it can be diagonalized by an orthogonal matrix
$Q$ (namely such that $Q^{-1}=Q^{T}$), therefore we have:

\begin{onehalfspace}
\begin{center}
{\small{
\begin{align*}
\overline{J}= & Q\widetilde{D}Q^{T}\\
\\
\widetilde{D}= & diag\left(\widetilde{d}_{1},\widetilde{d}_{2},...,\widetilde{d}_{N-1}\right)
\end{align*}
}}
\par\end{center}{\small \par}
\end{onehalfspace}

\noindent So we obtain:

\begin{onehalfspace}
\begin{center}
{\small{
\begin{align*}
A= & -\frac{1}{\tau}Id_{N}+\overline{J}S'\left(\mu\right)=Q\left[-\frac{1}{\tau}Id_{N}+\widetilde{D}S'\left(\mu\right)\right]Q^{T}\\
\\
\Phi\left(t\right)= & e^{At}=Qe^{\left[-\frac{1}{\tau}Id_{N}+\widetilde{D}S'\left(\mu\right)\right]t}Q^{T}=QD\left(t\right)Q^{T}
\end{align*}
}}
\par\end{center}{\small \par}
\end{onehalfspace}

\noindent having defined the diagonal matrix $D\left(t\right)$ as
follows:

\begin{onehalfspace}
\begin{center}
{\small{
\begin{align*}
D\left(t\right)= & e^{\left[-\frac{1}{\tau}Id_{N}+\widetilde{D}S'\left(\mu\right)\right]t}=diag\left(d_{1},d_{2},...,d_{N-1}\right)\\
\\
d_{i}= & e^{\left[-\frac{1}{\tau}+\widetilde{d}_{i}S'\left(\mu\right)\right]t}
\end{align*}
}}
\par\end{center}{\small \par}
\end{onehalfspace}

\noindent Moreover, also the matrix $A$ is symmetric in this case,
therefore:

\begin{onehalfspace}
\begin{center}
{\small{
\[
\Phi\left(t\right)\Phi^{T}\left(t\right)=e^{2At}=QD^{2}\left(t\right)Q^{T}
\]
}}
\par\end{center}{\small \par}
\end{onehalfspace}

\noindent so their components are:

\begin{onehalfspace}
\begin{center}
{\small{
\begin{align}
\Phi_{ij}\left(t\right)= & {\displaystyle \sum_{k=0}^{N-1}}e^{D_{k}\left(t-s\right)}Q_{ik}Q_{jk}\nonumber \\
\label{eq:Phi-matrix-symmetric-case}\\
\left[\Phi\left(t\right)\Phi^{T}\left(t\right)\right]_{ij}= & {\displaystyle \sum_{k=0}^{N-1}}e^{2D_{k}\left(t-s\right)}Q_{ik}Q_{jk}\nonumber 
\end{align}
}}
\par\end{center}{\small \par}
\end{onehalfspace}

\begin{figure}
\begin{centering}
\includegraphics[scale=0.3]{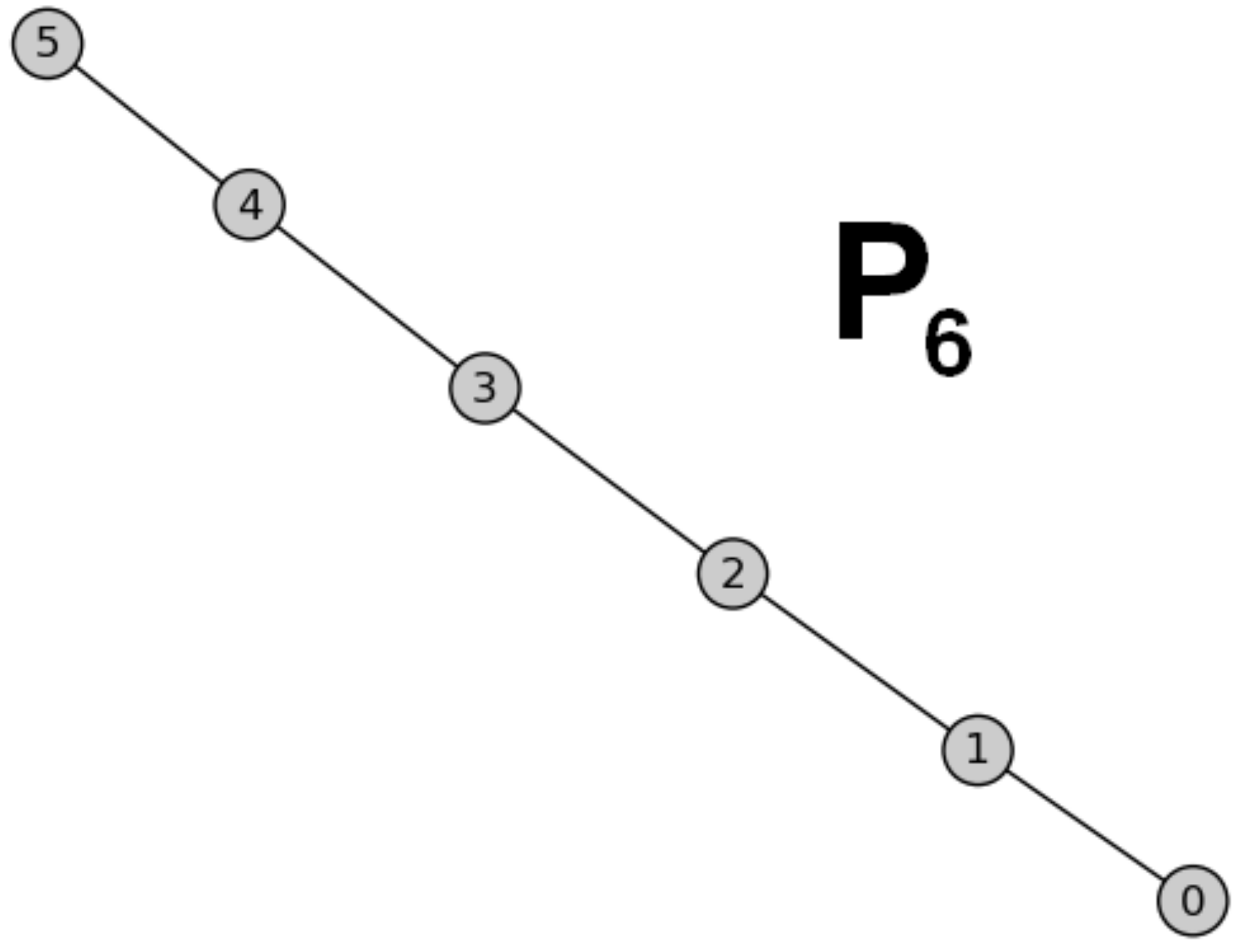}
\par\end{centering}

\caption[{\footnotesize{Path graph}}]{{\footnotesize{\label{fig:path-graph}Example of the graph $P_{N}$,
known as}}\textit{\footnotesize{ path on $N$ nodes}}{\footnotesize{.
Its connectivity matrix is tridiagonal without corner elements.}}}
\end{figure}

\noindent Again, now we need only the eigenquantities of $\overline{J}$,
but it is not possible to find explicit expressions for a general
symmetric connectivity matrix. Actually they can be calculated analytically
only if $\overline{J}$ has some special kind of structure. However,
since it is symmetric and all its non-zero entries have the same value
$\frac{\Lambda}{M}$ (as we said in Section \ref{sec:Description of the model}),
it can be interpreted as the adjacency matrix of an undirected unweighted
graph. Due to this correspondence, we can study the eigenquantities
of $\overline{J}$ using the powerful techniques already developed
in the context of Graph Theory for this kind of graphs. Lee and Yeh
\cite{LeeYeh1993} have proved that it is possible to perform binary
operations (in particular the Kronecker product $\otimes$ and the
Cartesian product $\times$) on pairs of graphs $G_{1}$ and $G_{2}$,
obtaining more complicated graphs, whose eigenvalues and eigenvectors
can be calculated easily from those of the graphs $G_{1}$ and $G_{2}$.
If $e_{G}$ and $\overrightarrow{v}_{G}$ represent respectively the
eigenvalues and eigenvectors of the graph $G$, then we obtain:

\begin{onehalfspace}
\begin{center}
{\small{
\begin{align}
 & \begin{cases}
e_{G_{1}\otimes G_{2}}^{i,j}= & e_{G_{1}}^{i}e_{G_{2}}^{j}\\
\\
\overrightarrow{v}_{G_{1}\otimes G_{2}}^{i,j}= & \overrightarrow{v}_{G_{1}}^{i}\otimes\overrightarrow{v}_{G_{2}}^{j}
\end{cases}\label{eq:eigenquantities-Kronecker-product}\\
\nonumber \\
 & \begin{cases}
e_{G_{1}\times G_{2}}^{i,j}= & e_{G_{1}}^{i}+e_{G_{2}}^{j}\\
\\
\overrightarrow{v}_{G_{1}\times G_{2}}^{i,j}= & \overrightarrow{v}_{G_{1}}^{i}\otimes\overrightarrow{v}_{G_{2}}^{j}
\end{cases}\label{eq:eigenquantities-Cartesian-product}
\end{align}
}}
\par\end{center}{\small \par}
\end{onehalfspace}

\noindent for $i,j=0,1,...,N-1$. In particular we can choose $G_{1}$
and $G_{2}$ to be $P_{N}$ and/or $Cy_{N}$, where $P_{N}$ is the
so called\textit{ path on $N$ nodes} (see Figure \ref{fig:path-graph}),
while $Cy_{N}$ is the cycle graph (see Figure \ref{fig:circulant-graphs}).

\noindent Their eigenquantities (in the case of unitary weights) are:

\begin{onehalfspace}
\begin{center}
{\small{
\begin{align}
 & \begin{cases}
e_{P_{N}}^{i}= & 2\cos\left[\frac{\left(i+1\right)\pi}{N+1}\right]\\
\\
\left[\overrightarrow{v}_{P_{N}}^{i}\right]_{j}= & \sin\left[\frac{\left(i+1\right)\left(j+1\right)\pi}{N+1}\right]
\end{cases}\label{eq:eigenquantities-path-graph}\\
\nonumber \\
 & \begin{cases}
e_{Cy_{N}}^{i}= & 2cos\left(\frac{2\pi i}{N}\right)\\
\\
\left[\overrightarrow{v}_{Cy_{N}}^{i}\right]_{j}= & e^{\frac{2\pi ij}{N}\iota}
\end{cases},\,\iota=\sqrt{-1}\label{eq:eigenquantities-cyclic-graph}
\end{align}
}}
\par\end{center}{\small \par}
\end{onehalfspace}

\begin{figure}
\begin{centering}
\includegraphics[scale=0.3]{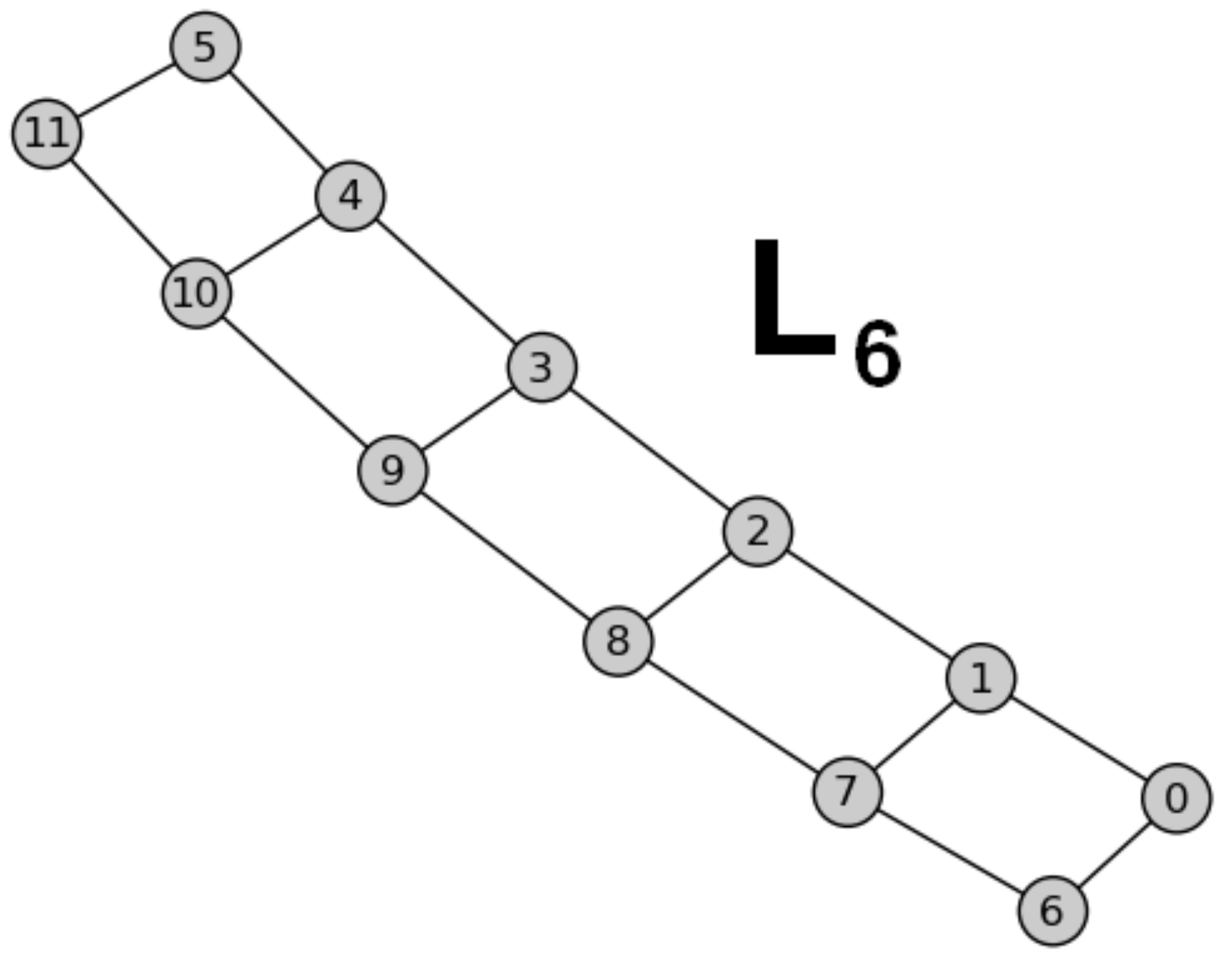}\includegraphics[scale=0.3]{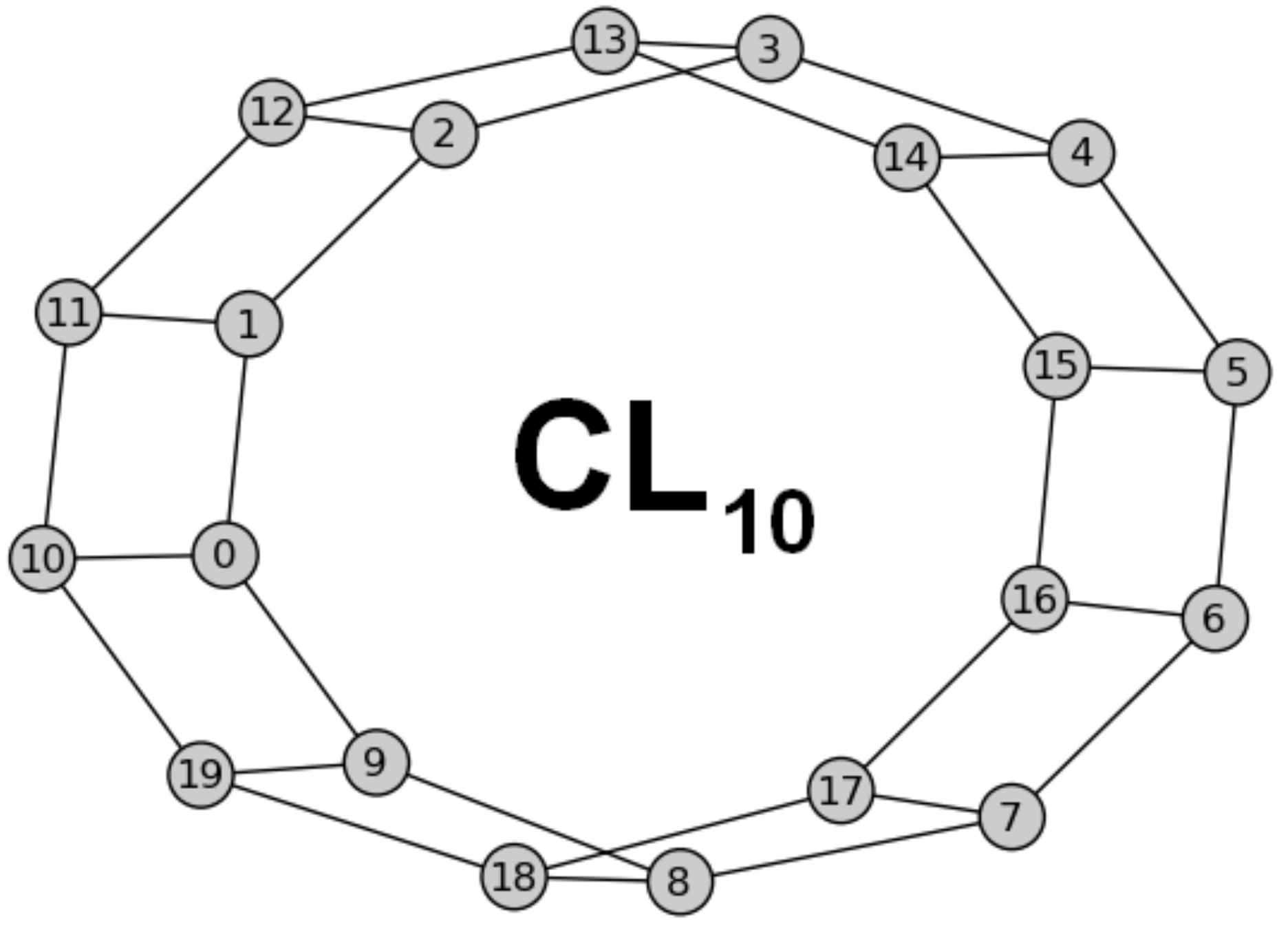}
\par\end{centering}

\begin{centering}
\includegraphics[scale=0.3]{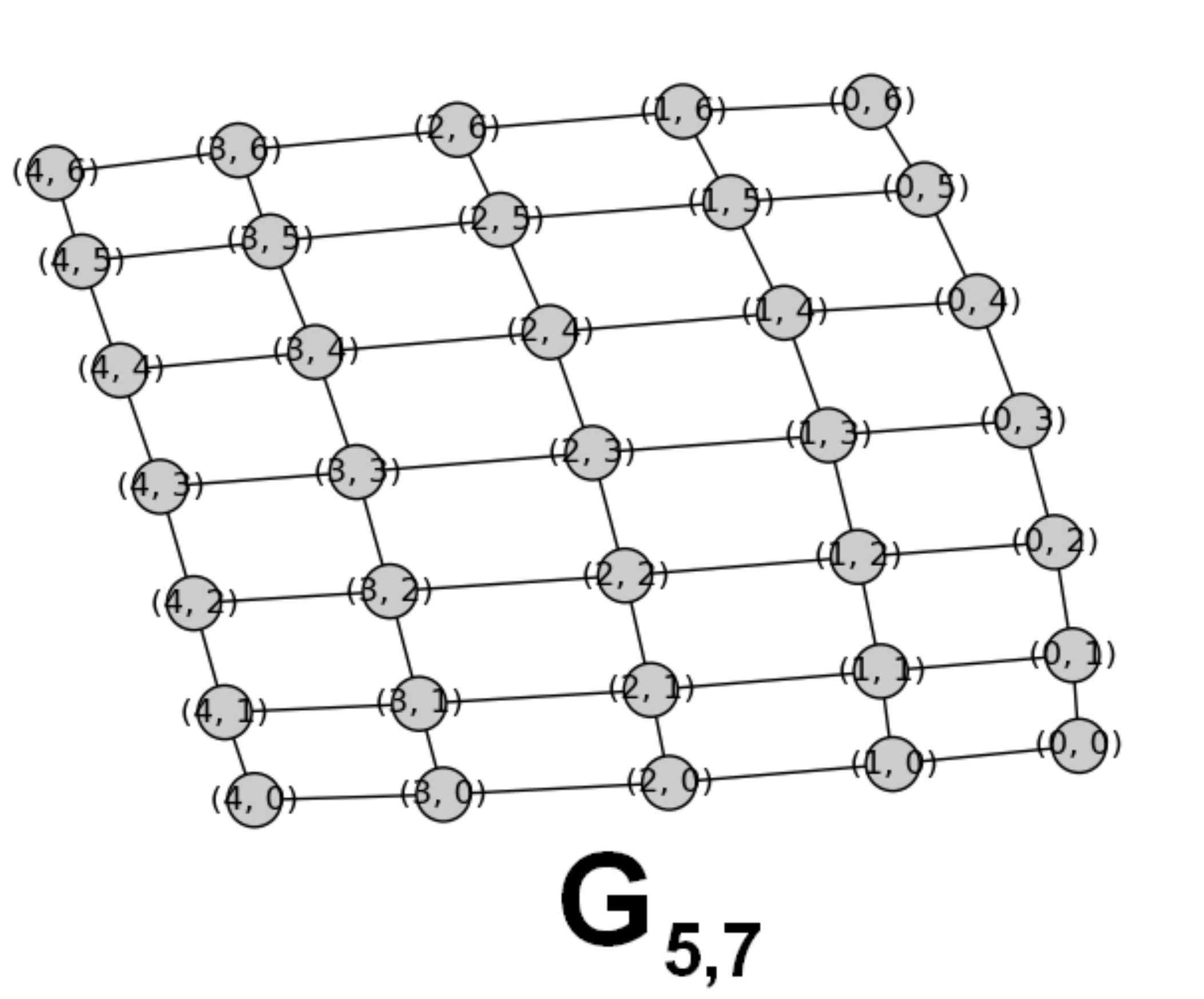}\includegraphics[scale=0.3]{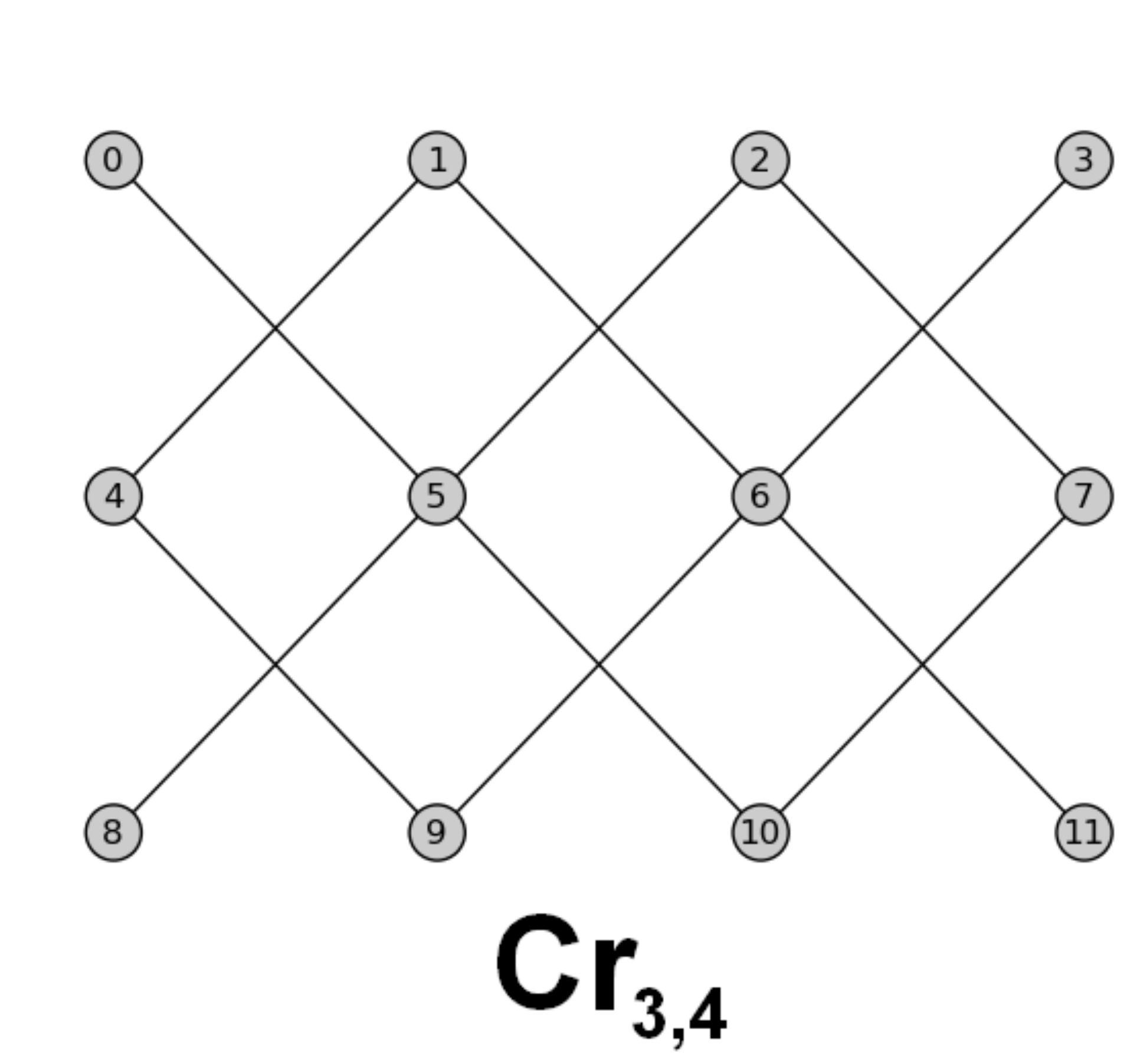}
\par\end{centering}

\caption[{\footnotesize{Examples of graphs}}]{{\footnotesize{\label{fig:graph-examples}Some examples of graphs:
the Ladder $L_{n}=P_{n}\times P_{2}$ (top left), the Circular Ladder
$CL_{n}=Cy_{n}\times P_{2}$ (top right), the Grid $G_{m,n}=P_{m}\times P_{n}$
(bottom left) and the Cross $Cr_{m,n}=P_{m}\otimes P_{n}$ (bottom
right).}}}
\end{figure}

\begin{figure}
\begin{centering}
\includegraphics[scale=0.3]{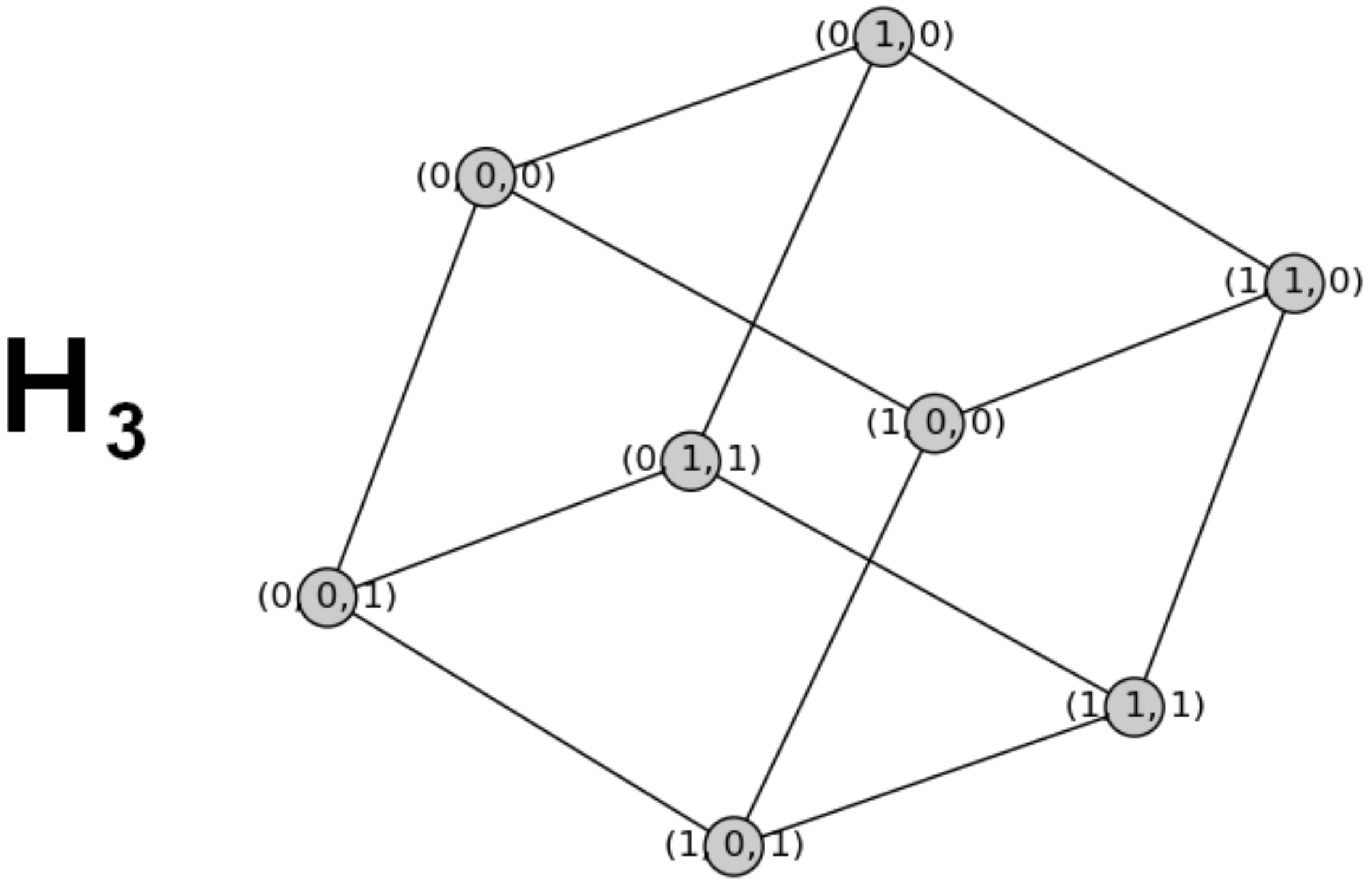}\includegraphics[scale=0.3]{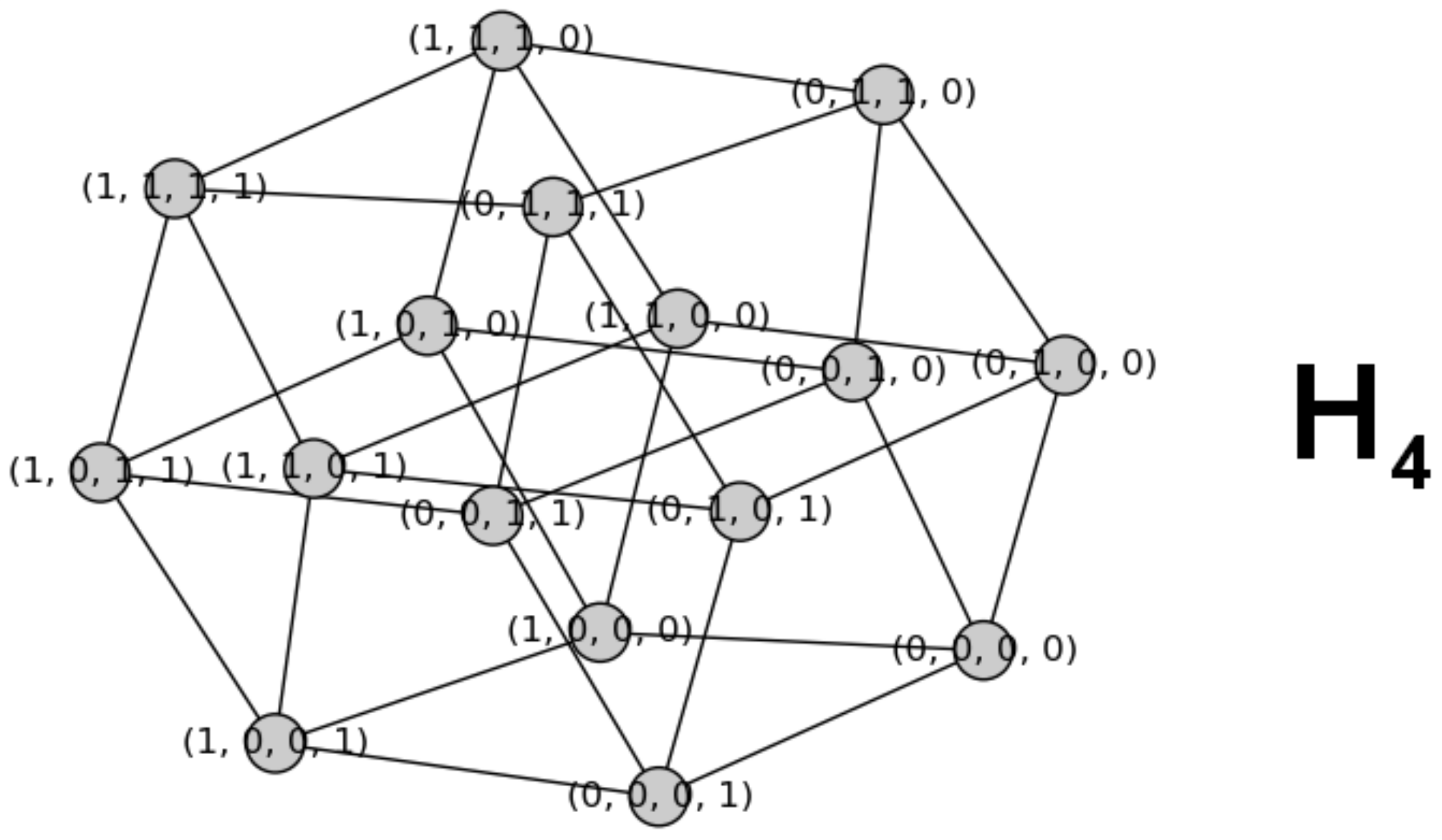}
\par\end{centering}

\begin{centering}
\includegraphics[scale=0.28]{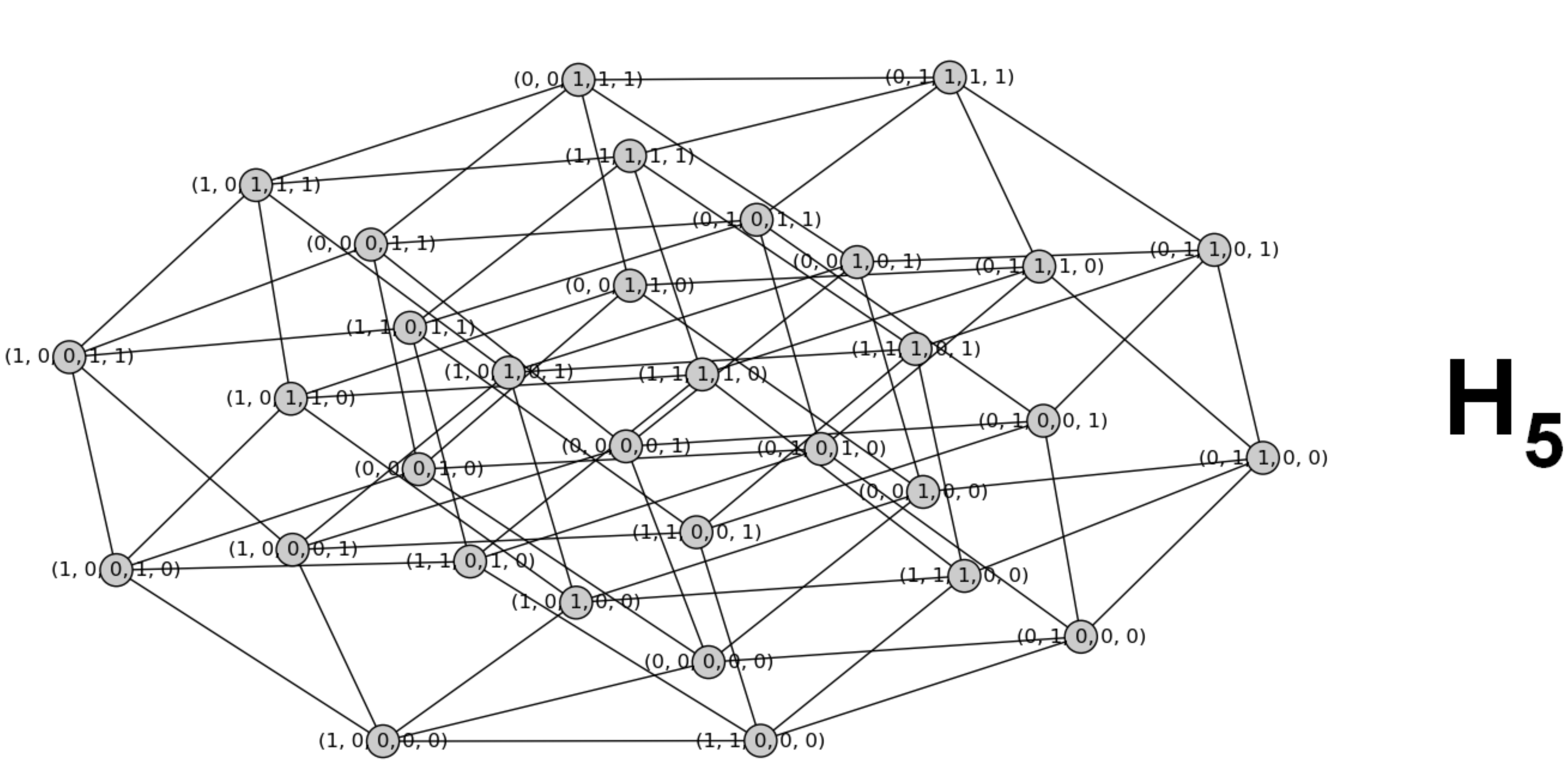}
\par\end{centering}

\caption[{\footnotesize{Hypercube graph}}]{{\footnotesize{\label{fig:hypercube-graph} Three examples of the
Hypercube graph $H_{n}$.}}}
\end{figure}

\noindent Combining them through the binary operations $\otimes$
and $\times$, we can create several classes of well-known graphs,
like:
\begin{itemize}
\item \noindent Ladder: $L_{n}=P_{n}\times P_{2}$, with $2n=N$;
\item \noindent Circular Ladder (also known as Annulus or Prism): $CL_{n}=Cy_{n}\times P_{2}$,
with $2n=N$;
\item \noindent Grid: $G_{m,n}=P_{m}\times P_{n}$, with $mn=N$;
\item \noindent Cylinder: $Cl_{m,n}=P_{m}\times Cy_{n}$, with $mn=N$;
\item \noindent Torus: $T_{m,n}=Cy_{m}\times Cy_{n}$, with $mn=N$;
\item \noindent Cross: $Cr_{m,n}=P_{m}\otimes P_{n}$, with $mn=N$;
\item \noindent Hypercube: $H_{n}=\underset{n-times}{\underbrace{P_{2}\times P_{2}\times...\times P_{2}}}$,
with $2^{n}=N$;
\end{itemize}
\noindent and so on and so forth. Some of these examples are shown
in the Figures \ref{fig:graph-examples} and \ref{fig:hypercube-graph}.
Even much more complicated graphs can be created in this way, like
$Cr_{m,n}\otimes T_{p,q}$, or $G_{m,n}\otimes T_{p,q}\times H_{r}\otimes Cl_{x,y}$,
and so on.

\noindent Using the mixed-product property of the Kronecker product,
we obtain:

\begin{onehalfspace}
\begin{center}
{\small{
\[
\overrightarrow{v}_{P_{N_{1}}\otimes P_{N_{2}}}^{i,j}\cdot\overrightarrow{v}_{P_{N_{1}}\otimes P_{N_{2}}}^{k,l}=\left(\overrightarrow{v}_{P_{N_{1}}}^{i}\cdot\overrightarrow{v}_{P_{N_{1}}}^{k}\right)\left(\overrightarrow{v}_{P_{N_{2}}}^{j}\cdot\overrightarrow{v}_{P_{N_{2}}}^{l}\right)=0
\]
}}
\par\end{center}{\small \par}
\end{onehalfspace}

\noindent if $i\neq k$ and/or $j\neq l$, since the eigenvectors
of the path are orthogonal (of course the same is true for $\overrightarrow{v}_{P_{N_{1}}\times P_{N_{2}}}^{i,j}\cdot\overrightarrow{v}_{P_{N_{1}}\times P_{N_{2}}}^{k,l}$).
Therefore the eigenvectors $\overrightarrow{v}_{P_{N_{1}}\otimes P_{N_{2}}}^{i,j}$
(or equivalently $\overrightarrow{v}_{P_{N_{1}}\times P_{N_{2}}}^{i,j}$
) are orthogonal. Moreover $\overrightarrow{v}_{P_{N_{1}}\otimes P_{N_{2}}}^{i,j}$
have real entries, therefore they form an orthogonal matrix, that
can be used directly to compute $\Phi\left(t\right)$ and $\Phi\left(t\right)\Phi^{T}\left(t\right)$
through formula \ref{eq:Phi-matrix-symmetric-case}.

\noindent For the eigenvectors $\overrightarrow{v}_{Cy_{N_{1}}\otimes Cy_{N_{2}}}^{i,j}$
the procedure is slightly more complicated, since the eigenvectors
$\overrightarrow{v}_{Cy_{N}}^{i}$ have in general complex entries
(with only the exception of the cases $i=0$ and $i=\frac{N}{2}$
for $N$ even) and therefore we cannot use them to form an orthogonal
matrix. However, if $\overline{J}$ is the connectivity matrix corresponding
to the graph $Cy_{N}$, we have:

\begin{onehalfspace}
\begin{center}
\[
\overline{J}\overrightarrow{v}_{Cy_{N}}^{i}=e_{Cy_{N}}^{i}\overrightarrow{v}_{Cy_{N}}^{i}
\]

\par\end{center}
\end{onehalfspace}

\noindent which implies:

\begin{onehalfspace}
\begin{center}
\[
\overline{J}^{*}\left(\overrightarrow{v}_{Cy_{N}}^{i}\right)^{*}=\left(e_{Cy_{N}}^{i}\right)^{*}\left(\overrightarrow{v}_{Cy_{N}}^{i}\right)^{*}
\]

\par\end{center}
\end{onehalfspace}

\noindent where $*$ is the element-by-element complex conjugation.
Since $\overline{J}$ and $e_{Cy_{N}}^{i}$ are real, we obtain that
$\overrightarrow{v}_{Cy_{N}}^{i}$ and $\left(\overrightarrow{v}_{Cy_{N}}^{i}\right)^{*}$
are both eigenvectors of $\overline{J}$, corresponding to the same
eigenvalue $e_{Cy_{N}}^{i}$. Therefore, if for all the complex eigenvectors
$\overrightarrow{v}_{Cy_{N}}^{i}$ we define the new vectors:

\begin{onehalfspace}
\begin{center}
{\small{
\begin{align}
\overrightarrow{V}_{Cy_{N}}^{i}= & \frac{1}{2}\left(\overrightarrow{v}_{Cy_{N}}^{i}+\left[\overrightarrow{v}_{Cy_{N}}^{i}\right]^{*}\right)\nonumber \\
\label{eq:eigenvectors-linear-combination}\\
\overrightarrow{W}_{Cy_{N}}^{i}= & \frac{1}{2\iota}\left(\overrightarrow{v}_{Cy_{N}}^{i}-\left[\overrightarrow{v}_{Cy_{N}}^{i}\right]^{*}\right)\nonumber 
\end{align}
}}
\par\end{center}{\small \par}
\end{onehalfspace}

\noindent we conclude that they are eigenvectors of $\overline{J}$
with eigenvalue $e_{Cy_{N}}^{i}$. Now, it is easy to see that $\overrightarrow{V}_{Cy_{N}}^{i}\cdot\overrightarrow{W}_{Cy_{N}}^{j}=0$
$\forall i,j$. Moreover $\overrightarrow{V}_{Cy_{N}}^{i}$ and $\overrightarrow{W}_{Cy_{N}}^{i}$
are orthogonal also to $\overrightarrow{V}_{Cy_{N}}^{0}$ and $\overrightarrow{V}_{Cy_{N}}^{\frac{N}{2}}$
in the case of $N$ even, and their entries are real. Therefore, if
we use this set of real eigenvectors with the rules \ref{eq:eigenquantities-Kronecker-product}
or \ref{eq:eigenquantities-Cartesian-product}, we obtain a set of
eigenvectors for $Cy_{N_{1}}\otimes Cy_{N_{2}}$ or $Cy_{N_{1}}\times Cy_{N_{2}}$
which are orthogonal and real (the proof is similar to the case $P_{N_{1}}\otimes P_{N_{2}}$
seen before). So they can be used to form an orthogonal matrix, through
which we can compute $\Phi\left(t\right)$ and $\Phi\left(t\right)\Phi^{T}\left(t\right)$,
according to formula \ref{eq:Phi-matrix-symmetric-case}.

\noindent To conclude, taking for example the cases we have written
previously, it is important to observe that only the graphs $CL_{n}$,
$T_{m,n}$ and $H_{n}$ can be considered in our analysis. In effect
these three graphs have the same number of incoming connections per
neuron, a feature that is not shown by the ladder, the grid etc, due
to their boundaries. The latter graphs can be studied using this approach
only in the thermodynamic limit $N\rightarrow\infty$. In fact only
in this case the number of incoming connections per neuron is the
same for all the neurons in the network, because when $N\rightarrow\infty$
the system \textquotedbl{}loses\textquotedbl{} its boundaries, since
they are pushed to infinity. For example, in the graph $L_{n}$ all
the neurons have $3$ incoming connections (see neurons $1-4$ and
$7-10$ in the case of the graph $L_{6}$ shown in Figure \ref{fig:graph-examples}),
with only the exception of those at the boundaries, which have only
$2$ connections (see neurons $0$, $5$, $6$ and $11$ in Figure
\ref{fig:graph-examples}). When $N\rightarrow\infty$, if we start
to travel on the ladder from its center toward the boundaries, we
will never reach them, since they are at infinity, therefore during
the trip we meet only neurons with the same number (namely $3$) of
incoming connections. Therefore in the thermodynamic limit all the
neurons of the graphs with boundaries behave in the same way. This
means that we have obtained the invariance of the system under exchange
of the neural indices, which is what we need in order to apply the
perturbative approach introduced in this article.

\section{\noindent \label{sec:Numerical comparison}Numerical comparison}

\begin{figure}
\begin{centering}
{\footnotesize{\includegraphics[scale=0.3]{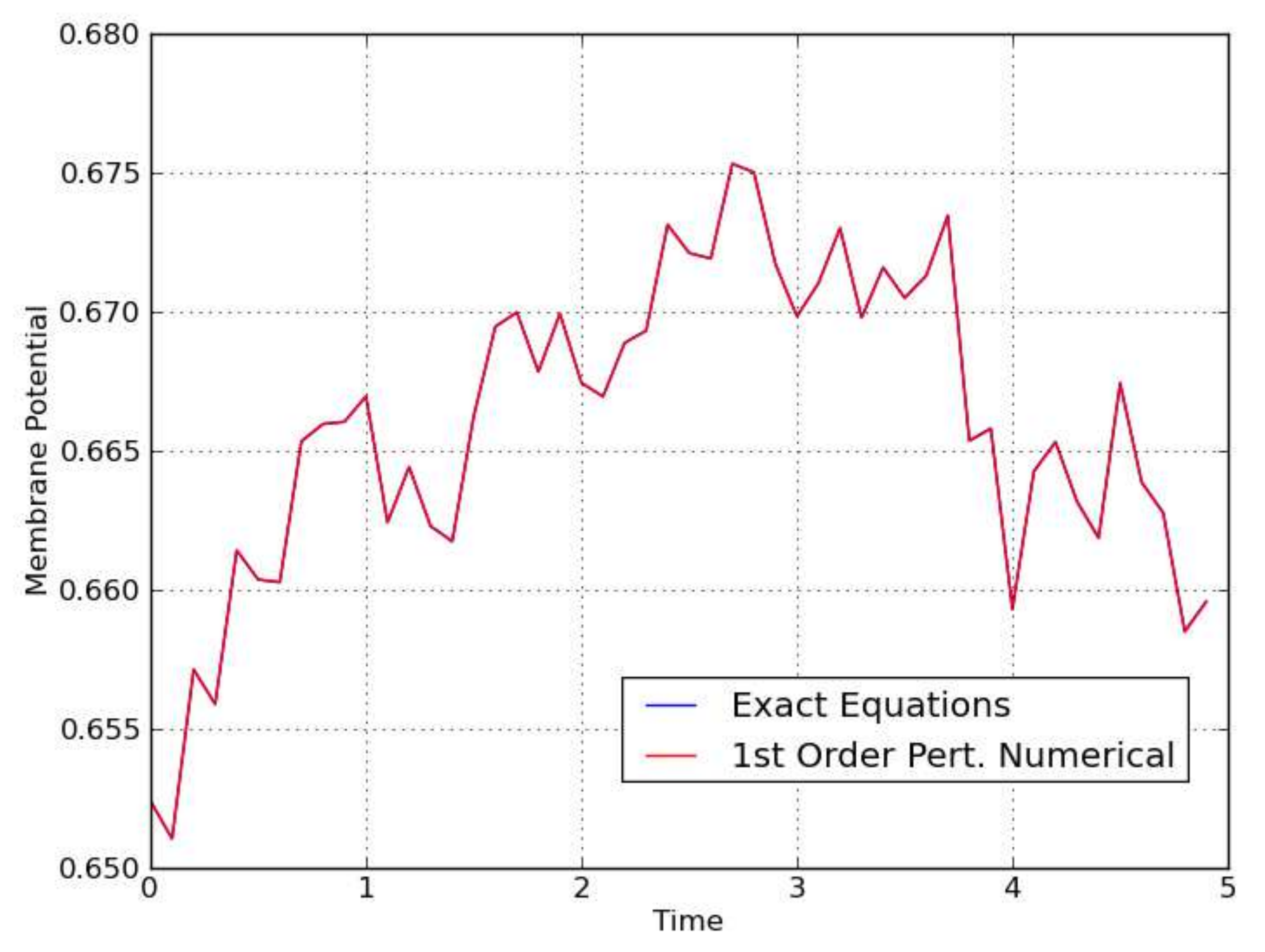}\includegraphics[scale=0.3]{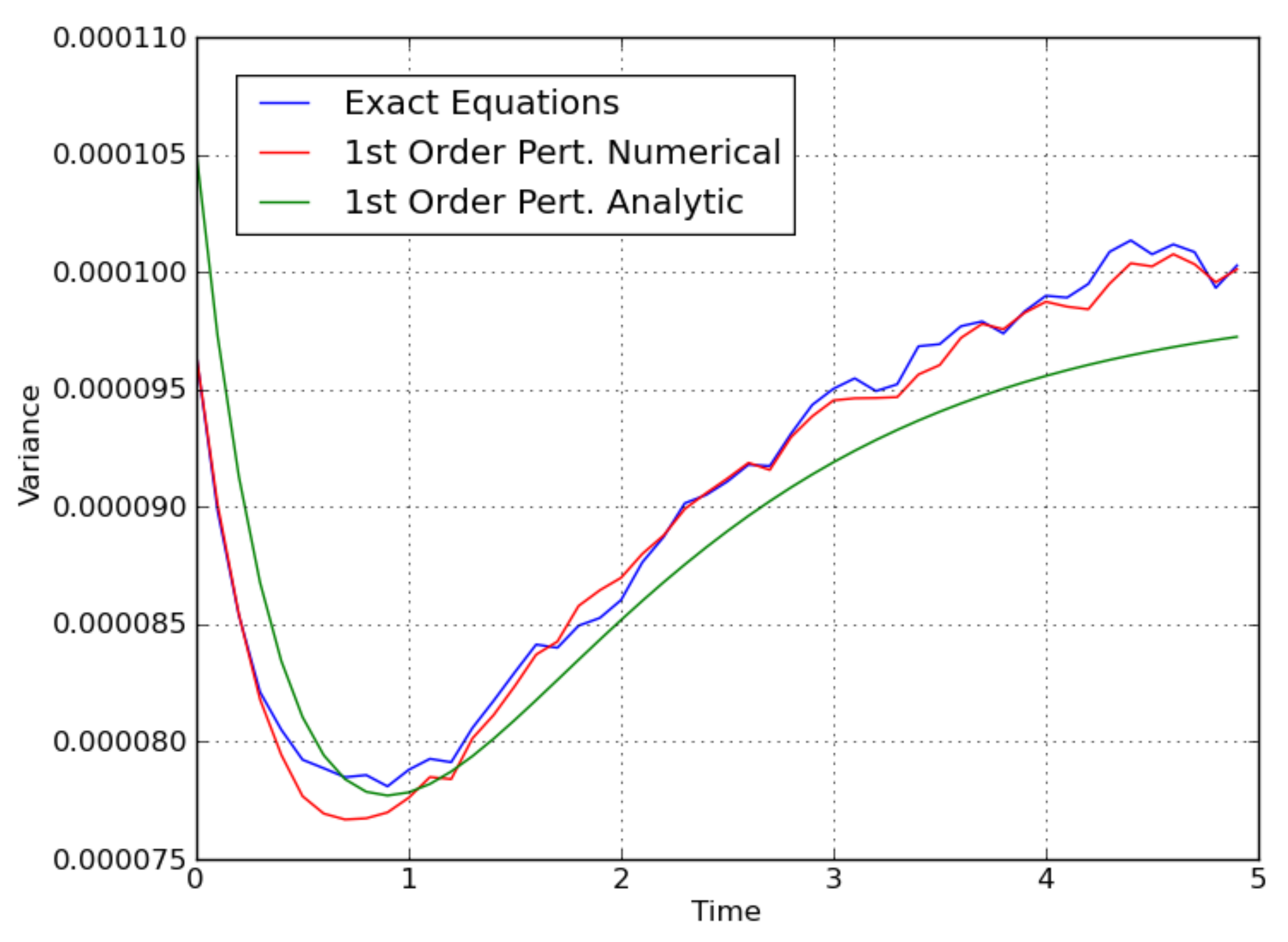}}}
\par\end{centering}{\footnotesize \par}

\begin{centering}
{\footnotesize{\includegraphics[scale=0.3]{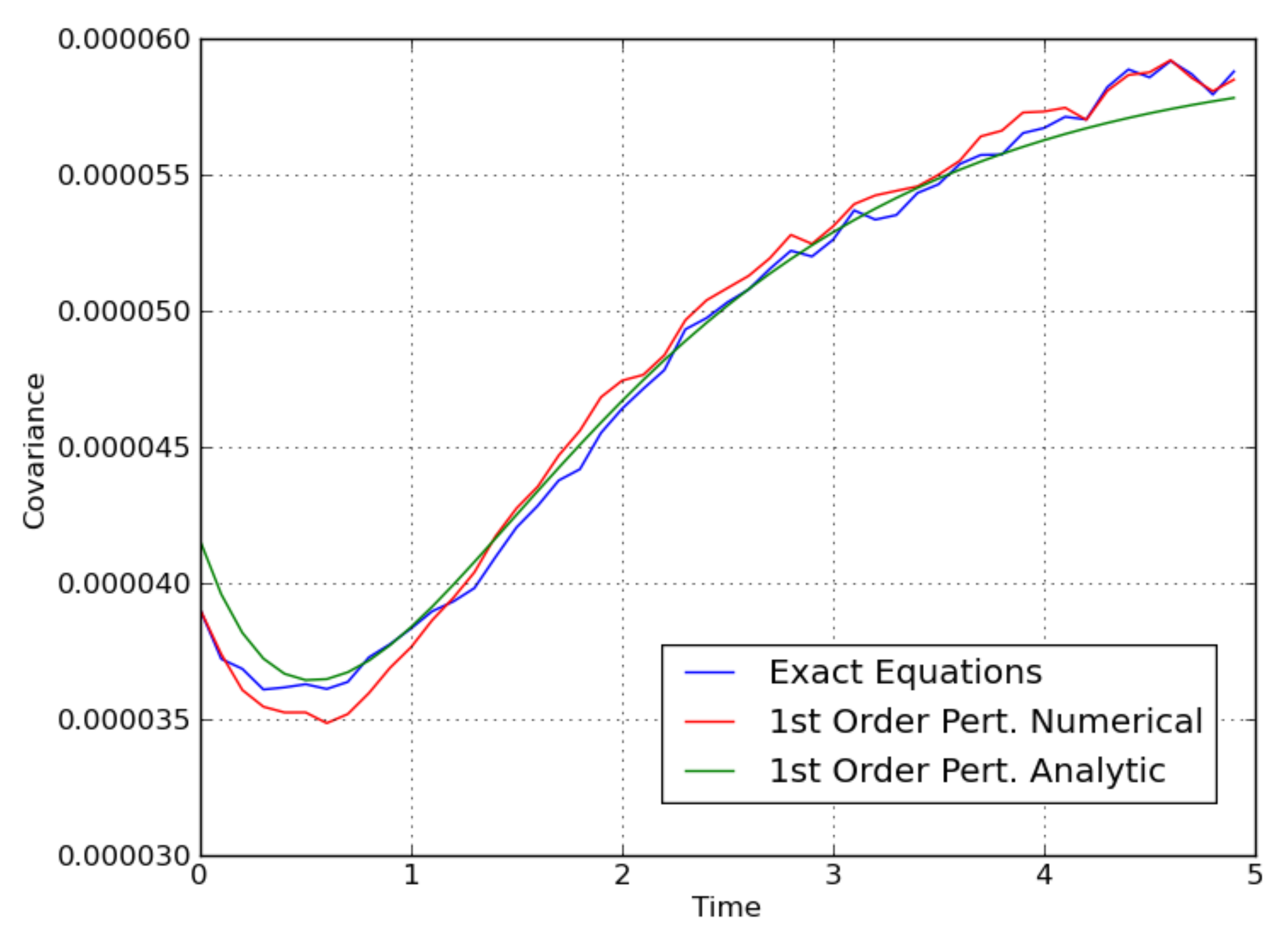}\includegraphics[scale=0.3]{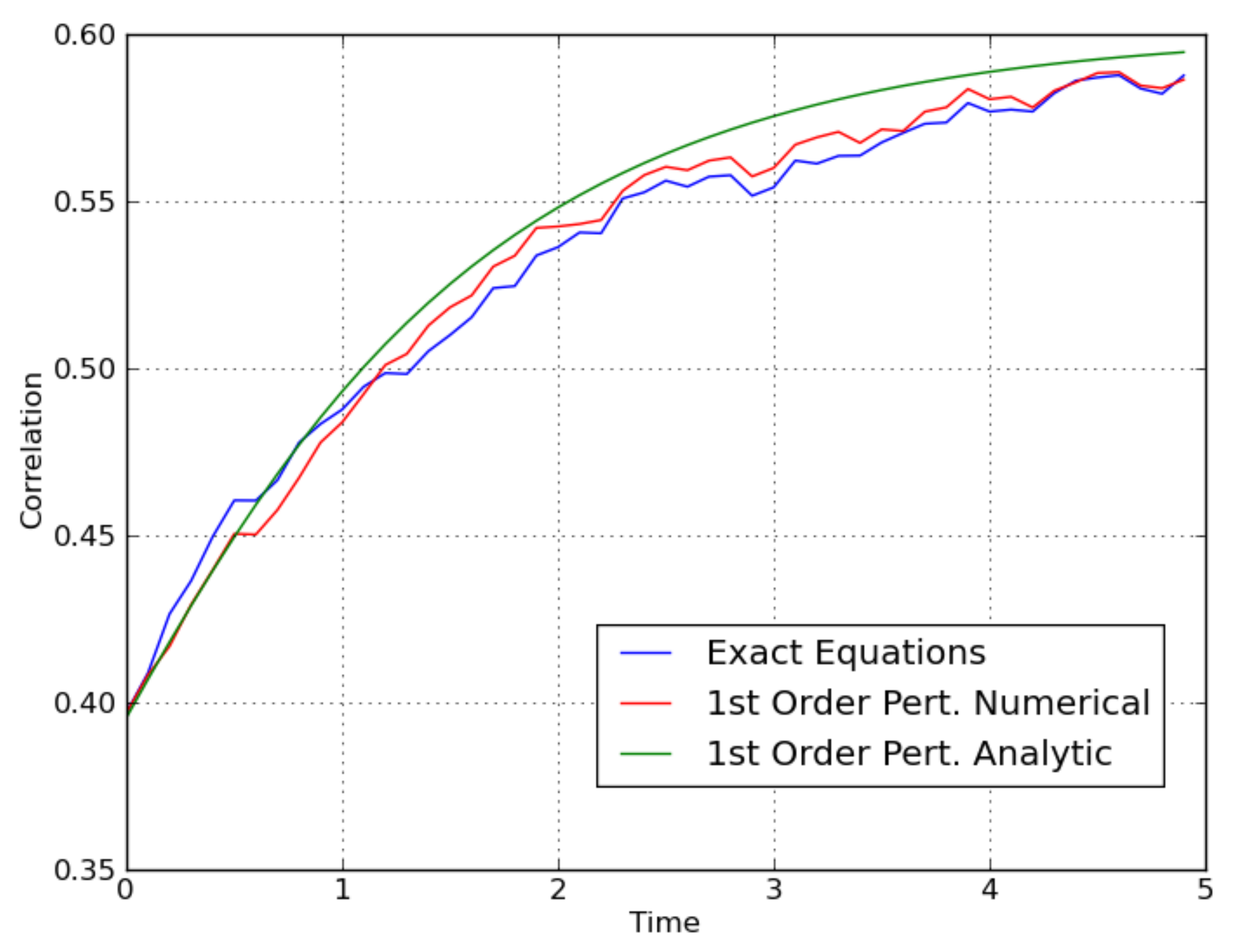}}}
\par\end{centering}{\footnotesize \par}

{\footnotesize{\caption[{\footnotesize{Numerical comparison of the perturbative expansion
with strong weights - 1}}]{{\footnotesize{\label{fig:circular-ladder-graph-simulation-1}First-order
perturbative expansion ($\sigma_{4}=\sigma_{5}=0$) for a network
with connectivity matrix $CL_{10}$ (namely $N=20$). These results
have been obtained for the values of the parameters reported in Table
\ref{tab:simulation-parameters-1}, for $\sigma_{1}=\sigma_{2}=\sigma_{3}=0.01$
and with the statistics evaluated through $10,000$ Monte Carlo simulations.
Since $\sigma_{1}$, $\sigma_{2}$ and $\sigma_{3}$ are small, in
the picture of the membrane potentials $V_{i}\left(t\right)$ (top-left)
there is a perfect agreement between the result obtained from the
exact network equations \ref{eq:rate-model-exact-equations-2} (blue
line) and that obtained from the first-order perturbative expansion,
namely from the equations \ref{eq:perturbative-equation-1} - \ref{eq:perturbative-equation-4}
(red line). Instead the comparison between the variances (top-right),
covariances (bottom-left) and correlations (bottom-right) is less
good because small values of $\sigma_{1}$, $\sigma_{2}$ and $\sigma_{3}$
determine small values of the variance and covariance, therefore a
higher number of Monte Carlo simulations is required in order to improve
the match. The green line represents the analytic result obtained
for the first-order perturbative expansion for an infinite number
of Monte Carlo simulations (formulae \ref{eq:covariance} - \ref{eq:variance}),
therefore it is the limit curve reached by the red line when the number
of simulations is increased indefinitely. The blue and red lines have
been obtained numerically by solving the corresponding equations with
the Euler-Maruyama scheme, while the integrals with respect to time
involved in the formulae for the evaluation of the green line have
been calculated with the trapezoidal rule. In all the cases the integration
time step is $\Delta t=0.1$.}}}
}}
\end{figure}

\begin{figure}
\begin{centering}
\includegraphics[scale=0.3]{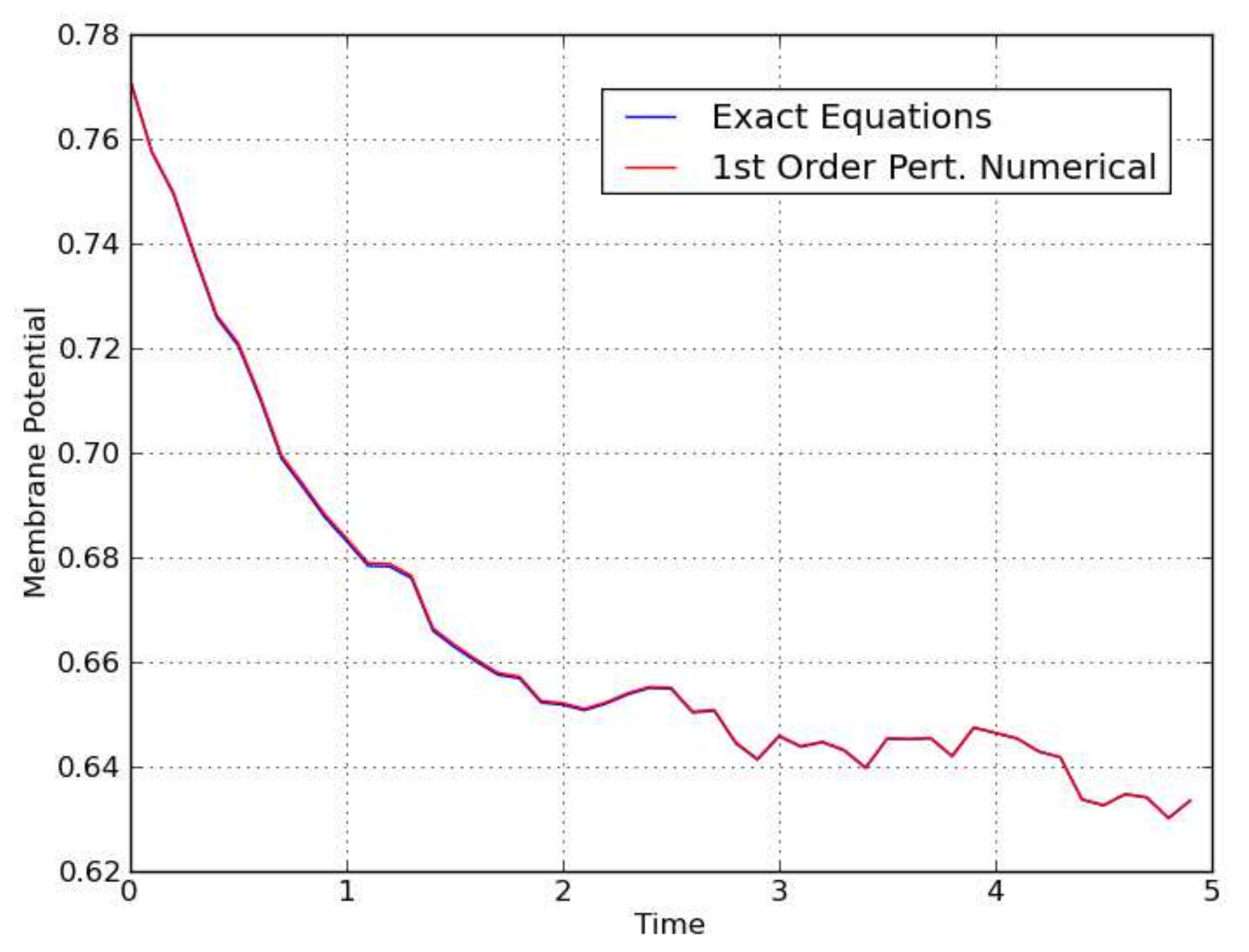}\includegraphics[scale=0.3]{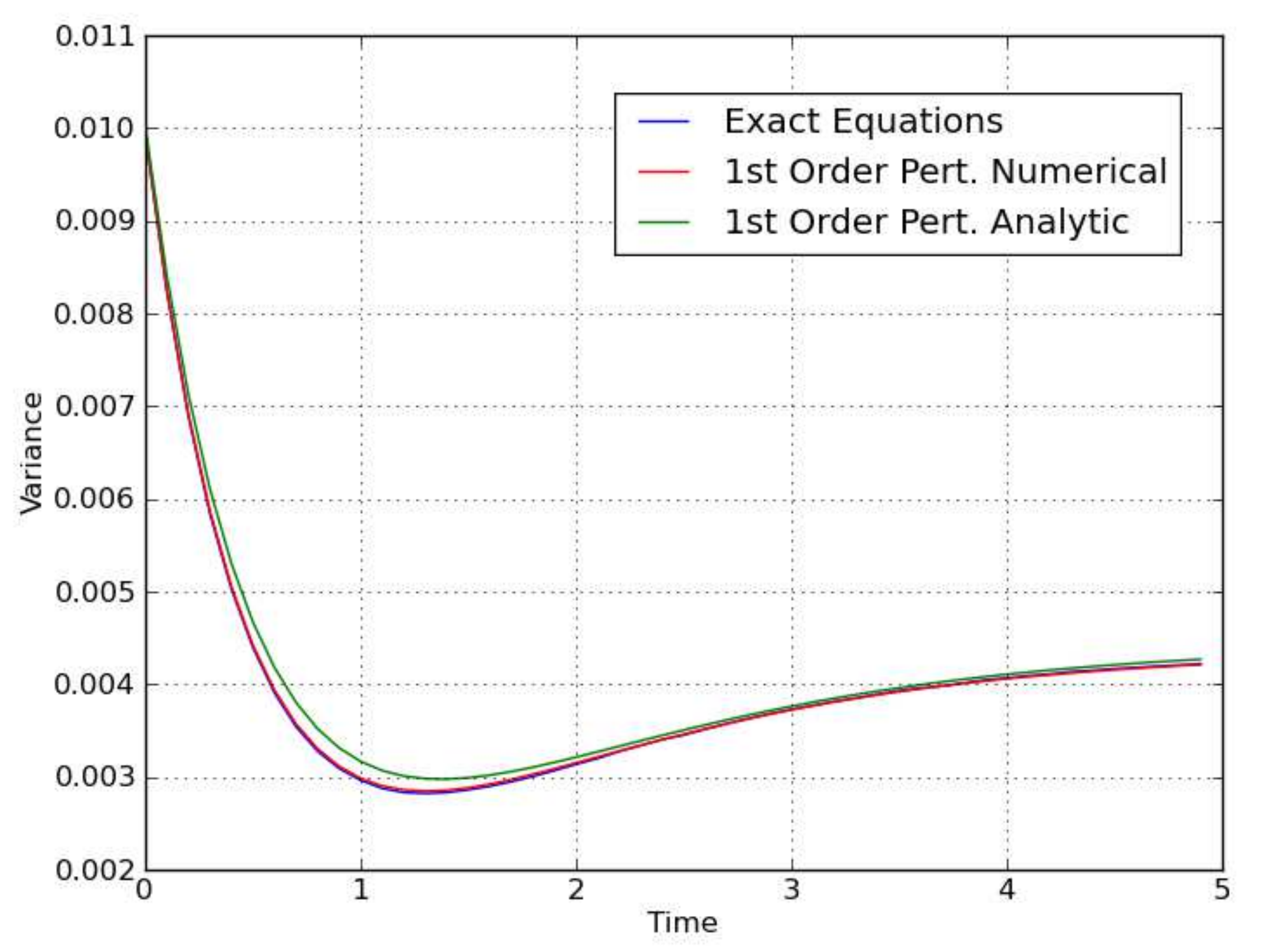}
\par\end{centering}

\begin{centering}
\includegraphics[scale=0.3]{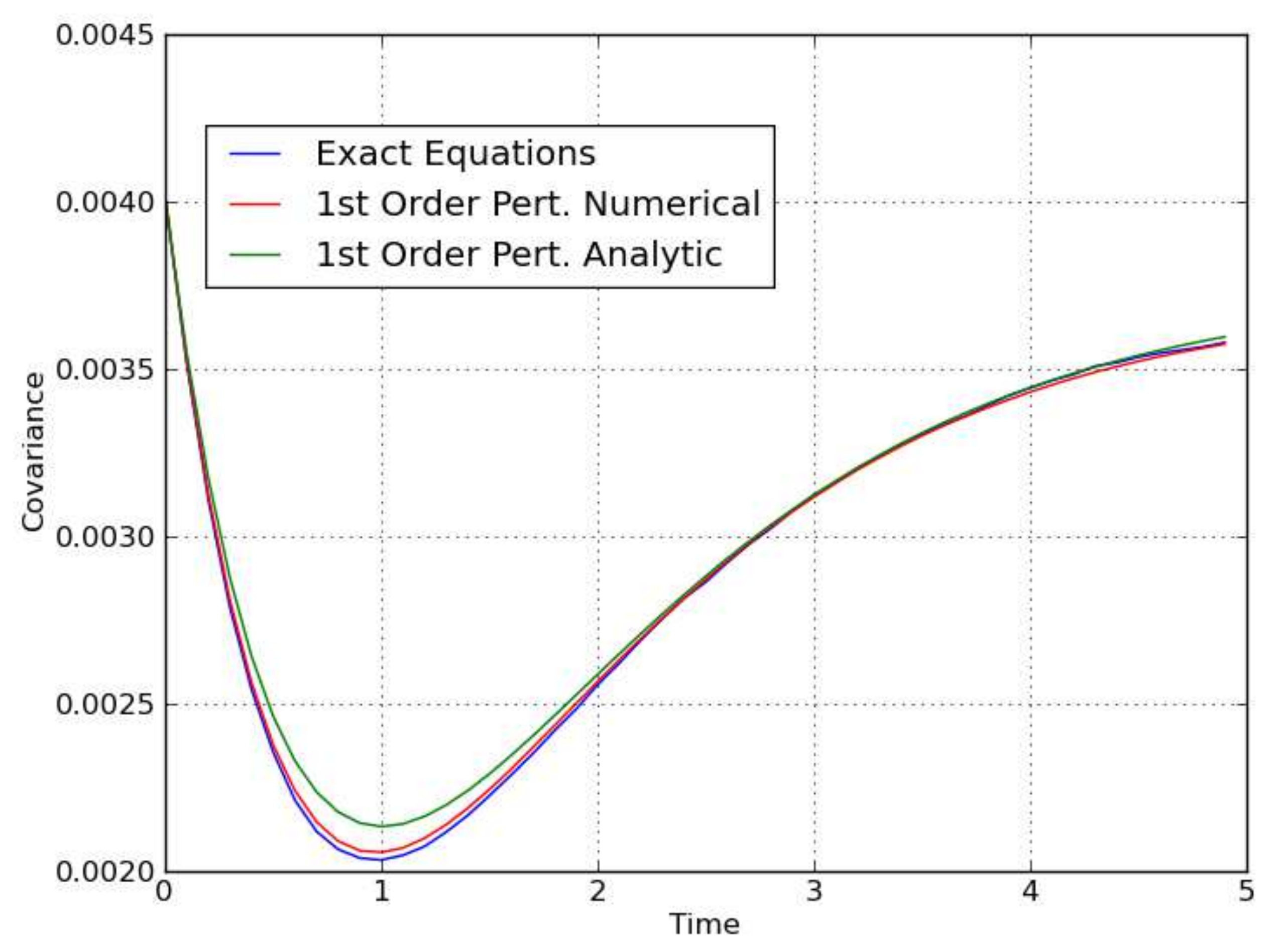}\includegraphics[scale=0.3]{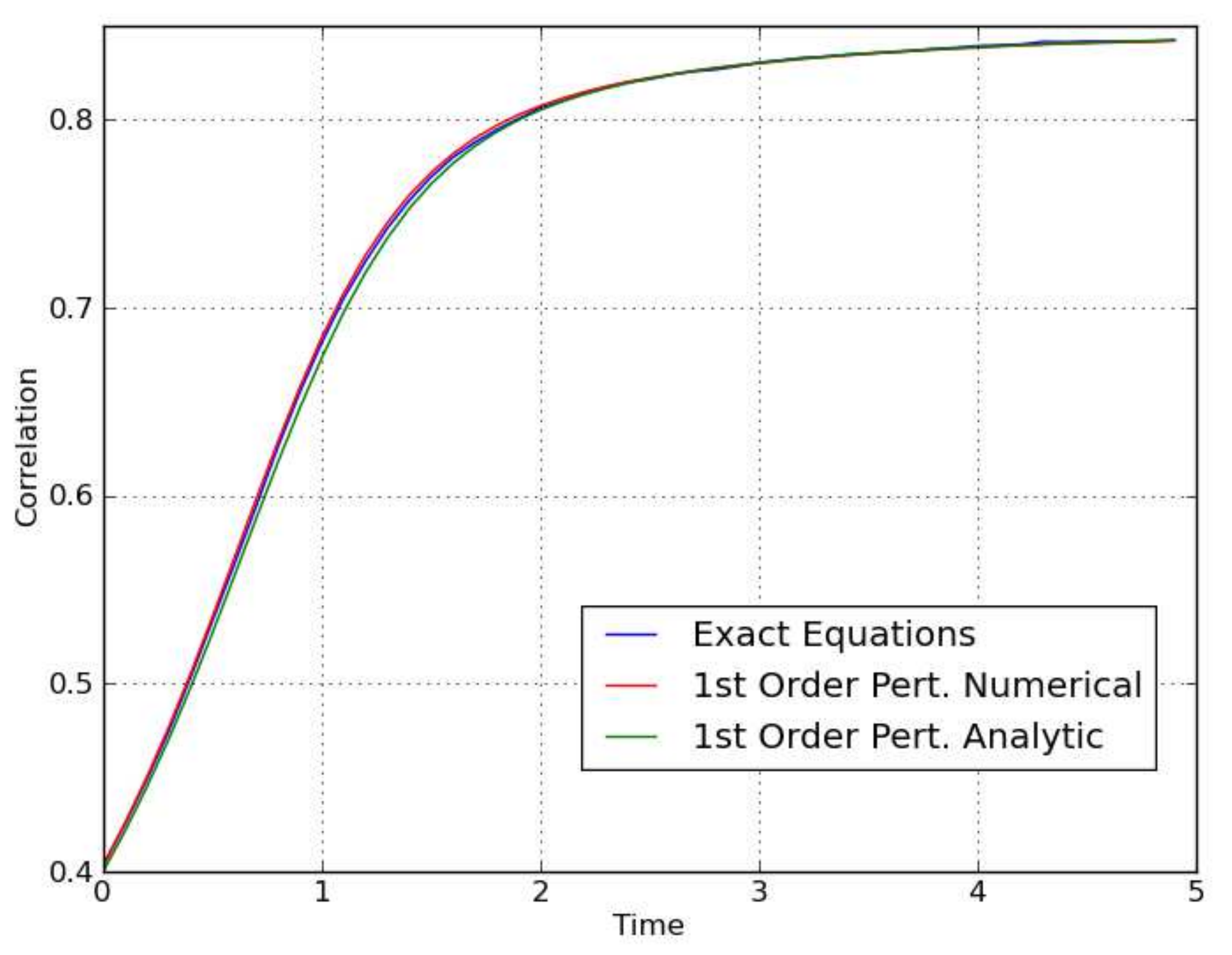}
\par\end{centering}

\caption[{\footnotesize{Numerical comparison of the perturbative expansion
with strong weights - 2}}]{{\footnotesize{\label{sec:circular-ladder-graph-simulation-2}First-order
perturbative expansion for a network with connectivity matrix $CL_{10}$.
These results have been obtained for the values of the parameters
reported in Table \ref{tab:simulation-parameters-1}, for $\sigma_{1}=0.01$,
$\sigma_{2}=\sigma_{3}=0.1$ and with the statistics evaluated through
$10,000$ Monte Carlo simulations. The parameter $\sigma_{1}$ is
small because high values determine large fluctuations of the variance
and covariance (see Figure \ref{fig:circular-ladder-graph-simulation-5}),
so in that case a higher number of Monte Carlo simulations is required
in order to obtain a good match.}}}
\end{figure}

\begin{figure}
\begin{centering}
\includegraphics[scale=0.3]{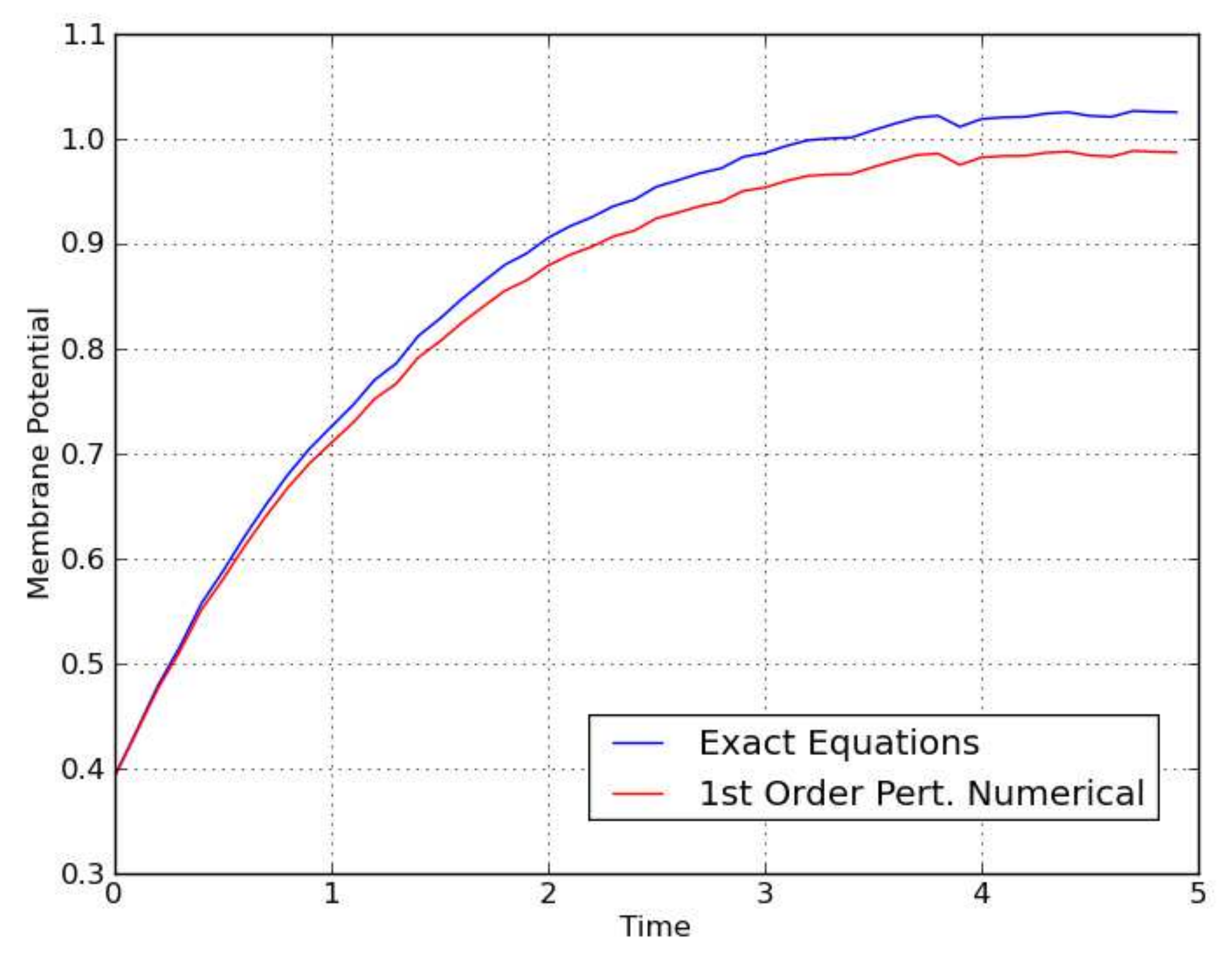}\includegraphics[scale=0.3]{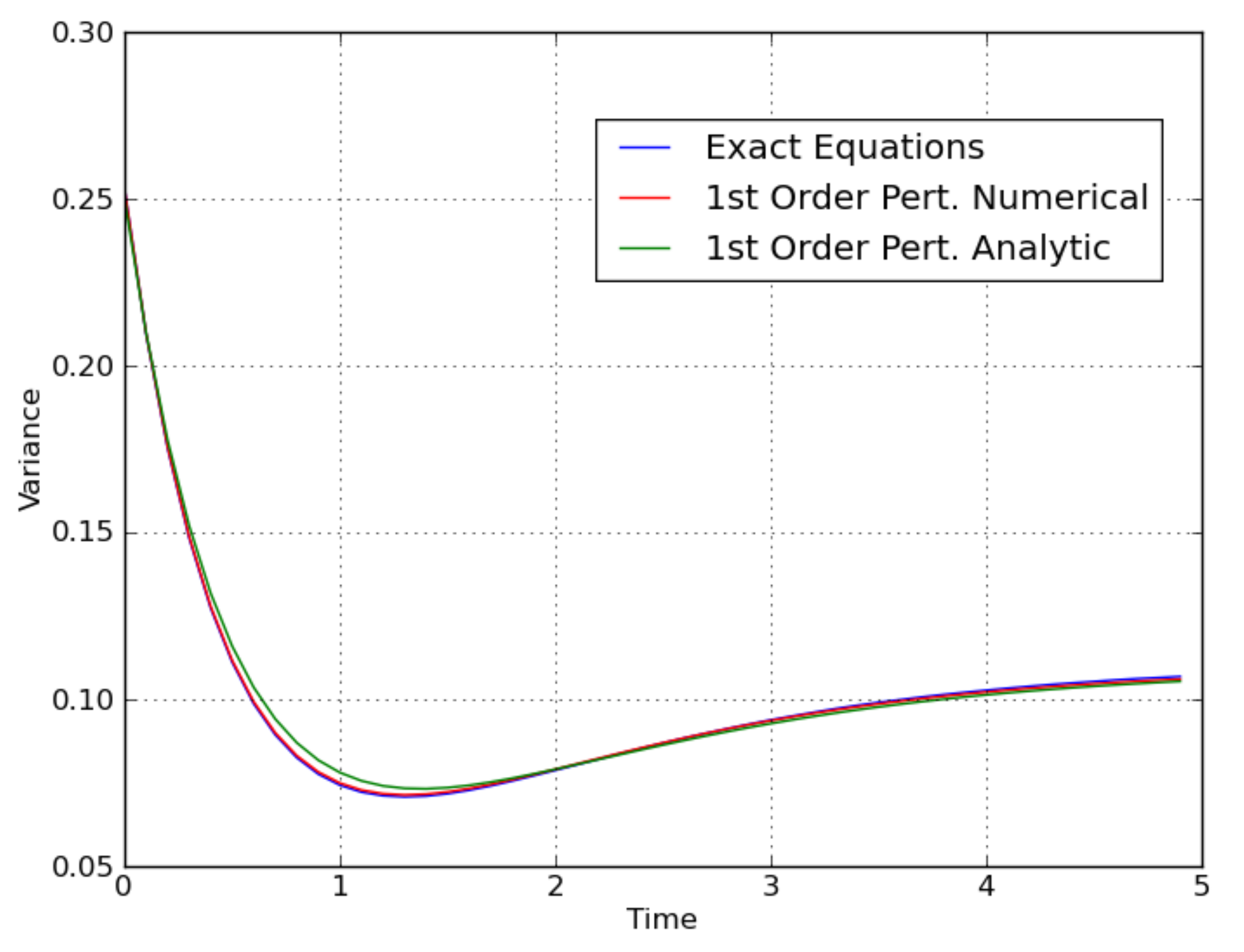}
\par\end{centering}

\begin{centering}
\includegraphics[scale=0.3]{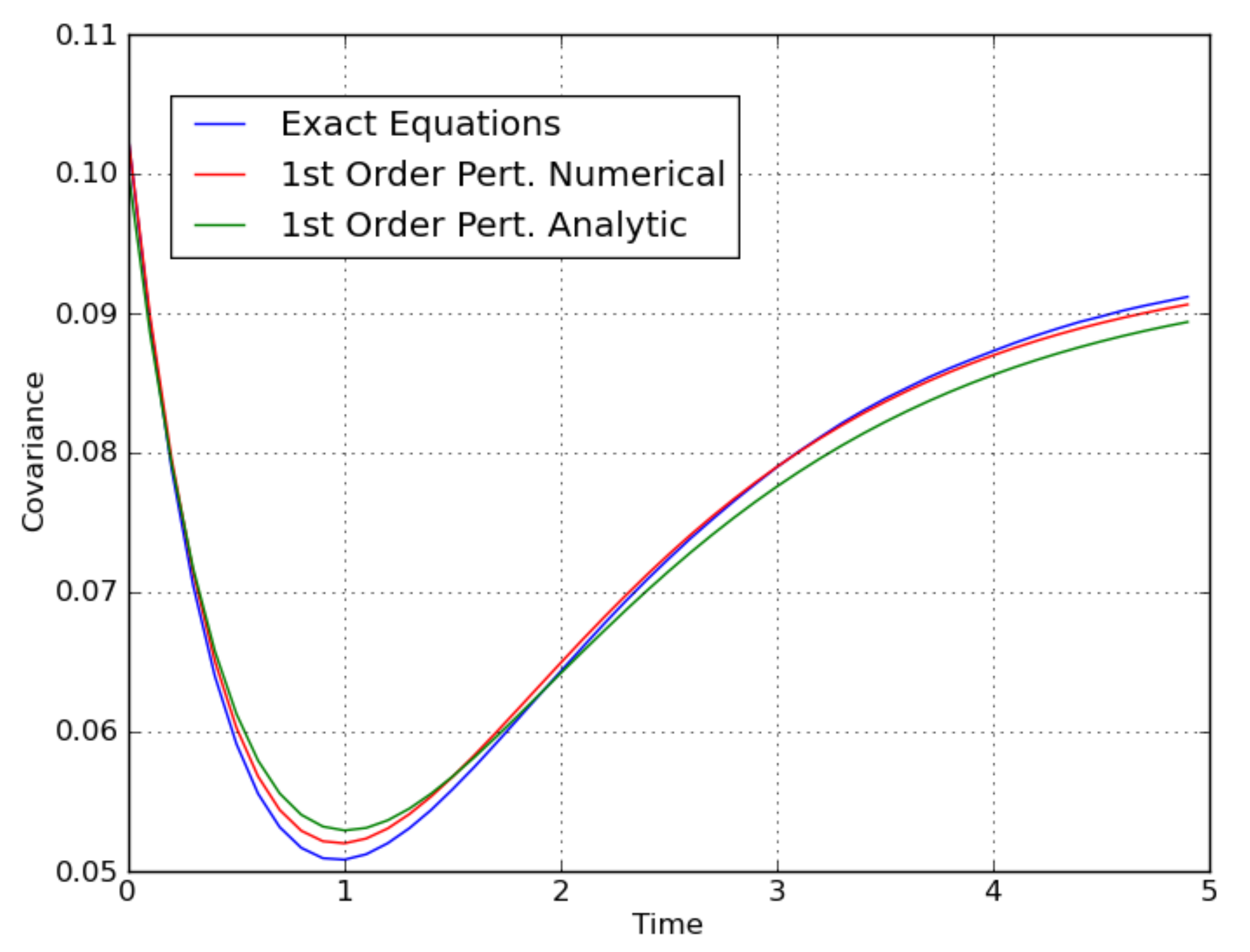}\includegraphics[scale=0.3]{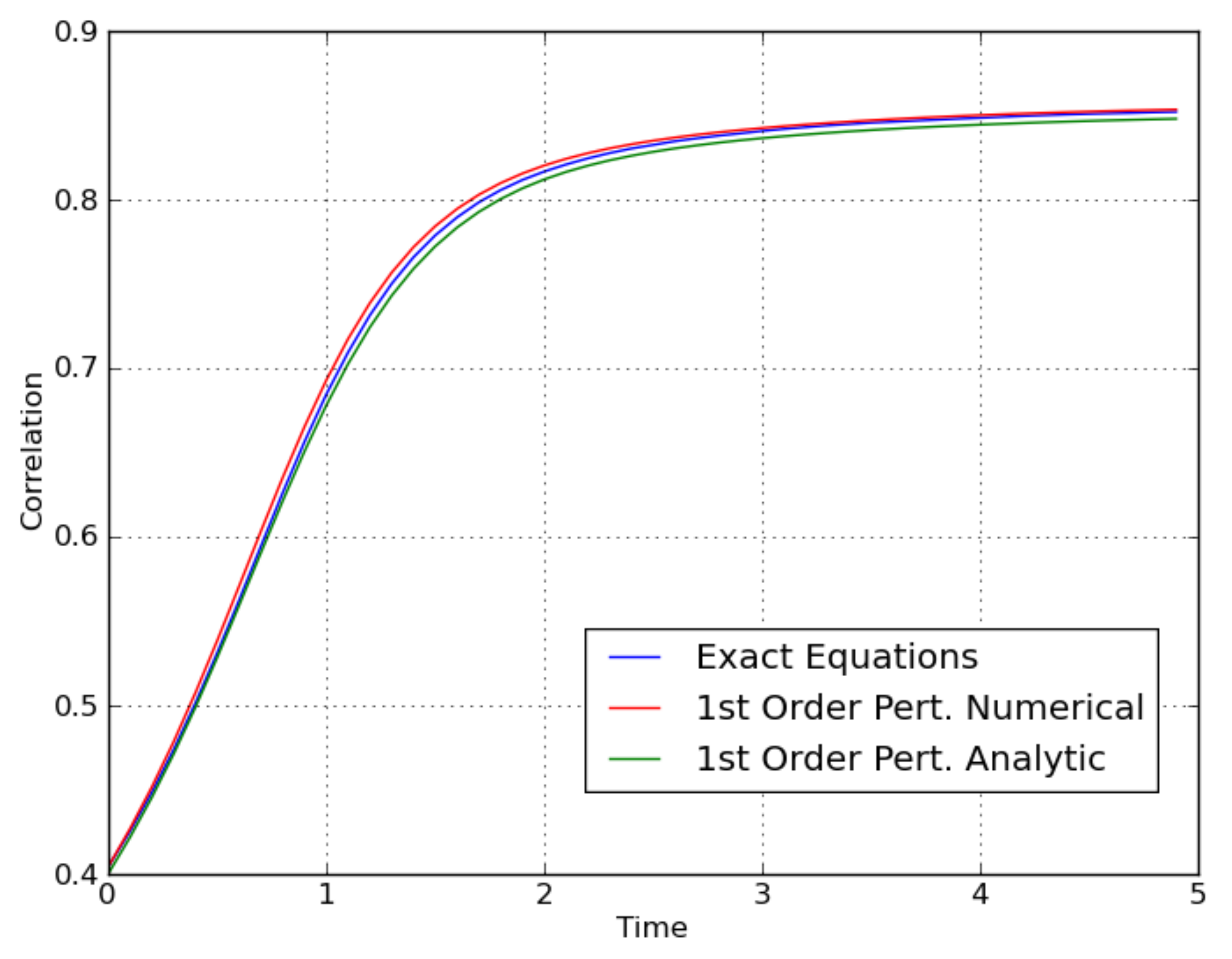}
\par\end{centering}

\caption[{\footnotesize{Numerical comparison of the perturbative expansion
with strong weights - 3}}]{{\footnotesize{\label{fig:circular-ladder-graph-simulation-3}First-order
perturbative expansion for a network with connectivity matrix $CL_{10}$.
These results have been obtained for the values of the parameters
reported in Table \ref{tab:simulation-parameters-1}, for $\sigma_{1}=0.01$,
$\sigma_{2}=\sigma_{3}=0.5$ and with the statistics evaluated through
$10,000$ Monte Carlo simulations.}}}
\end{figure}

\begin{figure}
\begin{centering}
\includegraphics[scale=0.3]{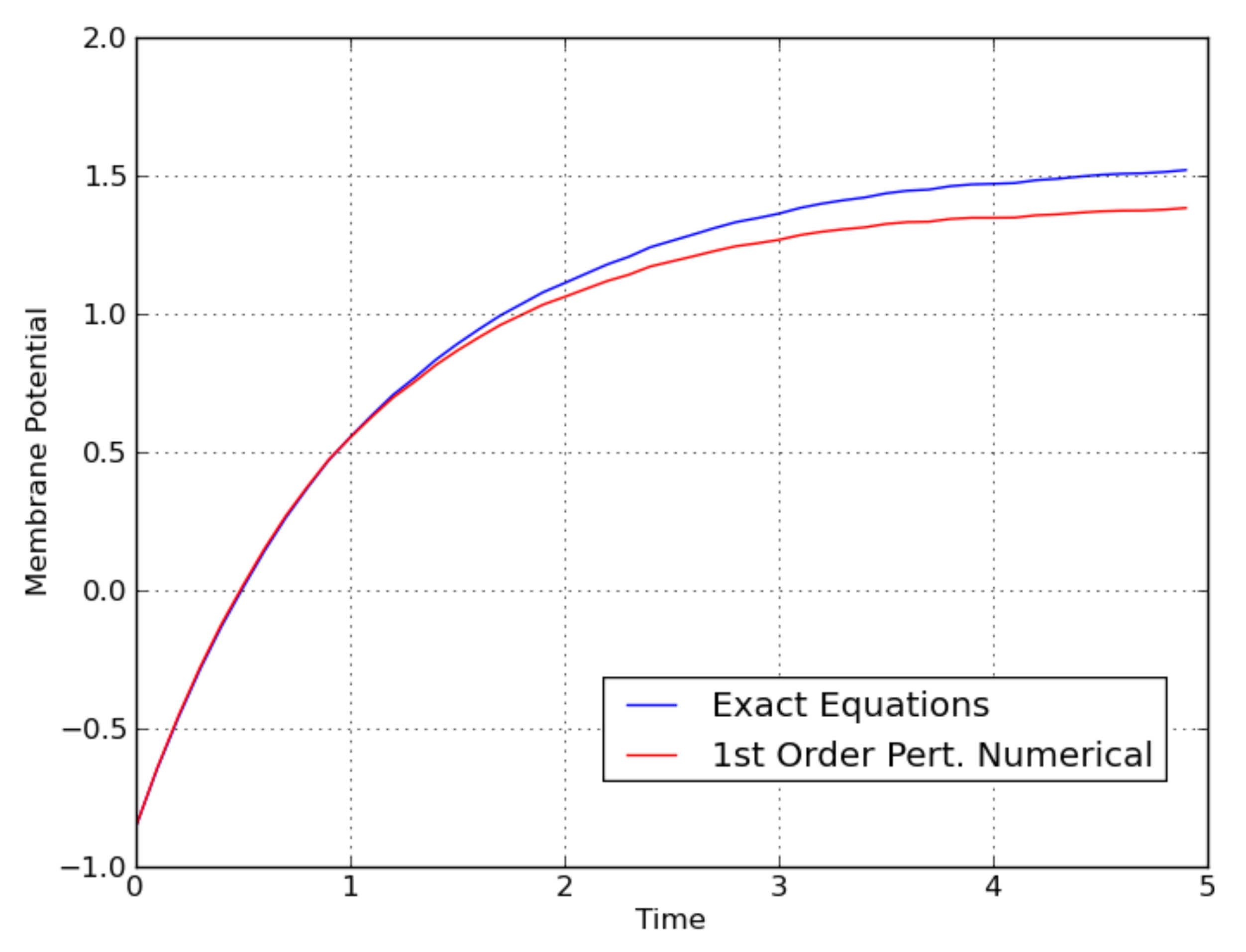}\includegraphics[scale=0.3]{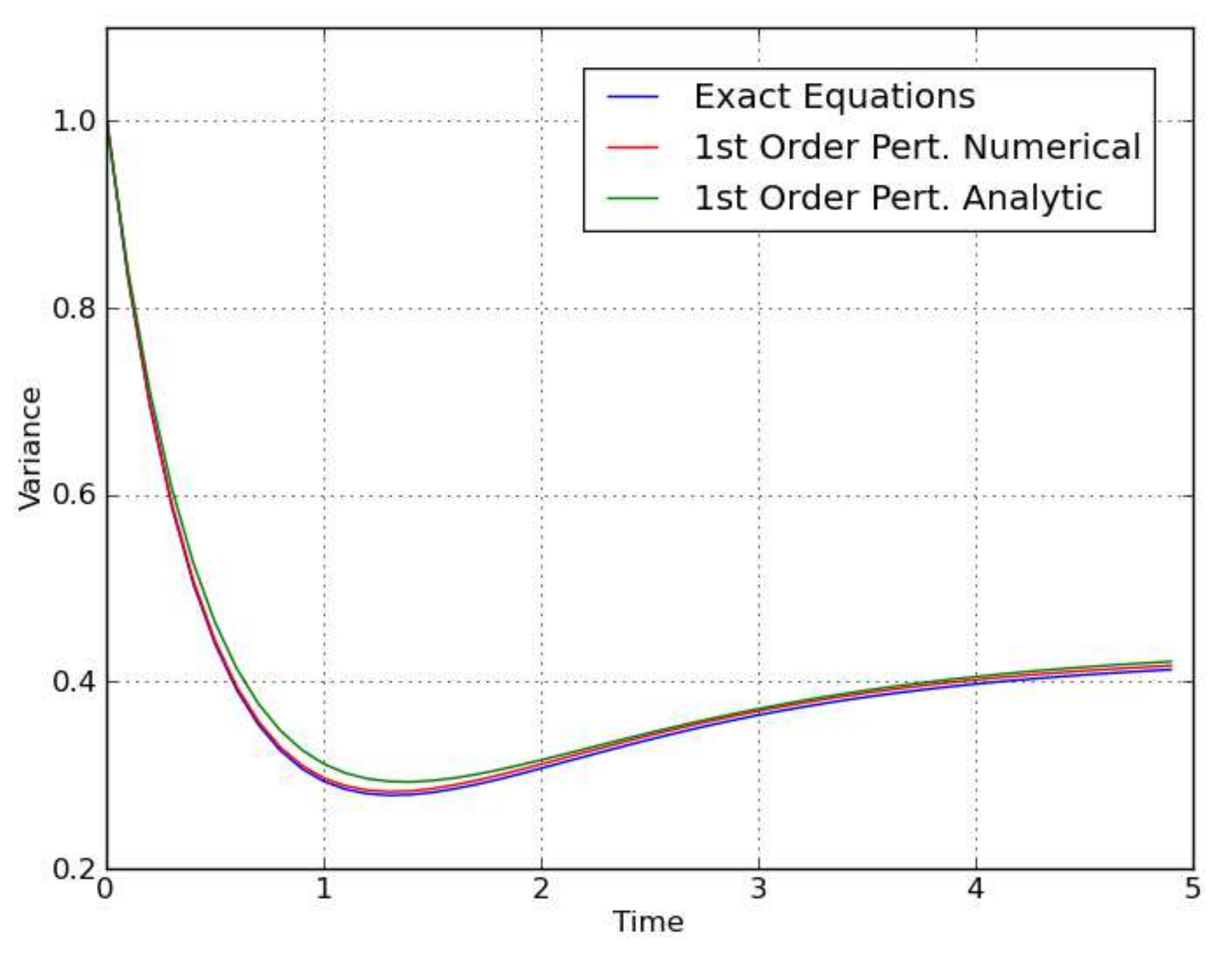}
\par\end{centering}

\begin{centering}
\includegraphics[scale=0.3]{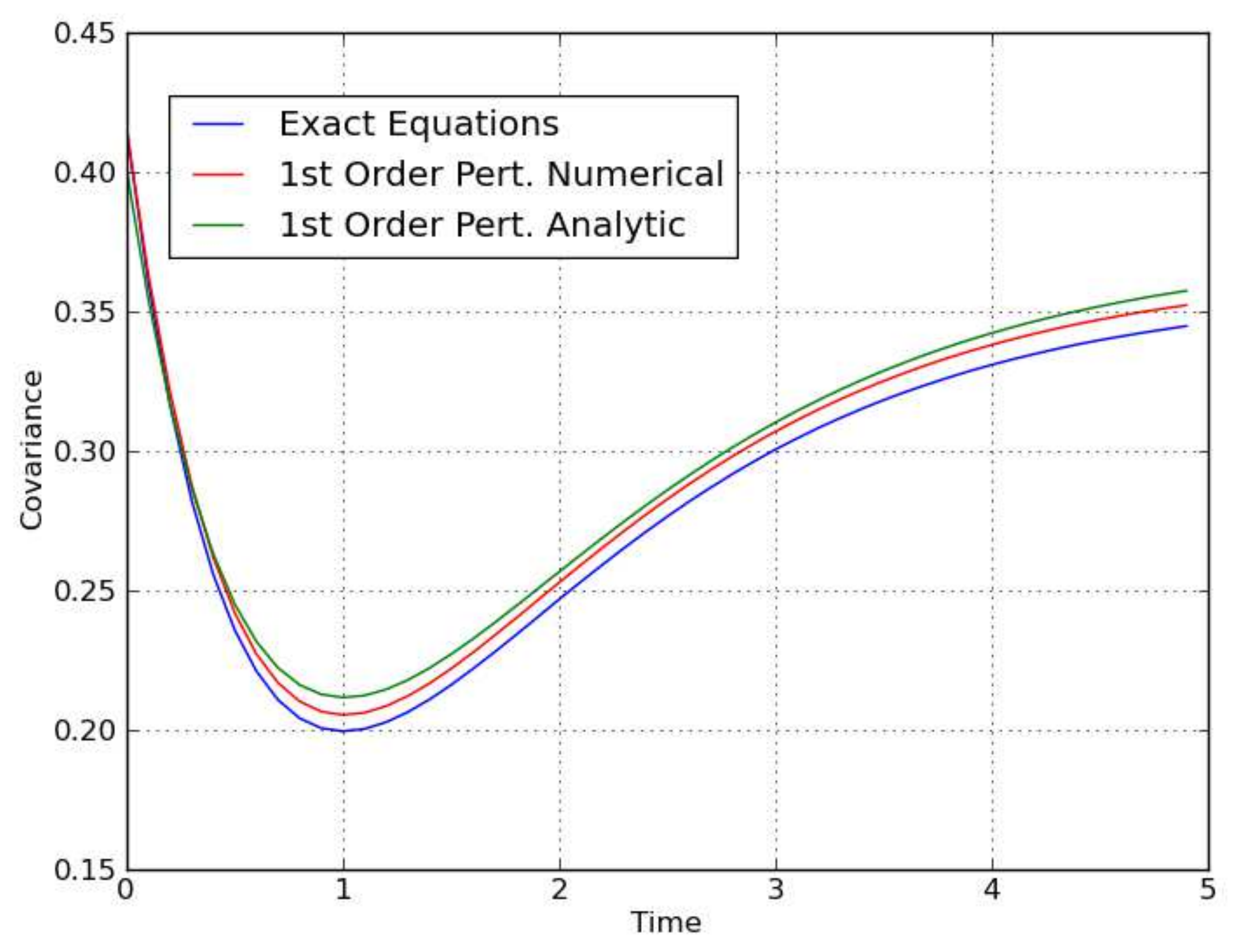}\includegraphics[scale=0.3]{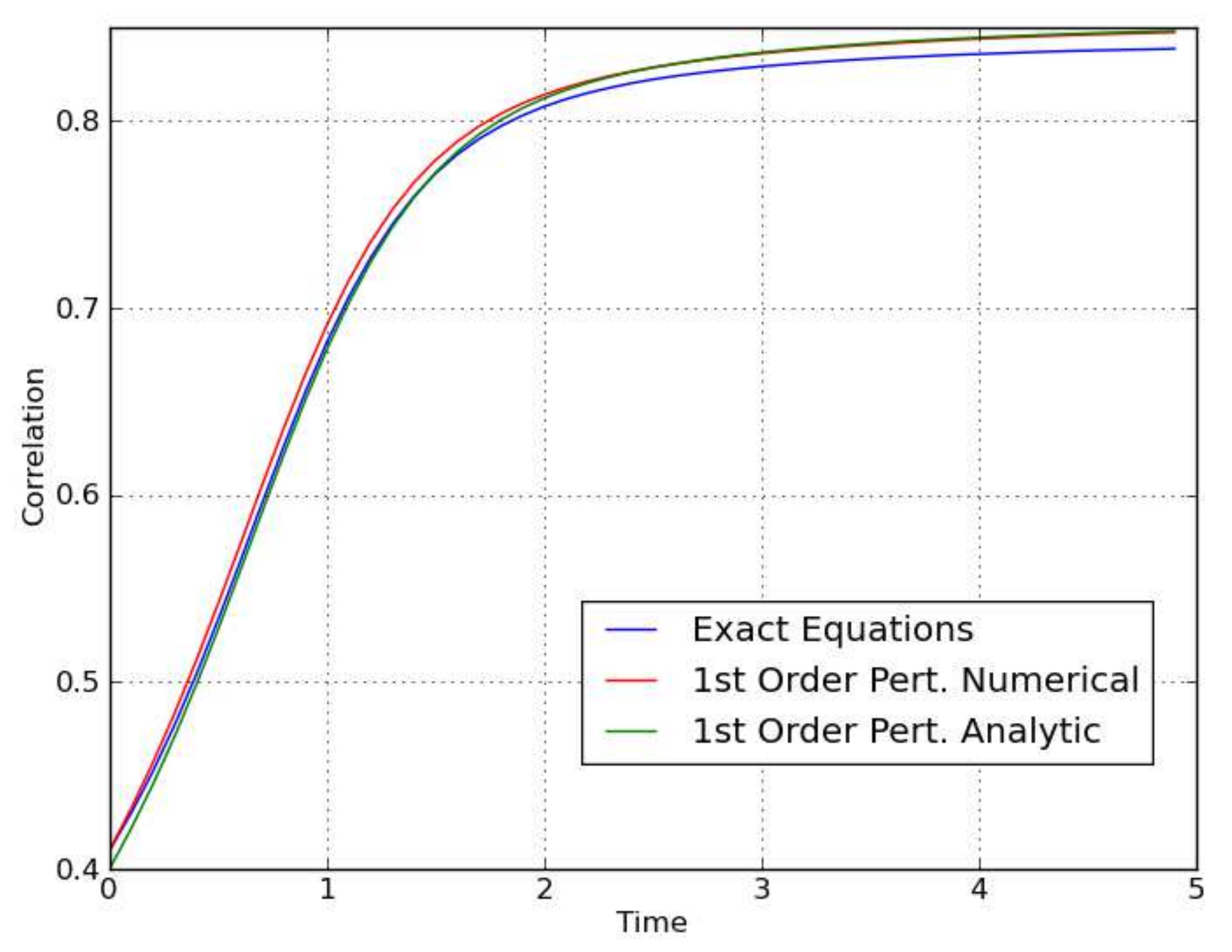}
\par\end{centering}

\caption[{\footnotesize{Numerical comparison of the perturbative expansion
with strong weights - 4}}]{{\footnotesize{\label{fig:circular-ladder-graph-simulation-4}First-order
perturbative expansion for a network with connectivity matrix $CL_{10}$.
These results have been obtained for the values of the parameters
reported in Table \ref{tab:simulation-parameters-1}, for $\sigma_{1}=0.01$,
$\sigma_{2}=\sigma_{3}=1$ and with the statistics evaluated through
$10,000$ Monte Carlo simulations. The match between the exact behavior
and the first-order perturbative expansion is still reasonably good,
even if $\sigma_{2}$ and $\sigma_{3}$ are large.}}}
\end{figure}

\begin{figure}
\begin{centering}
\includegraphics[scale=0.3]{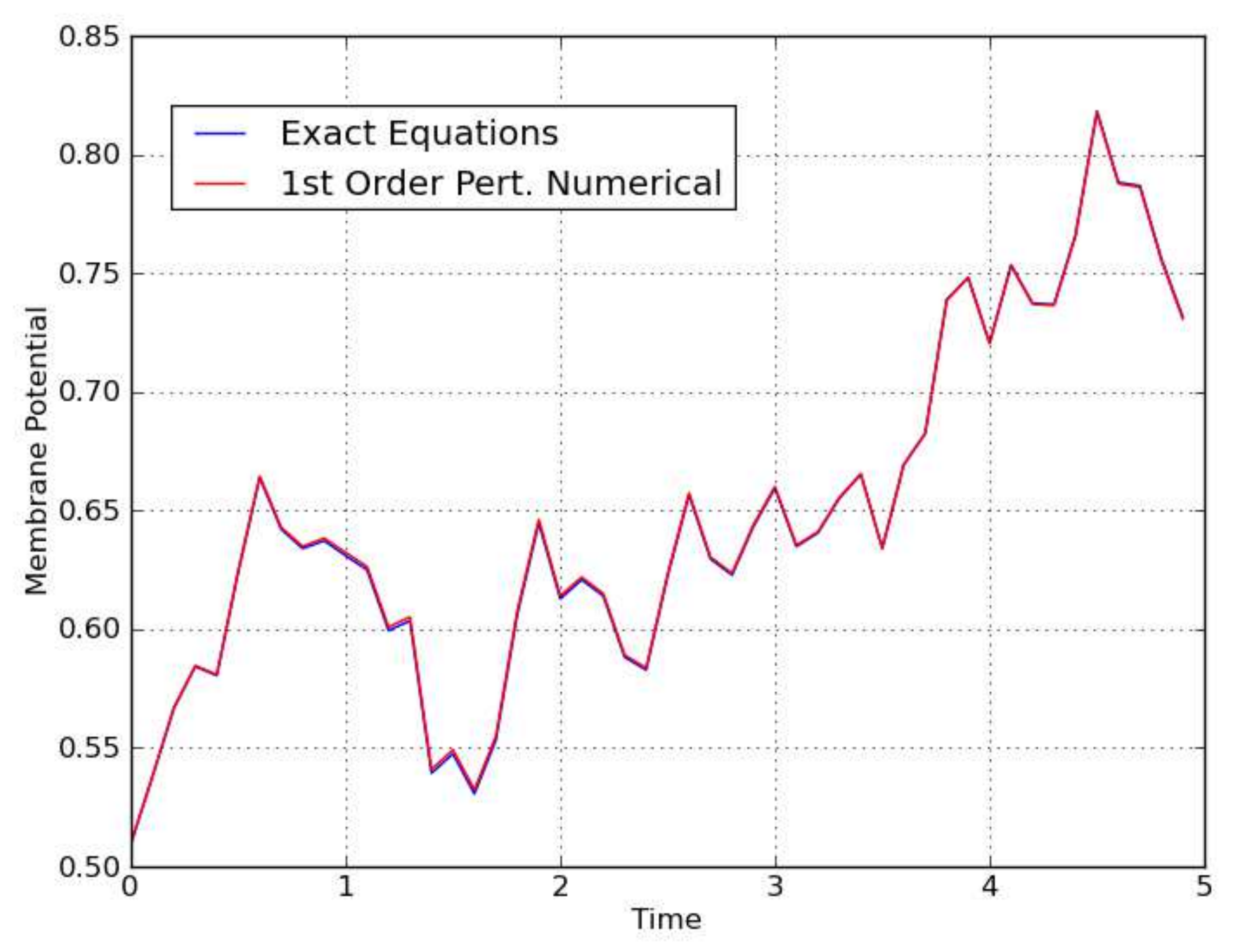}\includegraphics[scale=0.3]{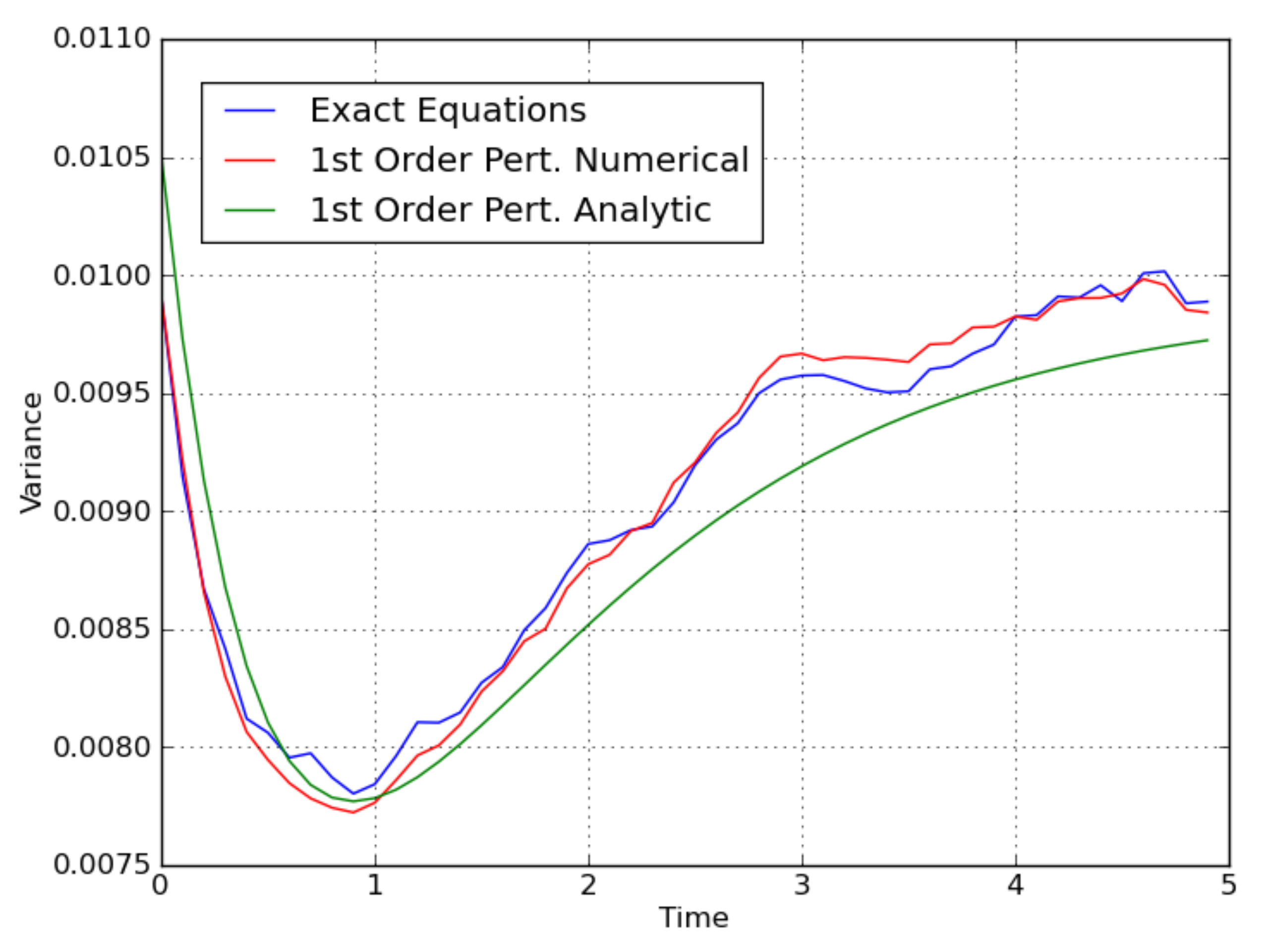}
\par\end{centering}

\begin{centering}
\includegraphics[scale=0.3]{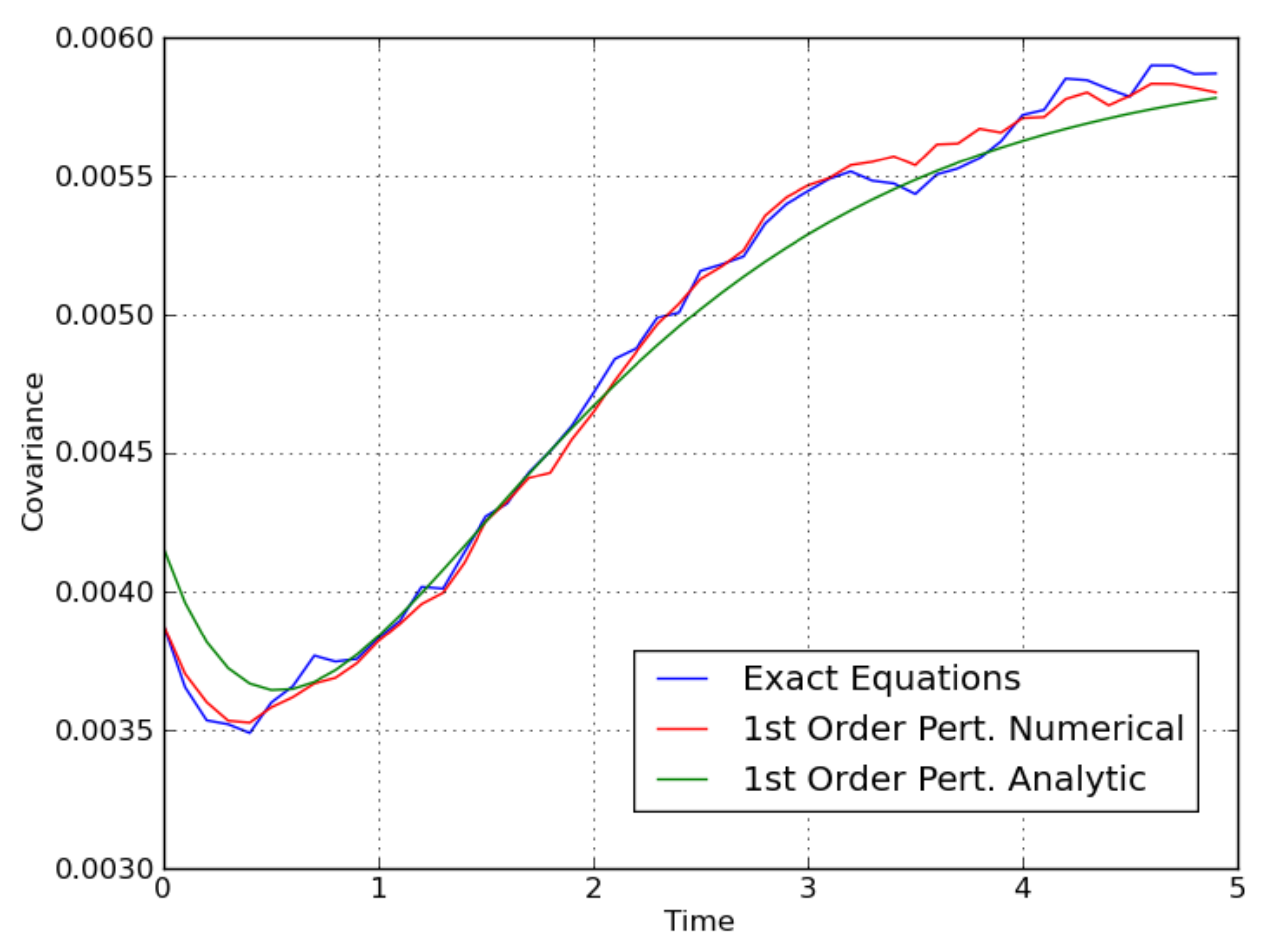}\includegraphics[scale=0.3]{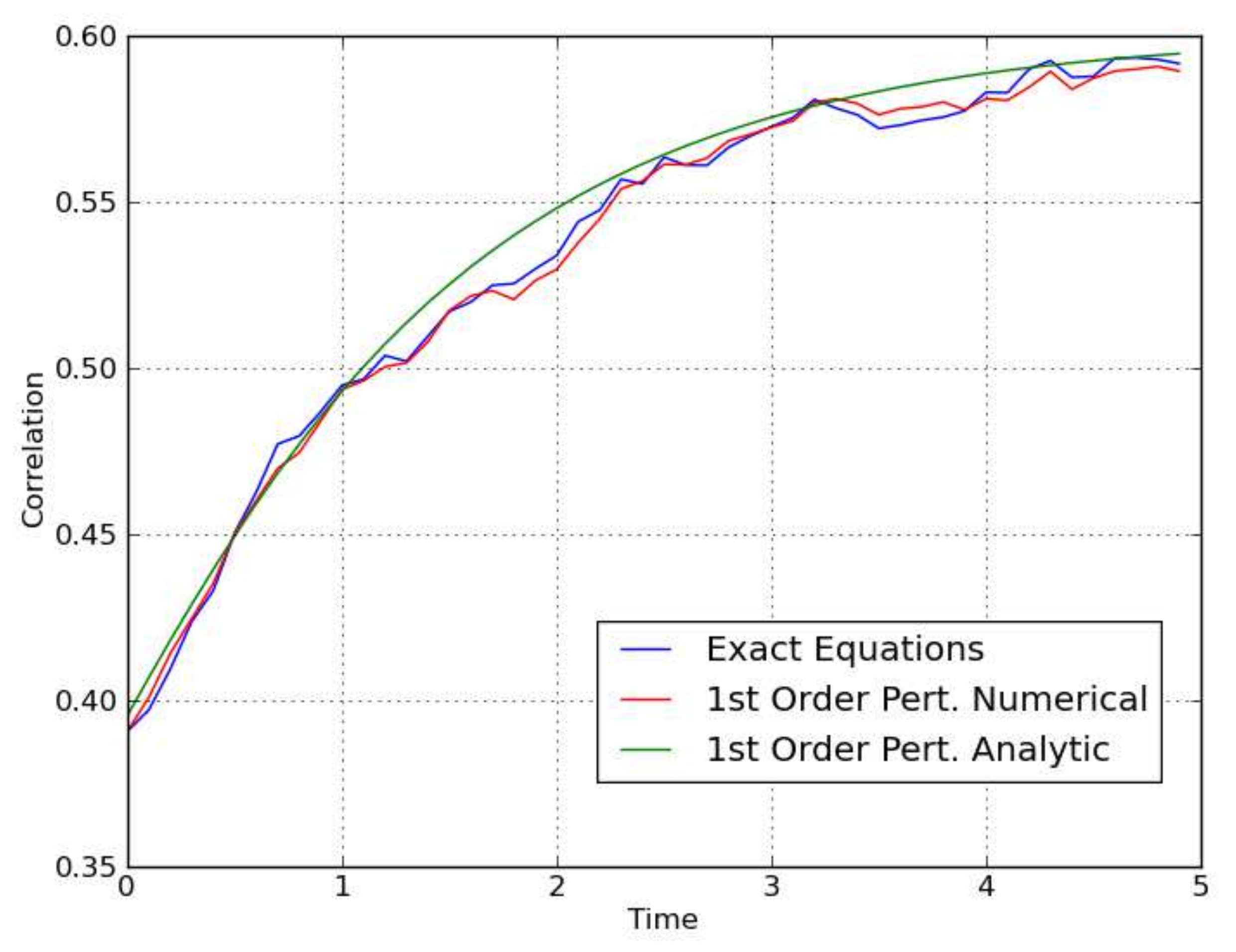}
\par\end{centering}

\caption[{\footnotesize{Numerical comparison of the perturbative expansion
with strong weights - 5}}]{{\footnotesize{\label{fig:circular-ladder-graph-simulation-5}First-order
perturbative expansion for a network with connectivity matrix $CL_{10}$.
These results have been obtained for the values of the parameters
reported in Table \ref{tab:simulation-parameters-1}, for $\sigma_{1}=\sigma_{2}=\sigma_{3}=0.1$
and with the statistics evaluated through $10,000$ Monte Carlo simulations.
The match is not as good as in the previous figures because large
values of $\sigma_{1}$ determine large fluctuations of the variance
and covariance. In other terms, the variance (over many repetitions
of groups made up of $10,000$ Monte Carlo simulations each) of the
variance and covariance is large if $\sigma_{1}$ is big.}}}
\end{figure}

\begin{figure}
\begin{centering}
\includegraphics[scale=0.3]{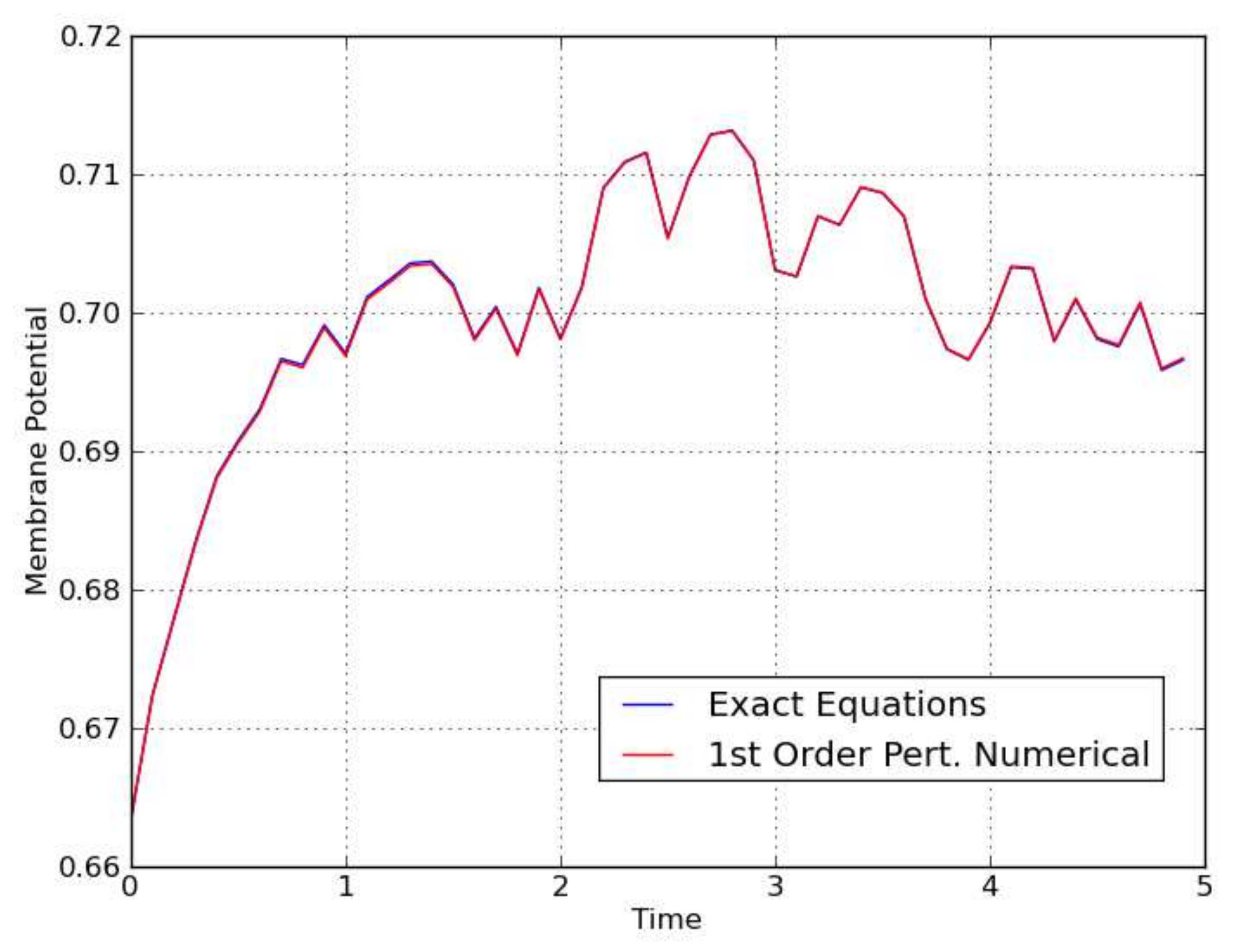}\includegraphics[scale=0.3]{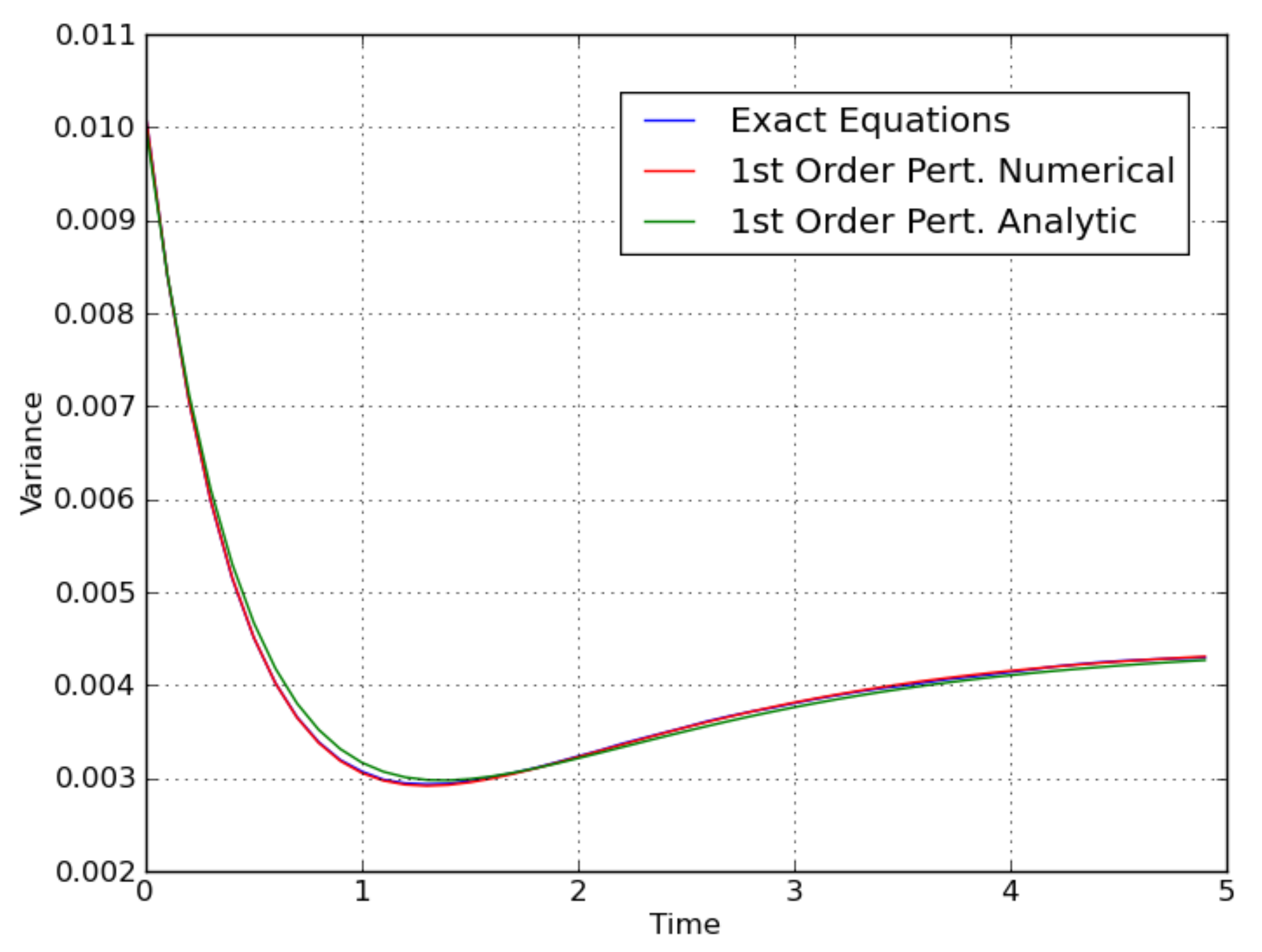}
\par\end{centering}

\begin{centering}
\includegraphics[scale=0.3]{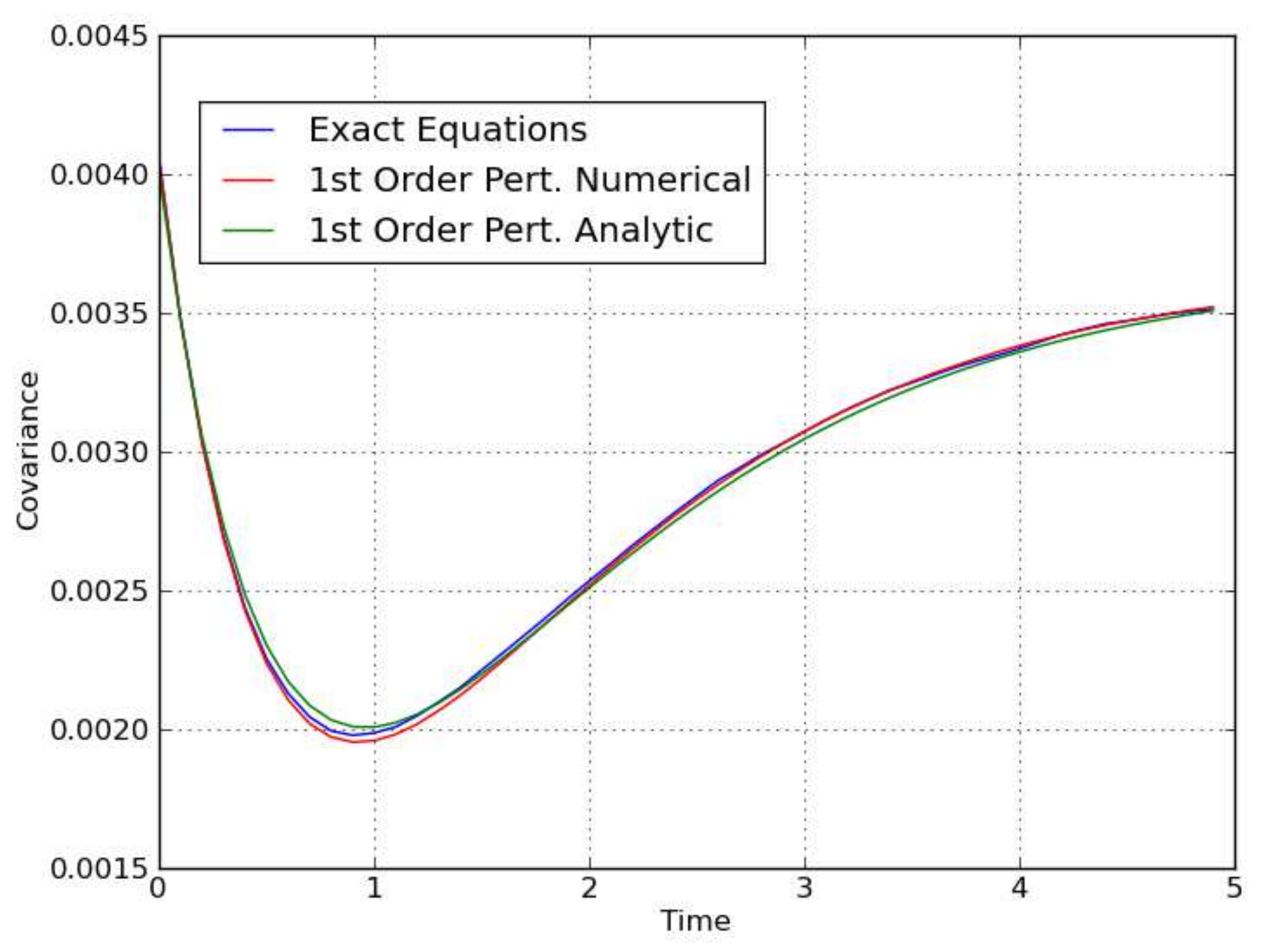}\includegraphics[scale=0.3]{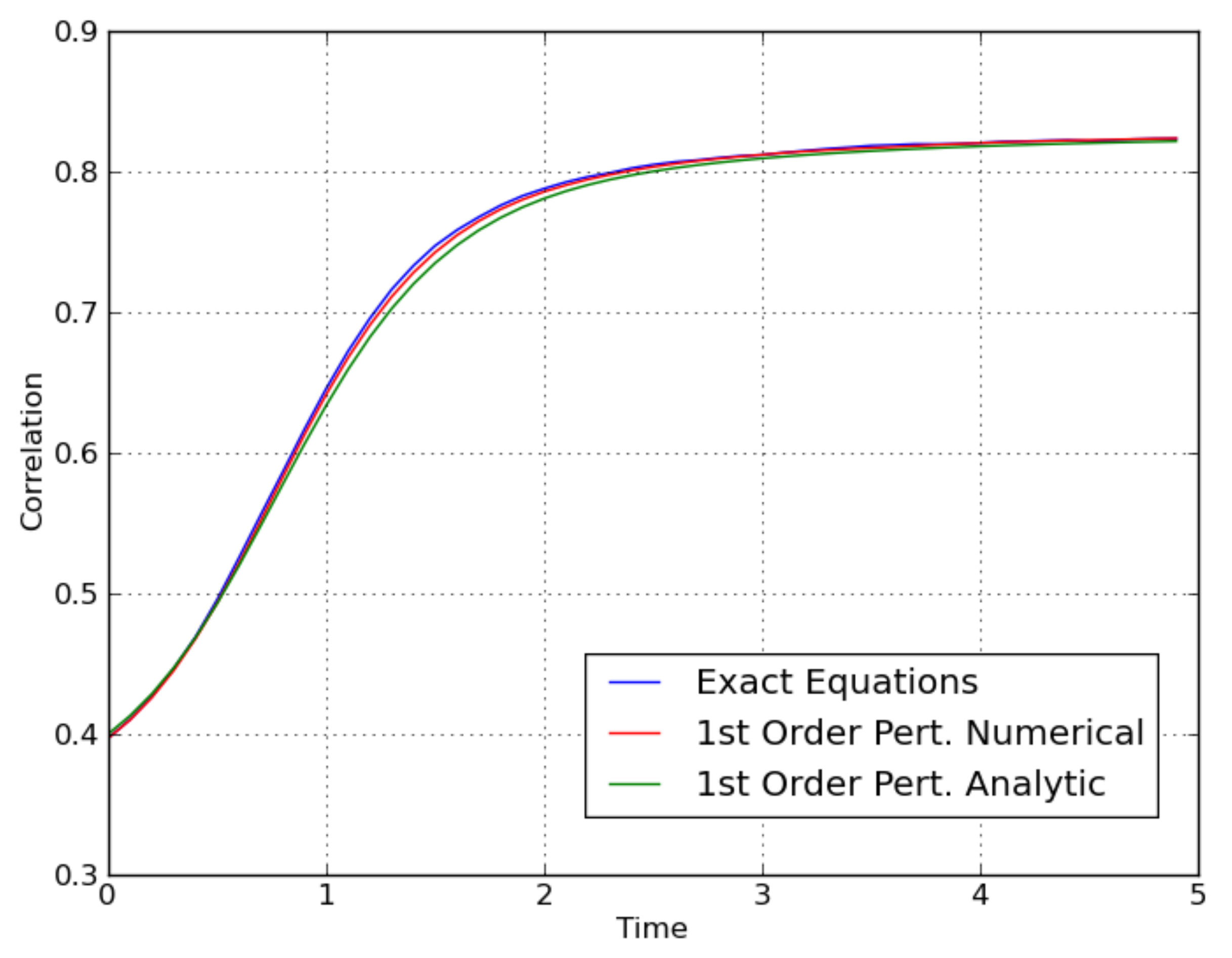}
\par\end{centering}

\caption[{\footnotesize{Numerical comparison of the perturbative expansion
with strong weights - 6}}]{{\footnotesize{\label{fig:cube-graph-simulation}First-order perturbative
expansion for a network with connectivity matrix $H_{3}$ (namely
$N=8$). These results have been obtained for the values of the parameters
reported in Table \ref{tab:simulation-parameters-1}, for $\sigma_{1}=0.01$,
$\sigma_{2}=\sigma_{3}=0.1$ and with the statistics evaluated through
$10,000$ Monte Carlo simulations.}}}
\end{figure}

\begin{figure}
\begin{centering}
\includegraphics[scale=0.3]{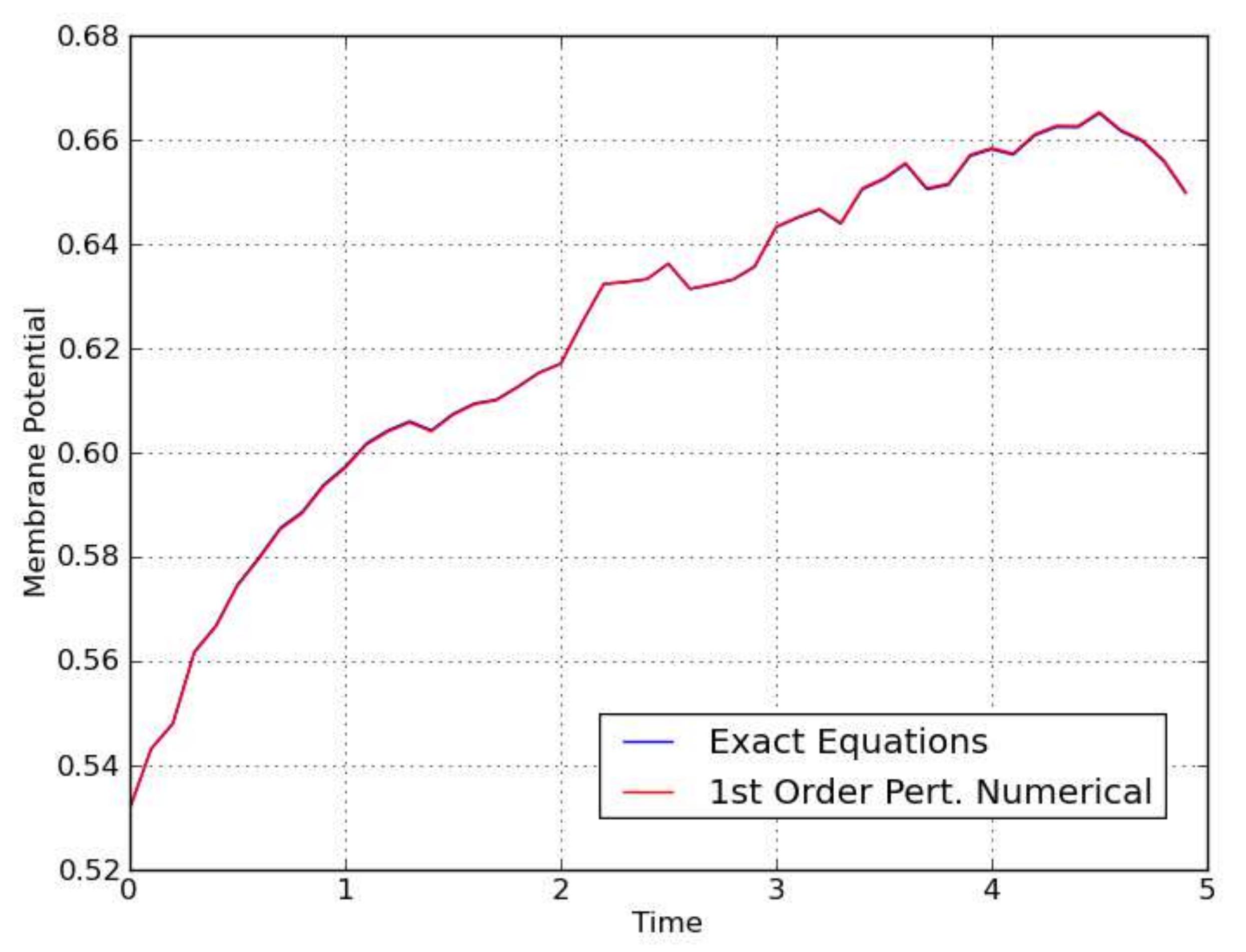}\includegraphics[scale=0.3]{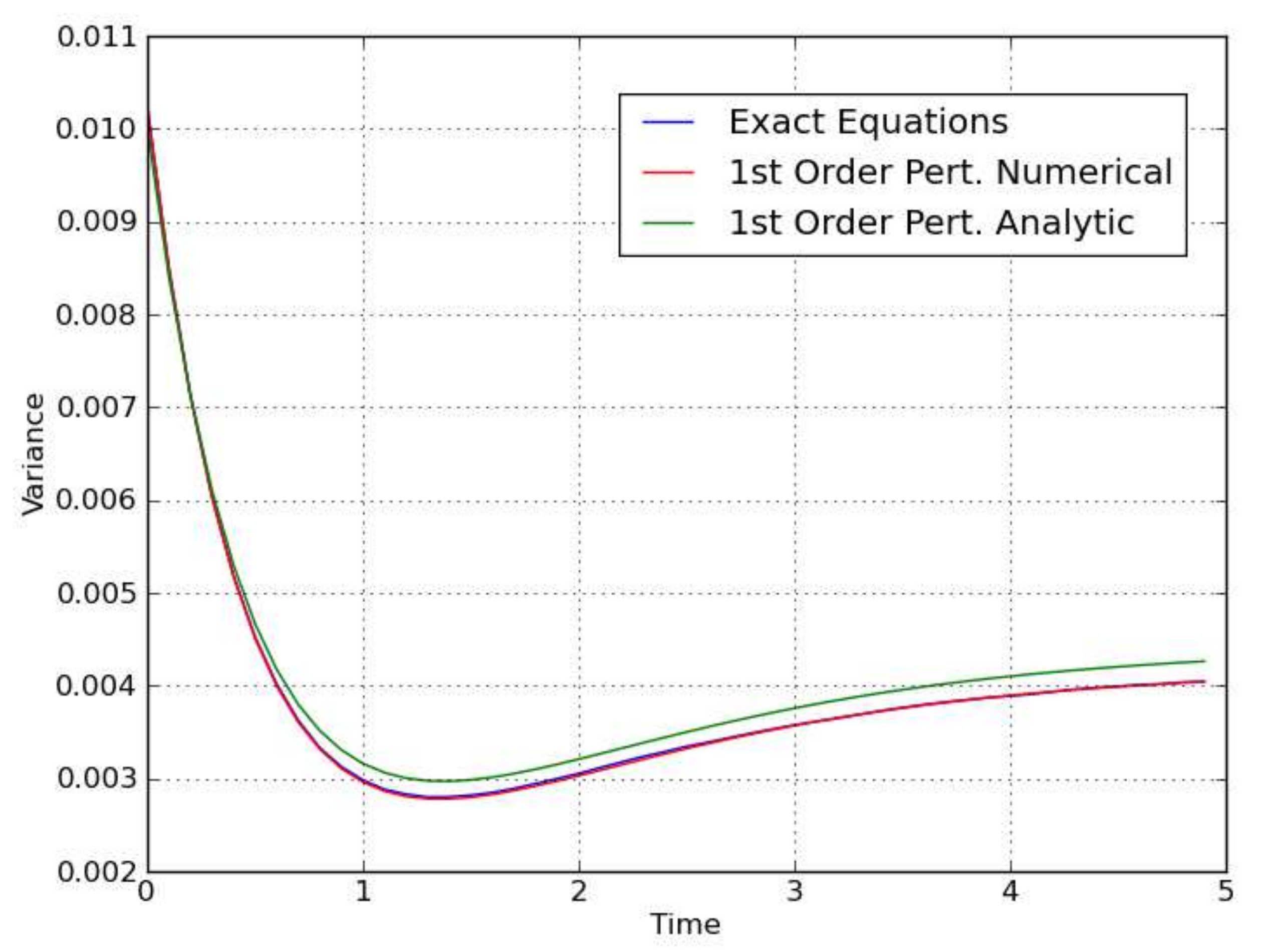}
\par\end{centering}

\begin{centering}
\includegraphics[scale=0.3]{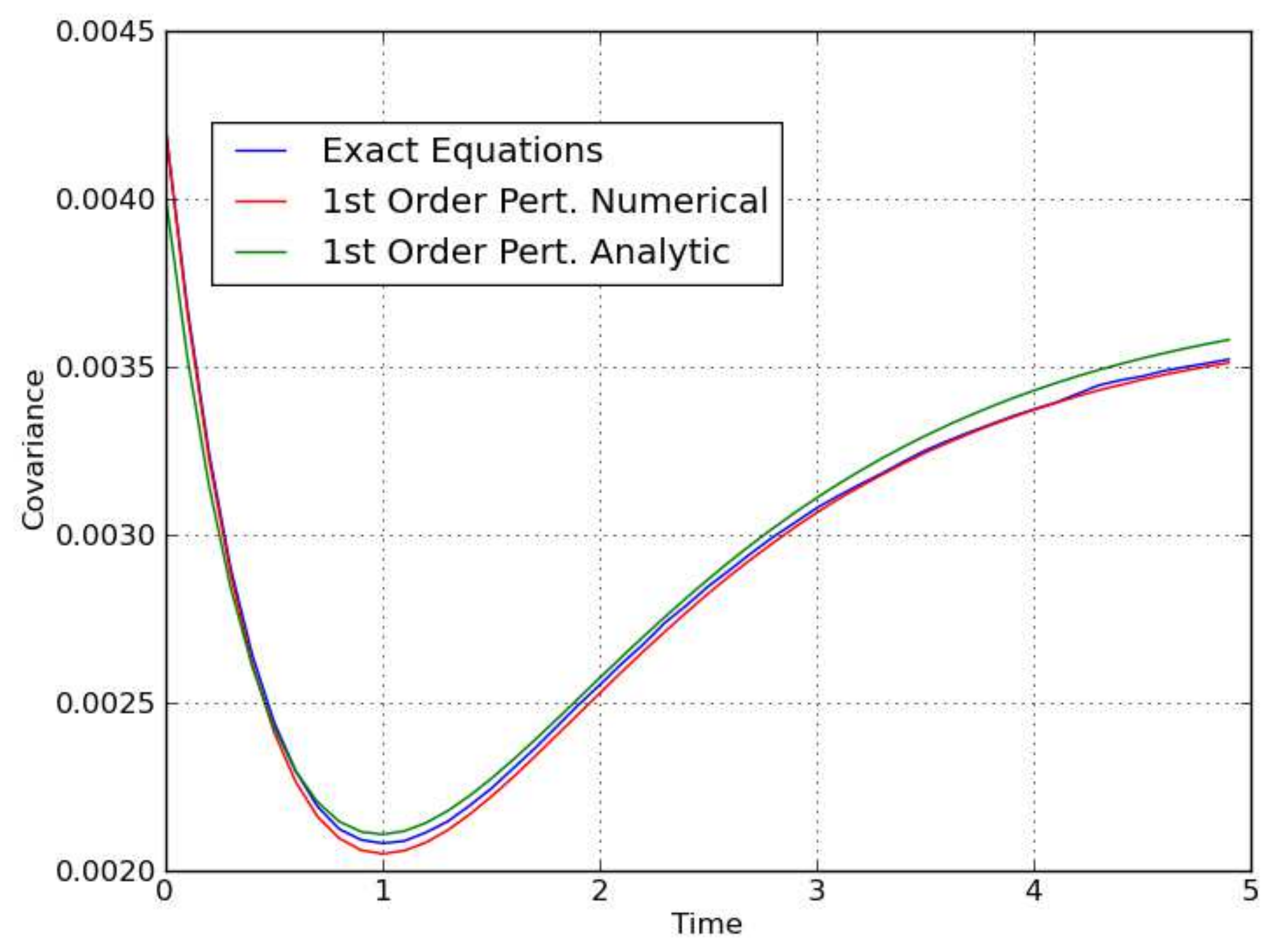}\includegraphics[scale=0.3]{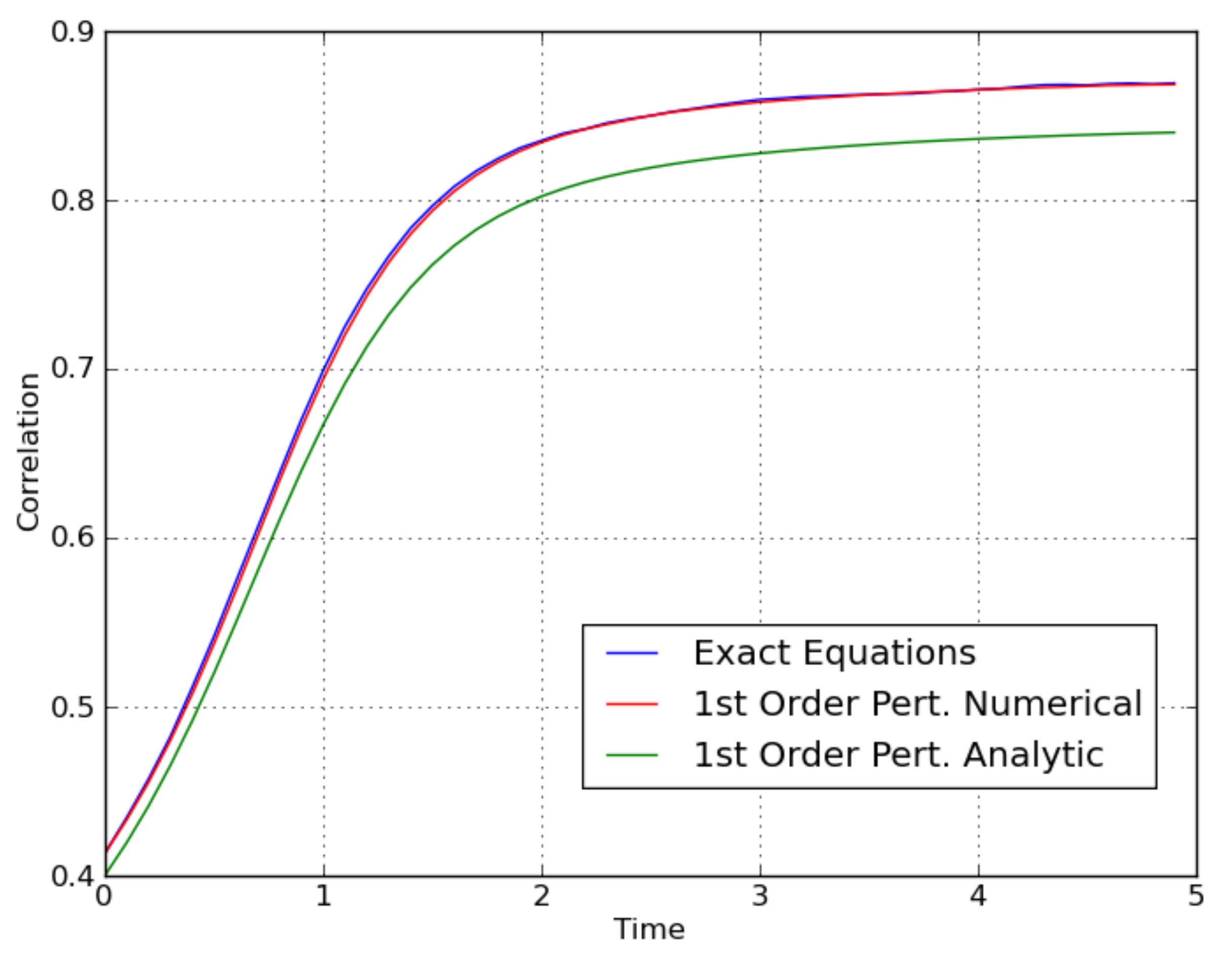}
\par\end{centering}

\caption[{\footnotesize{Numerical comparison of the perturbative expansion
with strong weights - 7}}]{{\footnotesize{}}\label{fig:circulant-graph-simulation}{\footnotesize{First-order
perturbative expansion for a network with connectivity matrix $C_{10}\left(1,2,0,...,0\right)$.
These results have been obtained for the values of the parameters
reported in Table \ref{tab:simulation-parameters-1}, for $\sigma_{1}=0.01$,
$\sigma_{2}=\sigma_{3}=0.1$ and with the statistics evaluated through
$10,000$ Monte Carlo simulations. This figure clearly shows that
the goodness of the match between the curves depends on the connectivity
matrix of the network, for a fixed number of Monte Carlo simulations.}}}
\end{figure}

\begin{figure}
\begin{centering}
\includegraphics[scale=0.3]{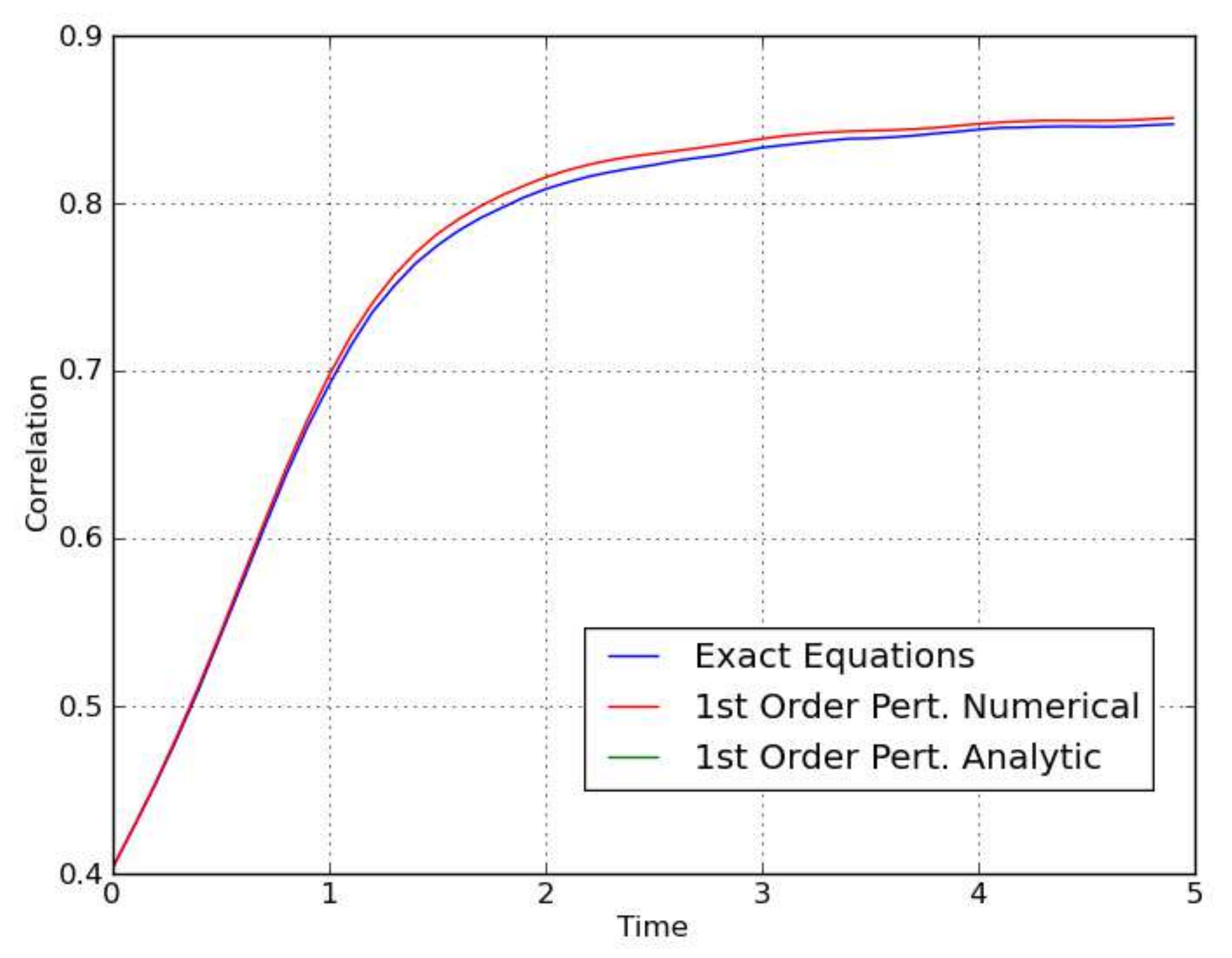}\includegraphics[scale=0.3]{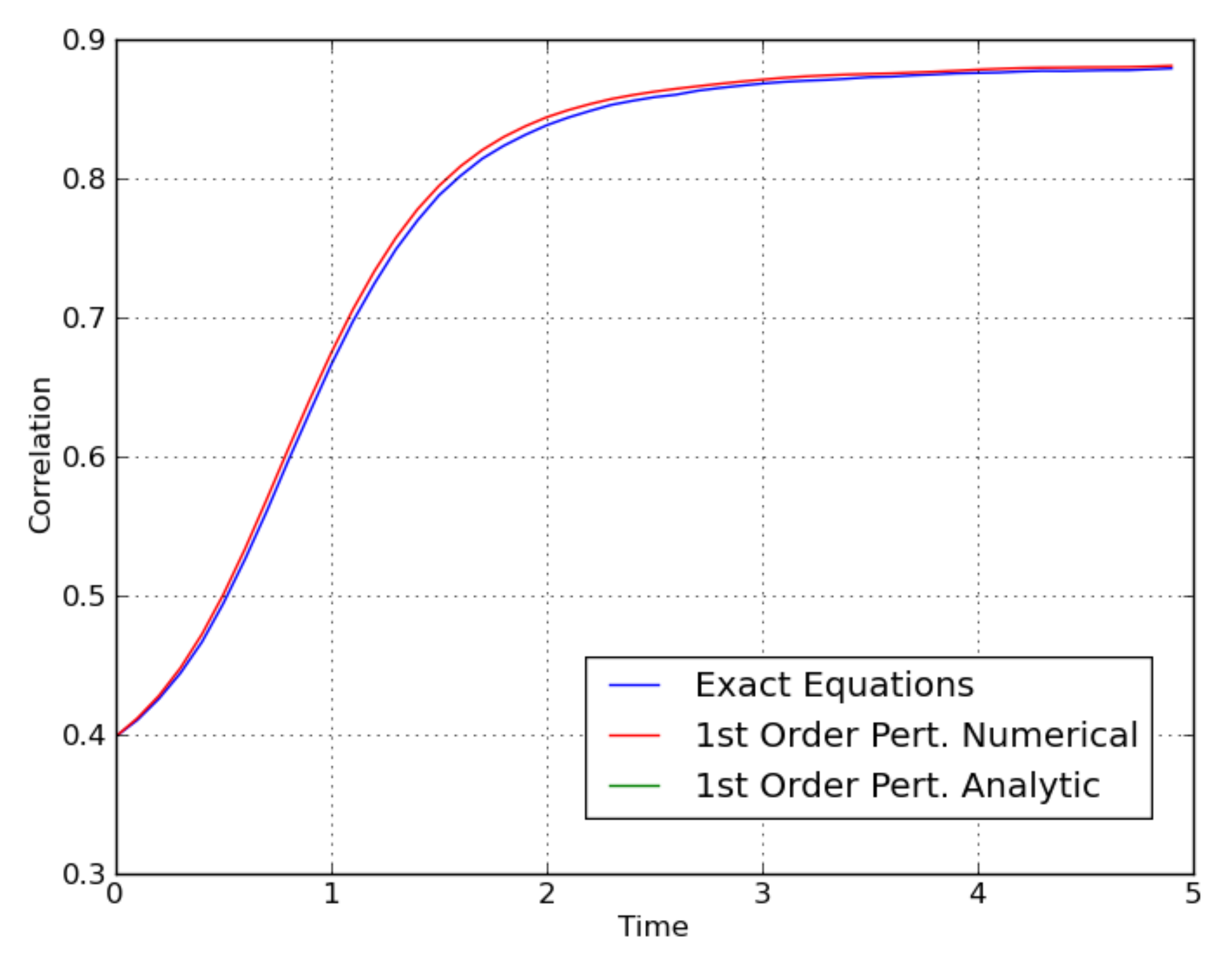}
\par\end{centering}

\begin{centering}
\includegraphics[scale=0.3]{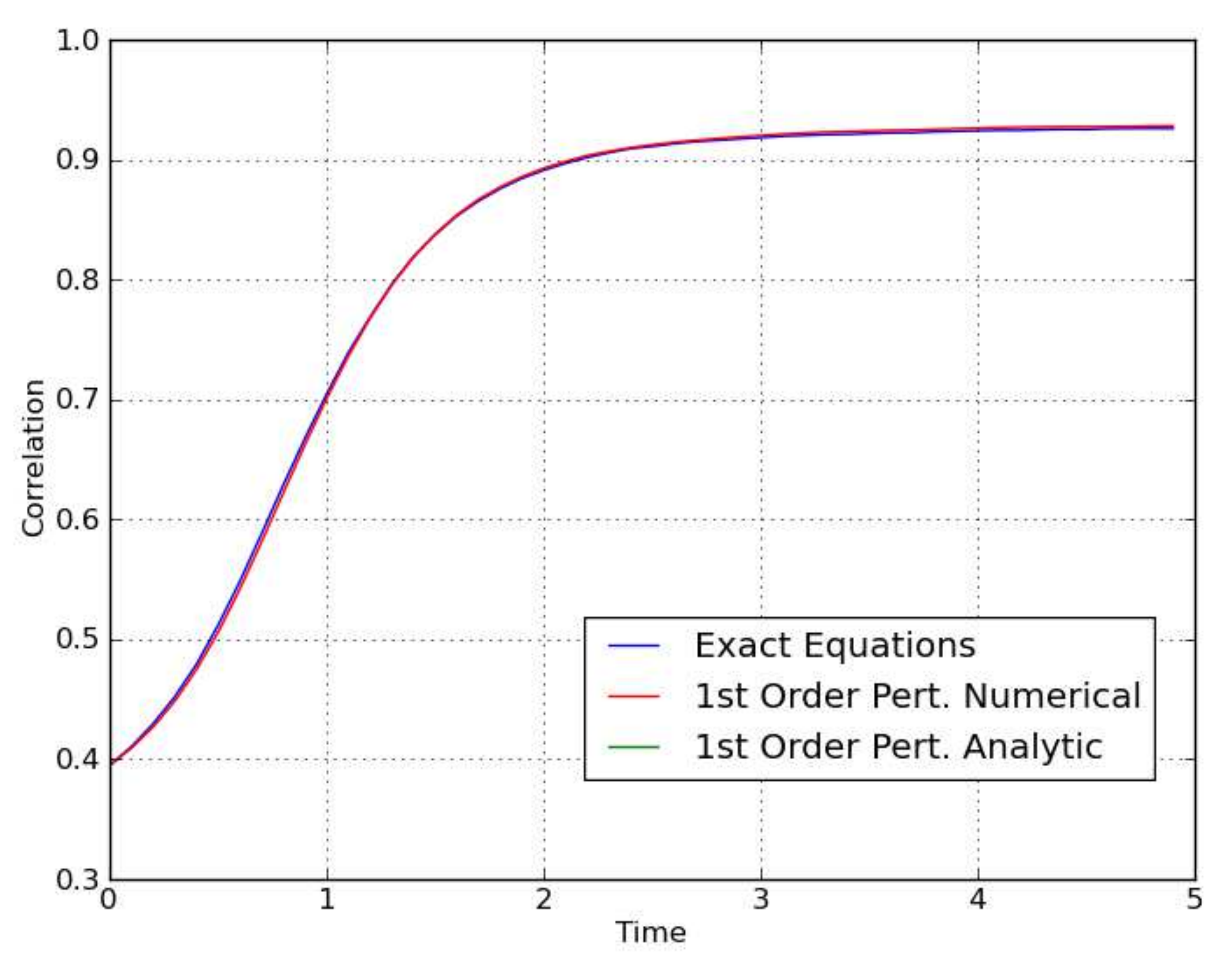}\includegraphics[scale=0.3]{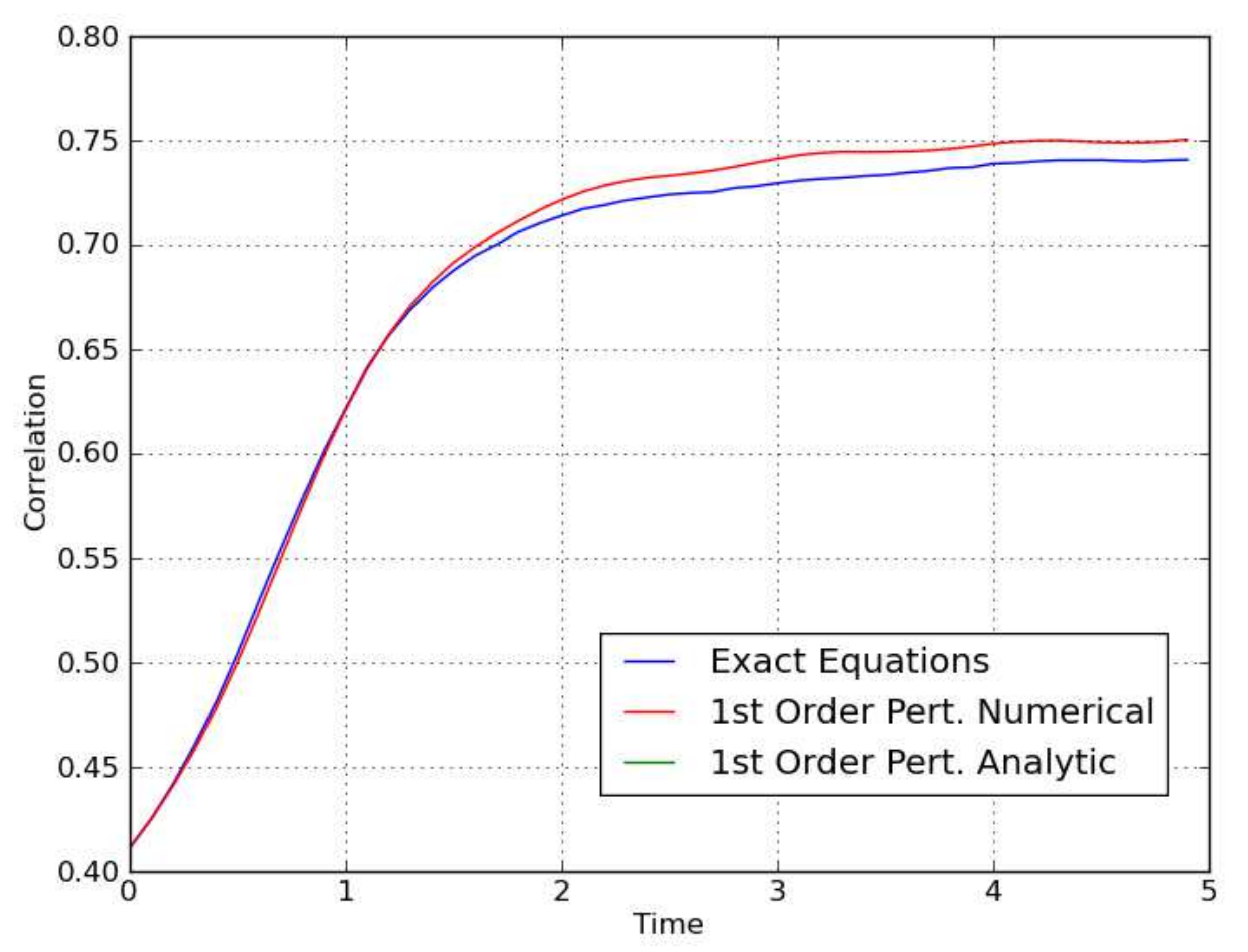}
\par\end{centering}

\caption[{\footnotesize{Numerical comparison of the perturbative expansion
with strong weights - 8}}]{{\footnotesize{\label{fig:second-order-simulation}Correlation function
obtained with the second-order perturbative expansion for a network
with connectivity matrix $CL_{10}$ (top-left), $H_{5}$ (top-right),
$K_{10}$ (bottom-left) and $Cy_{15}$ (bottom-right). These results
have been obtained for the values of the parameters reported in Table
\ref{tab:simulation-parameters-1}, for $\sigma_{1}=0.01$, $\sigma_{2}=\sigma_{3}=0.1$,
$\sigma_{4}=\sigma_{5}=1$, $Z\left(t\right)=e^{-t}\overline{J}$
and $\protect\overrightarrow{H}\left(t\right)=sin\left(2\pi t\right)\protect\overrightarrow{1}$
and with the statistics evaluated through $10,000$ Monte Carlo simulations.
The match is good even if $\sigma_{4}$ and $\sigma_{5}$ are large.}}}
\end{figure}

\noindent The Figures \ref{fig:circular-ladder-graph-simulation-1}
- \ref{fig:circulant-graph-simulation} show the numerical comparison
obtained with the first-order perturbative expansion. For simplicity
in this case we have chosen $\sigma_{4}=\sigma_{5}=0$, since according
to \ref{eq:covariance} these two parameters affect the covariance
only at a higher order. These figures report both the results obtained
from the exact network equations \ref{eq:rate-model-exact-equations-2}
(blue lines) and from the first-order perturbative expansion (red
lines), the latter being generated with the equations \ref{eq:perturbative-equation-1}
- \ref{eq:perturbative-equation-4}. Moreover we have shown the comparison
with the analytic results for the variance, covariance and correlation
generated by the formulae \ref{eq:covariance} - \ref{eq:variance}
(green lines). Instead the Figure \ref{fig:second-order-simulation}
shows the results for the second-order perturbative expansion, obtained
from the equations \ref{eq:perturbative-equation-1} - \ref{eq:perturbative-equation-12}.
For the sake of brevity, here we have reported only the results for
the correlation, but we have not shown the comparison with its analytic
formula (green lines), due to the complexity of the higher order terms
of the variance and covariance. In this case we have used $\sigma_{4}=\sigma_{5}=1$,
$Z\left(t\right)=e^{-t}\overline{J}$ and $\overrightarrow{H}\left(t\right)=sin\left(2\pi t\right)\overrightarrow{1}$,
where $\overrightarrow{1}$ is the vector whose entries are all ones.
For all these simulations we have used the parameters reported in
Table \ref{tab:simulation-parameters-1}, while the statistics have
been calculated with $10,000$ Monte Carlo simulations. Moreover the
equations \ref{eq:rate-model-exact-equations-2} and \ref{eq:perturbative-equation-1}
- \ref{eq:perturbative-equation-12} have been solved numerically
using the Euler-Maruyama scheme, while the time-integrals involved
in the formulae \ref{eq:covariance-part-1}, \ref{eq:covariance-part-2}
and \ref{eq:covariance-part-3} have been calculated with the trapezoidal
rule. The integration time step is $\Delta t=0.1$. The covariance
and correlation have always been calculated between the $0$th and
the $1$st neuron, while the potentials and the variances have been
reported only for the former. The general conclusion is that for small
enough values of the parameters $\sigma_{1}$-$\sigma_{5}$ there
is a very good agreement between the real network and the first-order
perturbative expansion and that for higher values of these parameters
the second-order expansion should be used. The match depends on the
dynamics of the neurons, on the synaptic connectivity and on the network
size, and in the case of the variance, covariance and correlation,
it also depends on the number of Monte Carlo simulations used to evaluate
the statistics.

\begin{table}
\begin{centering}
\begin{tabular}{|c||c||c||c|}
\hline 
Neuron & Input & Synaptic Weights & Sigmoid Function\tabularnewline
\hline 
\hline 
$\tau=1$ & $\overline{I}=0$ & $\Lambda=1$ & $T_{MAX}=1$\tabularnewline
\hline 
$C_{2}=0.4$ & $C_{1}=0.3$ & $C_{3}=0.5$ & $\lambda=1$\tabularnewline
\hline 
 &  &  & $V_{T}=0$\tabularnewline
\hline 
\end{tabular}
\par\end{centering}

\centering{}\caption[{\footnotesize{Parameters for the simulation of the rate model}}]{{\footnotesize{}}\label{tab:simulation-parameters-1}{\footnotesize{Parameters
used for all the numerical simulations of the Figures \ref{fig:circular-ladder-graph-simulation-1}
- \ref{fig:correlation-circulant-network}.}}}
\end{table}

{\footnotesize{}}
{\footnotesize \par}

\section{\noindent \label{sec:Correlation as a function of the input}Correlation
as a function of the input}

\begin{figure}
\noindent \begin{centering}
\includegraphics[scale=0.3]{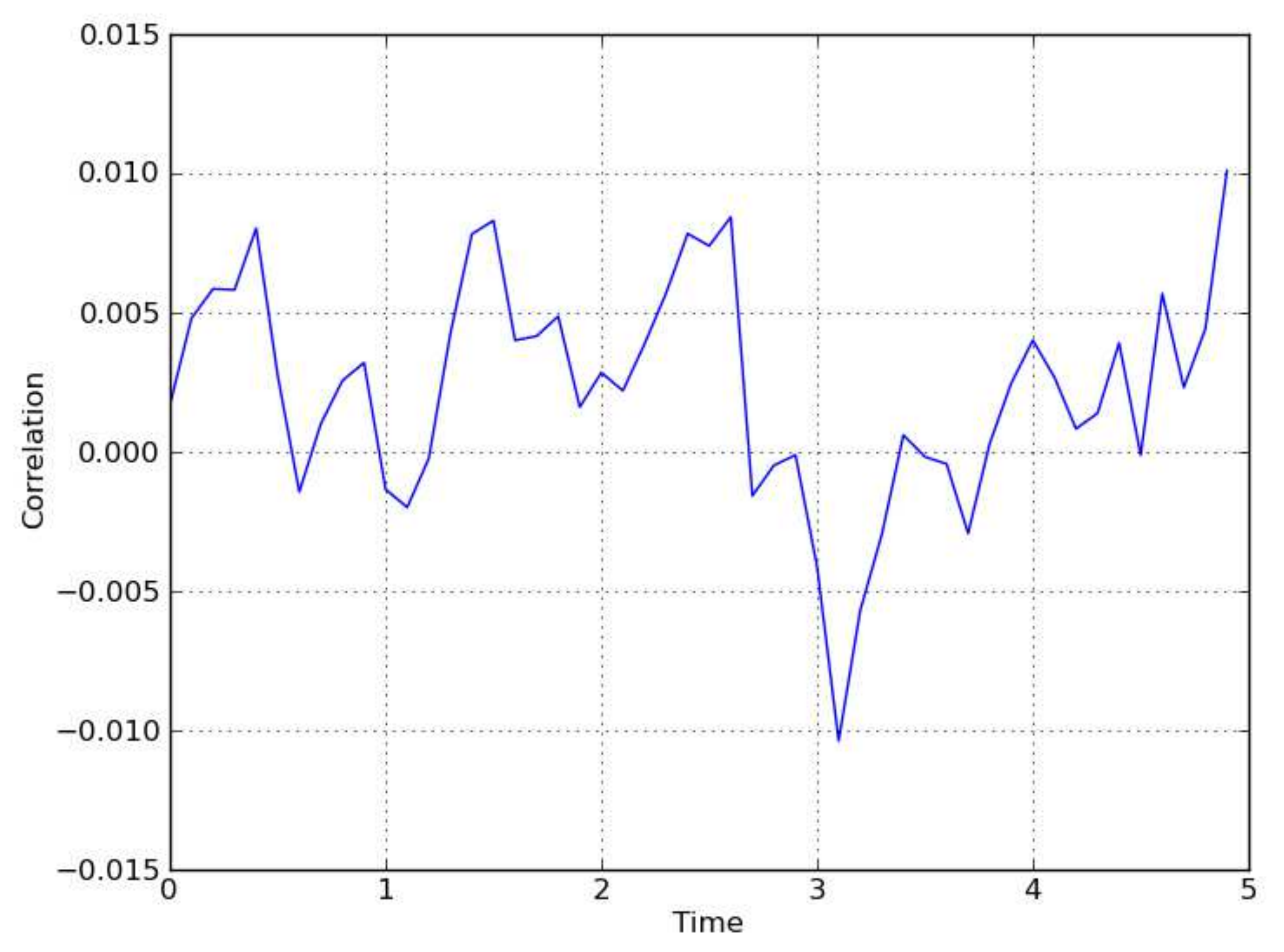}\includegraphics[scale=0.3]{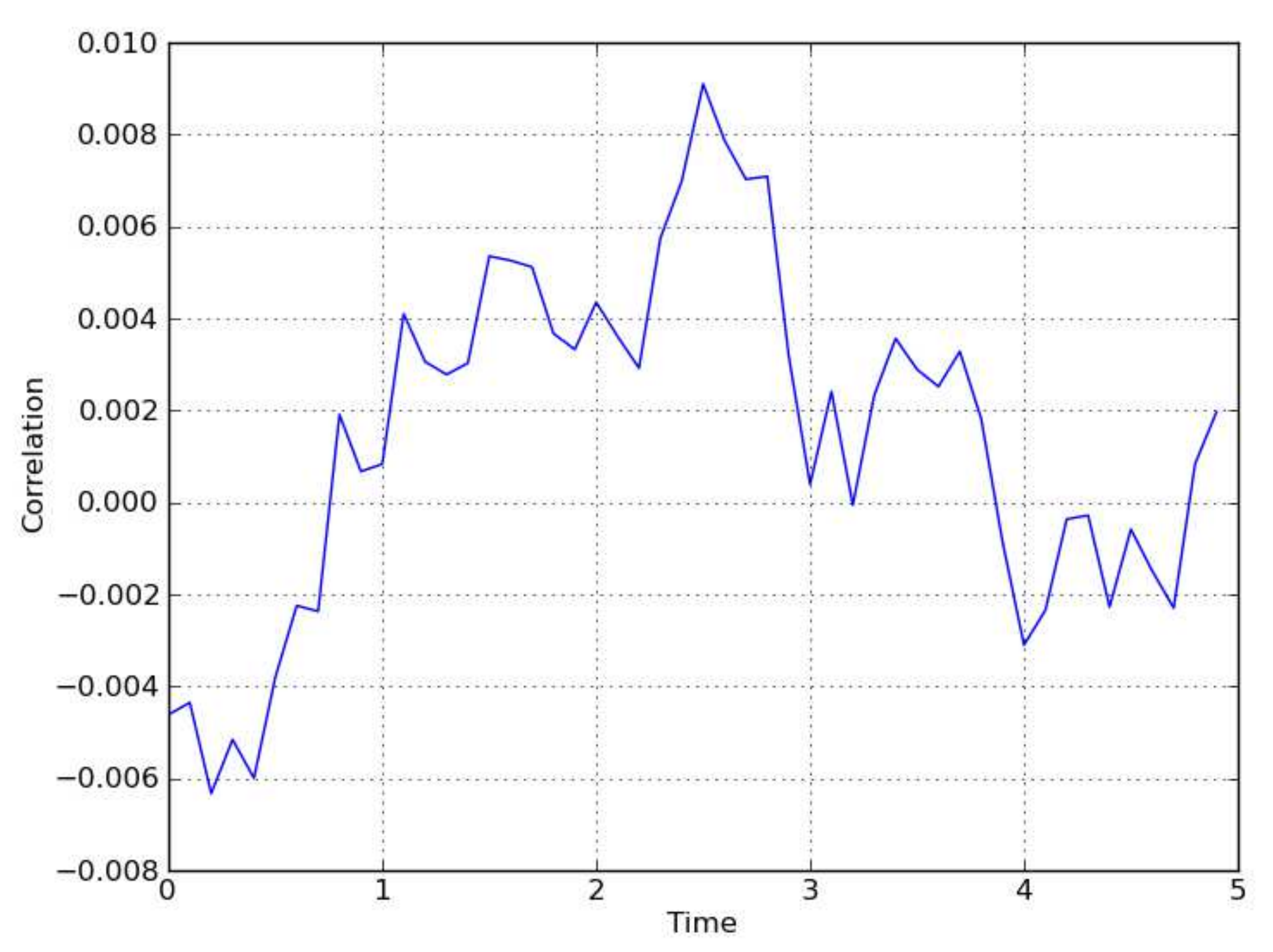}
\par\end{centering}

\noindent \begin{centering}
\includegraphics[scale=0.3]{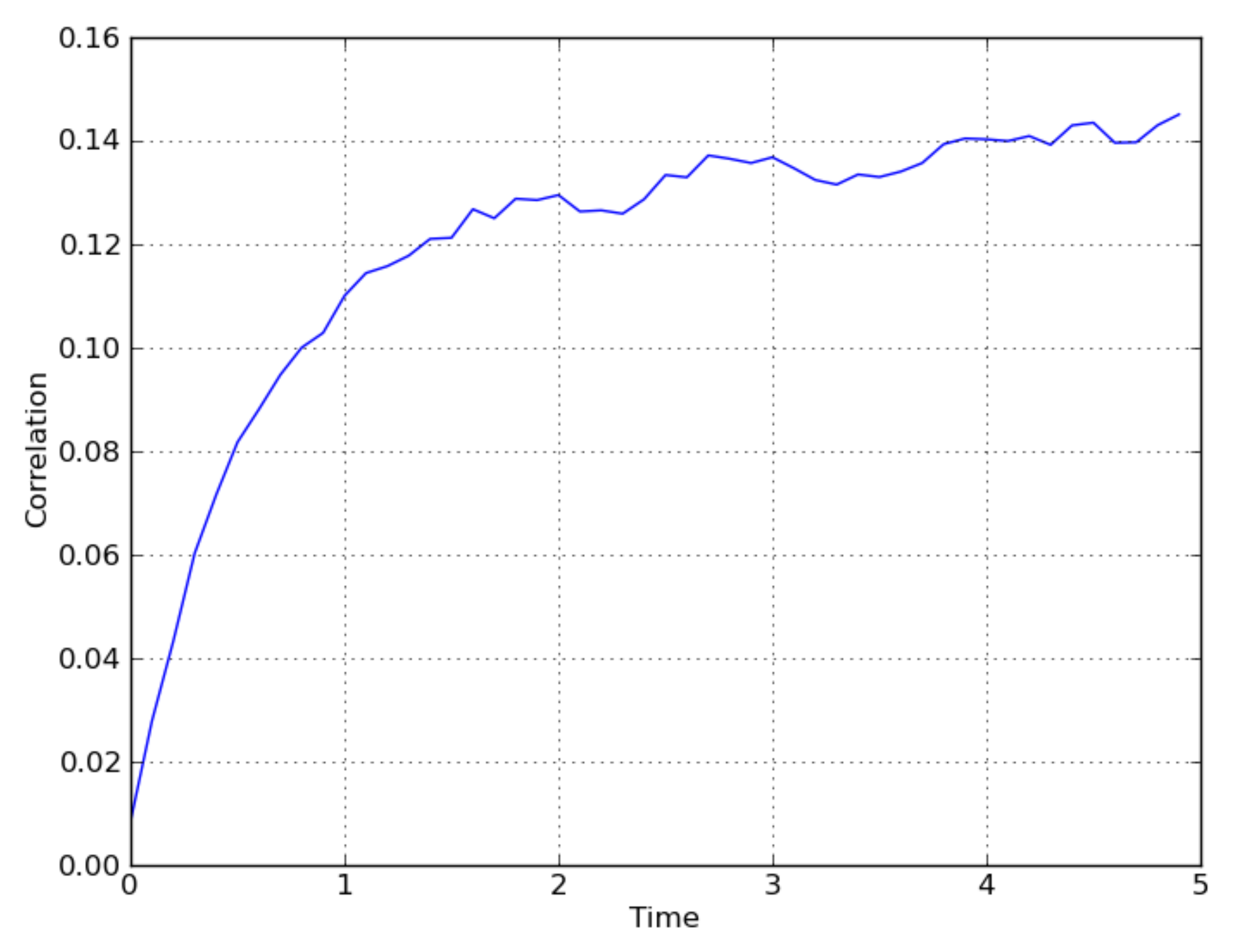}
\par\end{centering}

\caption[{\footnotesize{Correlation as a function of the input}}]{{\footnotesize{\label{fig:independence-induced-by-effective-connectivity}
Correlation function obtained for the graph $Cy_{5}$ and $\overline{I}=-5$
(top-left), $5$ (top-right) and $0$ (bottom). These results have
been obtained from the exact equations \ref{eq:rate-model-exact-equations-2},
numerically solved using the Euler-Maruyama scheme with integration
time step $\Delta t=0.1$ and with $50,000$ Monte Carlo simulations.
The parameters used are $C_{1}=C_{2}=C_{3}=0$, $\sigma_{1}=\sigma_{2}=\sigma_{3}=0.1$
and $\sigma_{4}=\sigma_{5}=0$, while all the remaining parameters
are those of Table \ref{tab:simulation-parameters-1}. From this figure
it is possible to see that the correlation between pairs of neurons
strongly decreases for high values of $\left|\overline{I}\right|$,
confirming its relation with the effective connectivity matrix $\overline{J}S'\left(\mu\right)$.}}}
\end{figure}

\noindent From the formulae \ref{eq:covariance-part-1}, \ref{eq:covariance-part-2},
\ref{eq:covariance-part-3} and \ref{eq:Phi-matrix-block-circulant-case}
the effect of the non-linearity introduced by the sigmoid function
$S\left(V\right)$ is evident. Through its slope, it generates an
effective connectivity matrix $\overline{J}S'\left(\mu\right)$, which
can be interpreted as the real connectivity matrix of the system if
it were linear. Now, the stationary solution $\mu$ depends on the
external input current $\overline{I}$ through the formula \ref{eq:perturbative-equation-1},
therefore the effective synaptic strength and the correlation structure
depend on $\overline{I}$ as well. In particular, it is interesting
to observe that if $\left|\overline{I}\right|$ is very large, then
$\left|\mu\right|$ is also very large, therefore $S'\left(\mu\right)$
and the entries of the effective connectivity matrix are small. In
other words, the neurons become (effectively) disconnected. An important
consequence of this phenomenon is that, for $C_{1}=C_{2}=C_{3}=0$
and for large values of $\left|\overline{I}\right|$, the neurons
become independent, even if the size of the network is finite. This
intuition is confirmed numerically in Figure \ref{fig:independence-induced-by-effective-connectivity},
which has been obtained for the graph $Cy_{5}$ (which is made of
$10$ neurons) simulated with the exact equations \ref{eq:rate-model-exact-equations-2},
for $\overline{I}=-5,0,5$ and $50,000$ Monte Carlo simulations.
The sources of randomness have intensities $\sigma_{1}=\sigma_{2}=\sigma_{3}=0.1$,
and moreover $\sigma_{4}=\sigma_{5}=0$, while all the remaining parameters
are those of Table \ref{tab:simulation-parameters-1}. As usual, the
numerical scheme is the Euler-Maruyama one, with integration time
step $\Delta t=0.1$.

\section{\noindent \label{sec:Failure of the mean field theory}Failure of
the mean-field theory}

\noindent In this section we show three different reasons that invalidate
the use of the mean-field theory for the mathematical analysis of
a neural network. A neural network is generally described by a large
set of stochastic differential equations, that makes it hard to understand
the underlying behavior of the system. However, if the neurons become
independent, their dynamics can be described with the mean-field theory
using a highly reduced set of equations, that are much simpler to
analyze. For this reason the mean-field theory is a powerful tool
that can be used to understand the network. One of the mechanisms
through which the independence of the neurons can be obtained is the
phenomenon known as \textit{propagation of chaos} \cite{journals/siamads/TouboulHF12}\cite{baladron:inserm-00732288}\cite{samuelides:inria-00529560}\cite{Touboul11}.
Propagation of chaos refers to the fact that, if we choose independent
initial conditions for the membrane potentials at $t=0$ (which may
be called \textit{initial chaos}), then the neurons are always perfectly
independent $\forall t>0$. Therefore the term \textit{propagation}
refers to the ``transfer'' of the independence of the membrane potentials
from $t=0$ to $t>0$. Under simplified assumptions about the nature
of the network (namely that the other sources of randomness in the
system, in our case the Brownian motions and the synaptic weights,
are independent), propagation of chaos does occur in the so called
\textit{thermodynamic limit }of the system, namely when the number
of neurons in the system grows to infinity. However in Sections \ref{sub:Independence does not occur if C1, C2 or C3 are not equal to zero},
\ref{sub:Propagation of chaos does not occur for a general connectivity matrix}
and \ref{sub:Stochastic synchronization} we show that for a system
with correlated Brownian motions, initial conditions and synaptic
weights, with a general connectivity matrix or with an arbitrarily
large (but still finite) size, the correlation between pairs of neurons
can be high. Therefore in general the neurons cannot be independent,
invalidating the use of the mean-field theory.

\subsection{\noindent \label{sub:Independence does not occur if C1, C2 or C3 are not equal to zero}Independence
does not occur for $N\rightarrow\infty$ if $C_{1}$, $C_{2}$ or
$C_{3}$ are not equal to zero}

\noindent Let us consider the case when at least one of the parameters
$C_{1}$, $C_{2}$ and $C_{3}$ (defined by \ref{eq:Brownian-covariance},
\ref{eq:initial-conditions-covariance} and \ref{eq:synaptic-weights-covariance})
is not equal to zero. For example we analyze the term proportional
to $C_{1}$ in the formula \ref{eq:covariance-part-1}, for a fully
connected network. Using the technique developed in Section \ref{sub:Block circulant matrices with circulant blocks},
it is easy to prove that this term for $i\neq j$ is:

\begin{onehalfspace}
\begin{center}
{\small{
\begin{align*}
 & C_{1}{\displaystyle \sum\limits _{\substack{k,l=0\\
k\neq l
}
}^{N-1}}\int_{0}^{t}\left[\Phi\left(t-s\right)\right]_{ik}\left[\Phi\left(t-s\right)\right]_{jl}ds\\
\\
 & =\frac{C_{1}}{2}\left\{ \left(1-\frac{1}{N}\right)\frac{1}{\frac{1}{\tau}-\Lambda S'\left(\mu\right)}\left[1-e^{-2\left(\frac{1}{\tau}-\Lambda S'\left(\mu\right)\right)t}\right]+\frac{1}{N}\frac{1}{\frac{1}{\tau}+\frac{\Lambda S'\left(\mu\right)}{N-1}}\left[1-e^{-2\left(\frac{1}{\tau}+\frac{\Lambda S'\left(\mu\right)}{N-1}\right)t}\right]\right\} 
\end{align*}
}}
\par\end{center}{\small \par}
\end{onehalfspace}

\noindent while for $i=j$ it is:

\begin{onehalfspace}
\begin{center}
{\small{
\begin{align*}
 & C_{1}{\displaystyle \sum\limits _{\substack{k,l=0\\
k\neq l
}
}^{N-1}}\int_{0}^{t}\left[\Phi\left(t-s\right)\right]_{ik}\left[\Phi\left(t-s\right)\right]_{il}ds\\
\\
 & =\frac{C_{1}}{2}\left(1-\frac{1}{N}\right)\left\{ \frac{1}{\frac{1}{\tau}-\Lambda S'\left(\mu\right)}\left[1-e^{-2\left(\frac{1}{\tau}-\Lambda S'\left(\mu\right)\right)t}\right]-\frac{1}{\frac{1}{\tau}+\frac{\Lambda S'\left(\mu\right)}{N-1}}\left[1-e^{-2\left(\frac{1}{\tau}+\frac{\Lambda S'\left(\mu\right)}{N-1}\right)t}\right]\right\} 
\end{align*}
}}
\par\end{center}{\small \par}
\end{onehalfspace}

\noindent So the covariance (and therefore also the correlation) does
not go to zero for $N\rightarrow\infty$, or in other words the neurons
are not independent, even in the thermodynamic limit.

\noindent The reader can easily check that the same result holds for
the terms of the covariance proportional to $C_{2}$ and $C_{3}$.

\subsection{\noindent \label{sub:Propagation of chaos does not occur for a general connectivity matrix}Propagation
of chaos does not occur for a general connectivity matrix}

\begin{figure}
\begin{centering}
\includegraphics[scale=0.3]{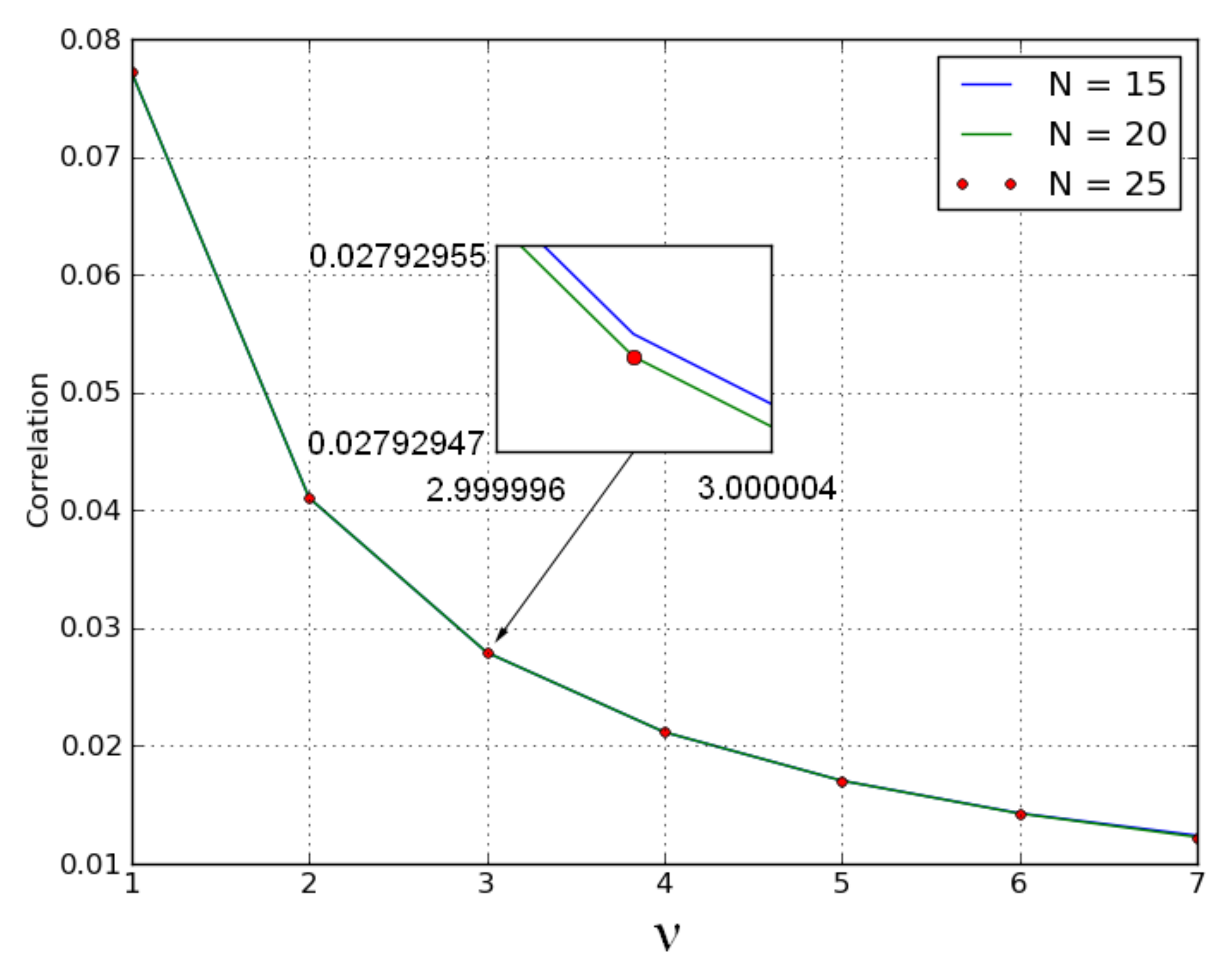}\includegraphics[scale=0.3]{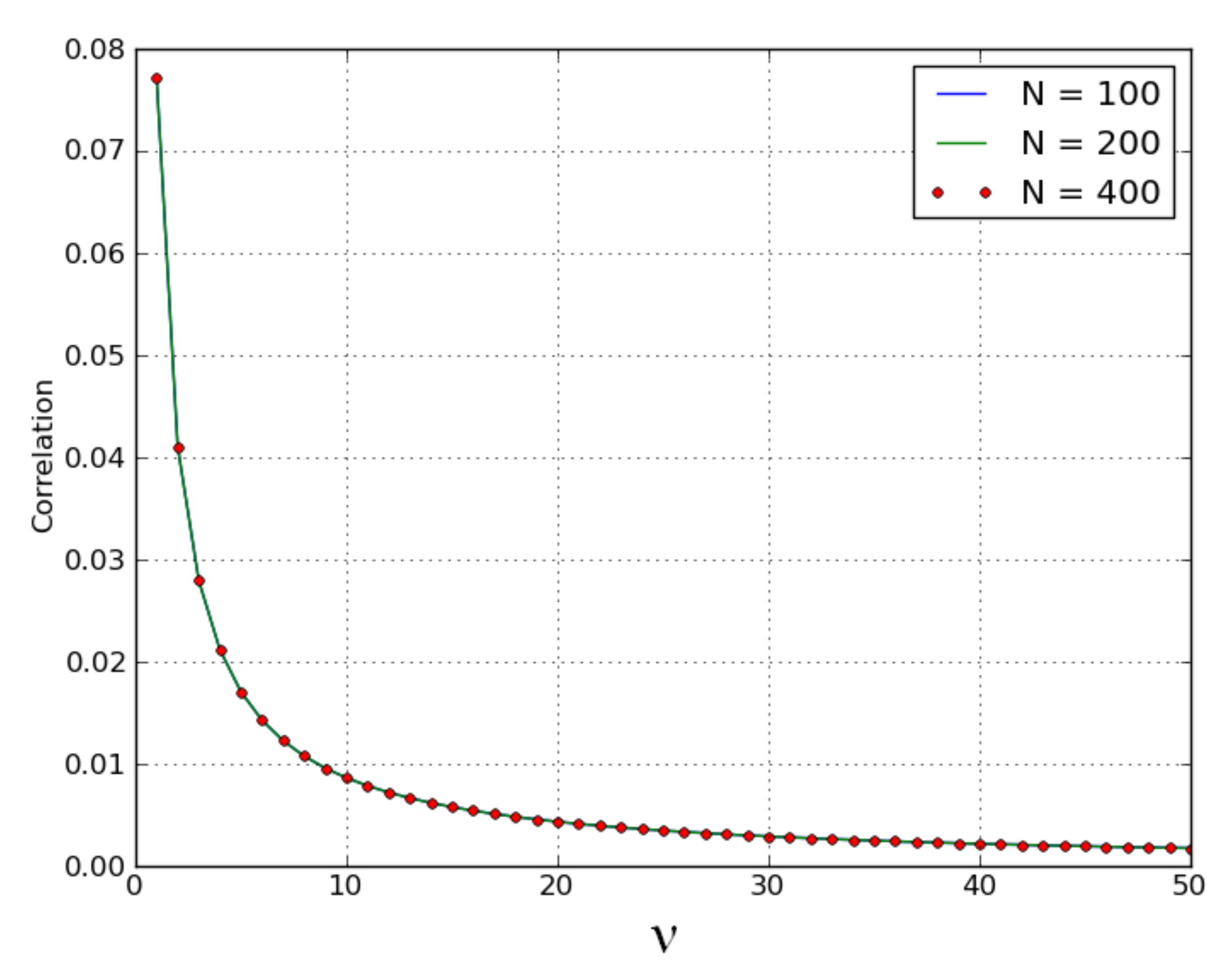}
\par\end{centering}

\caption[{\footnotesize{Correlation as a function of the number of incoming
connections per neuron}}]{{\footnotesize{\label{fig:correlation-circulant-network}Propagation
of chaos for $t=1$ as a function of $\nu=\frac{M}{2}$, in the case
of a circulant connectivity matrix. This result has been obtained
for $C_{1}=C_{2}=C_{3}=0$ (while all the remaining parameters are
those of Table \ref{tab:simulation-parameters-1}), using the analytic
formula \ref{eq:covariance-circulant-network} (normalized with the
variance).}}}
\end{figure}

\noindent We study propagation of chaos as a function of the number
of connections in the circulant network. To this purpose, we have
to set $C_{2}=0$ (initial chaos) and also $C_{1}=C_{3}=0$, because
otherwise the neurons cannot be independent, as explained in Section
\ref{sub:Independence does not occur if C1, C2 or C3 are not equal to zero}.
Using the formulae \ref{eq:covariance-part-1}, \ref{eq:covariance-part-2},
\ref{eq:covariance-part-3} and \ref{eq:Phi-matrix-block-circulant-case}
we obtain that in this case the covariance is:

\begin{onehalfspace}
\begin{center}
{\small{
\begin{align}
Cov\left(V_{i}\left(t\right),V_{j}\left(t\right)\right)= & \sigma_{1}^{2}\int_{0}^{t}\left[\Phi\left(t-s\right)\Phi^{T}\left(t-s\right)\right]_{ij}ds+\sigma_{2}^{2}\left[\Phi\left(t\right)\Phi^{T}\left(t\right)\right]_{ij}\nonumber \\
\nonumber \\
 & +\sigma_{3}^{2}\frac{S^{2}\left(\mu\right)}{M}\sum_{k=0}^{N-1}\left[\int_{0}^{t}\Phi_{ik}\left(t-s\right)ds\right]\left[\int_{0}^{t}\Phi_{jk}\left(t-s\right)ds\right]\label{eq:covariance-circulant-network}
\end{align}
}}
\par\end{center}{\small \par}
\end{onehalfspace}

\noindent where:

\begin{onehalfspace}
\begin{center}
{\small{
\begin{align*}
 & \int_{0}^{t}\left[\Phi\left(t-s\right)\Phi^{T}\left(t-s\right)\right]_{ij}ds=\frac{1}{2N}\sum_{k=0}^{N-1}\frac{cos\left[\frac{2\pi}{N}k\left(i-j\right)\right]}{-\frac{1}{\tau}+e_{k}S'\left(\mu\right)}\left\{ 1-e^{2\left[-\frac{1}{\tau}+S'\left(\mu\right)e_{k}\right]t}\right\} \\
\\
 & \left[\Phi\left(t\right)\Phi^{T}\left(t\right)\right]_{ij}=\frac{1}{N}\sum_{k=0}^{N-1}e^{2\left[-\frac{1}{\tau}+S'\left(\mu\right)e_{k}\right]t}cos\left[\frac{2\pi}{N}k\left(i-j\right)\right]\\
\\
 & \sum_{k=0}^{N-1}\left[\int_{0}^{t}\Phi_{ik}\left(t-s\right)ds\right]\left[\int_{0}^{t}\Phi_{jk}\left(t-s\right)ds\right]\\
 & =\frac{1}{N^{2}}\sum_{l,m=0}^{N-1}e^{\frac{2\pi}{N}li\iota}e^{\frac{2\pi}{N}mj\iota}\left[\sum_{k=0}^{N-1}e^{-\frac{2\pi}{N}\left(l+m\right)k\iota}\right]\left\{ \frac{1-e^{\left[-\frac{1}{\tau}+S'\left(\mu\right)e_{l}\right]t}}{-\frac{1}{\tau}+e_{l}S'\left(\mu\right)}\right\} \left\{ \frac{1-e^{\left[-\frac{1}{\tau}+S'\left(\mu\right)e_{m}\right]t}}{-\frac{1}{\tau}+e_{m}S'\left(\mu\right)}\right\} \\
 & =\frac{1}{N}\sum_{l=0}^{N-1}\left\{ \frac{1-e^{\left[-\frac{1}{\tau}+S'\left(\mu\right)e_{l}\right]t}}{-\frac{1}{\tau}+e_{l}S'\left(\mu\right)}\right\} ^{2}cos\left[\frac{2\pi}{N}l\left(i-j\right)\right]
\end{align*}
}}
\par\end{center}{\small \par}
\end{onehalfspace}

\noindent while the eigenvalues $e_{k}$ are given by formula \ref{eq:circulant-network-eigenvalues}
or by formula \ref{eq:complete-network-eigenvalues}. Now, for $N\rightarrow\infty$
the right-hand side of formula \ref{eq:covariance-circulant-network}
converges to a non-zero function (see Figure \ref{fig:correlation-circulant-network}),
therefore for every finite value of $\nu$ (which is the number of
incoming connections per neuron divided by $2$) propagation of chaos
does not occur.

\noindent Moreover correlation decreases with $\nu$, therefore propagation
of chaos occurs in the circulant network only in the thermodynamic
limit $N\rightarrow\infty$ and if $\nu$ is an increasing function
of $N$, namely if $\underset{N\rightarrow\infty}{lim}\nu=\infty$.
For example, in the fully connected network $\nu=\left\lfloor \frac{N}{2}\right\rfloor $,
so it explains why in this case correlation goes to zero in the thermodynamic
limit. Instead in a network described by a cycle graph, perfect decorrelation
is never possible, also for $N\rightarrow\infty$, since $\nu=1$.
In other words, having infinitely many neurons is not a sufficient
condition for getting propagation of chaos, because also infinite
connections per neuron are required.

\subsection{\noindent \label{sub:Stochastic synchronization}Stochastic synchronization}

\noindent In this section we show that for every finite and arbitrarily
large number of neurons $N$ in the network, it is possible to choose
special values of the parameters of the system such that, at some
finite and arbitrarily large time instant $\overline{t}$, the correlation
between pairs of neurons is (approximately) $1$. In general $\overline{t}$
increases with $N$.

\subsubsection{\noindent \label{sub:The general theory}The general theory}

\noindent We show that even when $C_{1}=C_{2}=C_{3}=0$, if the matrix
$A=-\frac{1}{\tau}Id_{N}+\overline{J}S'\left(\mu\right)$ has an eigenvalue
of multiplicity $1$ with non-negative real part, while all the other
eigenvalues have negative real parts, then correlation goes to $1$
for $t\rightarrow+\infty$, for every finite $N$. In other terms,
the stochastic components of the membrane potentials become perfectly
synchronized. From now on we refer to this phenomenon as \textit{stochastic
synchronization}. To prove this, we suppose that $A$ has an eigenvalue
$\overline{a}$ with non-negative real part and with a generic multiplicity
$m>0$, while all the other eigenvalues have negative real parts.
Now we recall that $e^{At}=Qe^{Dt}Q^{-1}$, where $D$ is the diagonal
matrix of the eigenvalues of $A$, and $Q$ is the matrix of its eigenvectors.
So for $t\rightarrow+\infty$ we have:

\begin{onehalfspace}
\begin{center}
{\small{
\[
e^{Dt}\prec diag(0,0,...,0,\underset{m\mathrm{-times}}{\underbrace{e^{\overline{a}t},e^{\overline{a}t},...,e^{\overline{a}t}}},0,0,...,0)
\]
}}
\par\end{center}{\small \par}
\end{onehalfspace}

\noindent where $\prec$ means \textit{dominated by}, because all
the eigenvalues have negative real part but $\overline{a}$. If the
$\overline{a}$s are the $r$-th, $\left(r+1\right)$-th, ..., $\left(r+m-1\right)$-th
eigenvalues of $A$ and if we call $Q^{-1}=B$ in order to simplify
the notation, we obtain:

\begin{onehalfspace}
\begin{center}
{\small{
\begin{align*}
Qe^{Dt}B\prec & e^{\overline{a}t}\left[\begin{array}{cccccccccccc}
0 & 0 & \cdots & 0 & Q_{0,r} & Q_{0,r+1} & \cdots & Q_{0,r+m-1} & 0 & 0 & \cdots & 0\\
0 & 0 & \cdots & 0 & Q_{1,r} & Q_{1,r+1} & \cdots & Q_{1,r+m-1} & 0 & 0 & \cdots & 0\\
\vdots & \vdots & \ddots & \vdots & \vdots & \vdots & \ddots & \vdots & \vdots & \vdots & \ddots & \vdots\\
0 & 0 & \cdots & 0 & Q_{N-1,r} & Q_{N-1,r+1} & \cdots & Q_{N-1,r+m-1} & 0 & 0 & \cdots & 0
\end{array}\right]\\
\\
 & \times\left[\begin{array}{cccc}
B_{0,0} & B_{0,1} & ... & B_{0,N-1}\\
B_{1,0} & B_{1,1} & ... & B_{1,N-1}\\
\vdots & \vdots & \ddots & \vdots\\
B_{N-1,0} & B_{N-1,1} & ... & B_{N-1,N-1}
\end{array}\right]
\end{align*}
}}
\par\end{center}{\small \par}
\end{onehalfspace}

\noindent {\small{and therefore:}}{\small \par}

\begin{onehalfspace}
\begin{center}
{\small{
\begin{align*}
e^{At}= & Qe^{Dt}B\prec e^{\overline{a}t}E\\
\\
E_{pq}= & \sum_{k=0}^{m-1}Q_{p,r+k}B_{r+k,q}
\end{align*}
}}
\par\end{center}{\small \par}
\end{onehalfspace}

\noindent This means that:

\begin{onehalfspace}
\begin{center}
{\small{
\begin{align}
Cov\left(V_{i}\left(t\right),V_{j}\left(t\right)\right)= & \sigma_{1}^{2}{\displaystyle \sum_{k=0}^{N-1}}\int_{0}^{t}\left[e^{A\left(t-s\right)}\right]_{ik}\left[e^{A\left(t-s\right)}\right]_{jk}ds+\sigma_{2}^{2}{\displaystyle \sum_{k=0}^{N-1}}\left[e^{At}\right]_{ik}\left[e^{At}\right]_{jk}\nonumber \\
 & +\sigma_{3}^{2}\frac{S^{2}\left(\mu\right)}{M}{\displaystyle \sum_{k=0}^{N-1}}\left\{ \int_{0}^{t}\left[e^{A\left(t-s\right)}\right]_{ik}ds\right\} \left\{ \int_{0}^{t}\left[e^{A\left(t-s\right)}\right]_{jk}ds\right\} \nonumber \\
\prec & \left[\frac{\sigma_{1}^{2}}{2\overline{a}}+\sigma_{2}^{2}+\frac{\sigma_{3}^{2}}{\overline{a}^{2}}\frac{S^{2}\left(\mu\right)}{M}\right]e^{2\overline{a}t}{\displaystyle \sum_{k=0}^{N-1}}E_{ik}E_{jk}\label{eq:stochastic-synchronization-covariance}
\end{align}
}}
\par\end{center}{\small \par}
\end{onehalfspace}

\noindent so the variance is:

\begin{onehalfspace}
\begin{center}
{\small{
\begin{equation}
Var\left(V_{i}\left(t\right)\right)=Cov\left(V_{i}\left(t\right),V_{i}\left(t\right)\right)\prec\left[\frac{\sigma_{1}^{2}}{2\overline{a}}+\sigma_{2}^{2}+\frac{\sigma_{3}^{2}}{\overline{a}^{2}}\frac{S^{2}\left(\mu\right)}{M}\right]e^{2\overline{a}t}\sum_{k=0}^{N-1}\left(E_{ik}\right)^{2}\label{eq:stochastic-synchronization-variance}
\end{equation}
}}
\par\end{center}{\small \par}
\end{onehalfspace}

\noindent Therefore the correlation is:

\begin{onehalfspace}
\begin{center}
{\small{
\begin{equation}
Corr\left(V_{i}\left(t\right),V_{j}\left(t\right)\right)=\frac{Cov\left(V_{i}\left(t\right),V_{j}\left(t\right)\right)}{\sqrt{Var\left(V_{i}\left(t\right)\right)Var\left(V_{j}\left(t\right)\right)}}\rightarrow\frac{{\displaystyle \sum_{k=0}^{N-1}}E_{ik}E_{jk}}{\sqrt{\left[{\displaystyle \sum_{k=0}^{N-1}}\left(E_{ik}\right)^{2}\right]\left[{\displaystyle \sum_{k=0}^{N-1}}\left(E_{jk}\right)^{2}\right]}}\label{eq:stochastic-synchronization-correlation}
\end{equation}
}}
\par\end{center}{\small \par}
\end{onehalfspace}

\noindent when $t\rightarrow+\infty$. Now, in the special case $m=1$
we obtain:

\begin{onehalfspace}
\begin{center}
{\small{
\begin{align*}
E_{pq}= & Q_{pr}B_{rq}\\
\\
\sum_{k=0}^{N-1}E_{ik}E_{jk}= & Q_{ir}Q_{jr}\sum_{k=0}^{N-1}\left(B_{rk}\right)^{2}\\
\\
\sum_{k=0}^{N-1}\left(E_{ik}\right)^{2}= & \left(Q_{ir}\right)^{2}\sum_{k=0}^{N-1}\left(B_{rk}\right)^{2}
\end{align*}
}}
\par\end{center}{\small \par}
\end{onehalfspace}

\noindent so we conclude that $Corr\left(V_{i}\left(t\right),V_{j}\left(t\right)\right)\rightarrow1$
when $t\rightarrow+\infty$. This proves that if $C_{1}=C_{2}=C_{3}=0$
and the matrix $A$ has an eigenvalue of multiplicity $1$ with non-negative
real part while all the other eigenvalues have negative real parts,
then propagation of chaos does not occur. For continuity, for every
finite $N$ we have that $Corr\left(V_{i}\left(t\right),V_{j}\left(t\right)\right)\rightarrow1$
also for $\mathcal{R}\left(\overline{a}\right)\rightarrow0^{-}$ (where
$\mathcal{R}$ means \textit{the real part of}), i.e. correlation
is very big also when the system is stable but close to the instability
region $\mathcal{R}\left(\overline{a}\right)>0$. It is also interesting
to observe that, due to the Perron-Frobenius theorem \cite{1406483},
if $\Lambda>0$ and if $\overline{J}$ is an irreducible matrix (namely
if its corresponding directed graph is \textit{strongly connected},
which means that it is possible to reach each vertex in the graph
from any other vertex, by moving on the edges according to their connectivity
directions), then it has a unique largest positive eigenvalue, which
can be used to generate stochastic synchronization. We conclude that
in general propagation of chaos does not always occur, even if $C_{1}=C_{2}=C_{3}=0$,
therefore this invalidates the use of the mean-field theory, at least
in this special case.

\subsubsection{\noindent \label{sub:The example of the fully connected network}The
example of the fully connected network}

\noindent We show how to set the parameters of the system such that
the phenomenon of stochastic synchronization does occur. For simplicity
we assume a fully connected network. In this case, from formula \ref{eq:complete-network-eigenvalues},
we know that the matrix $A$ has eigenvalues:

\begin{onehalfspace}
\begin{center}
{\small{
\begin{equation}
\begin{array}{ccc}
a_{0}=-\frac{1}{\tau}+\Lambda S'\left(\mu\right), &  & a_{1}=-\frac{1}{\tau}-\frac{\Lambda S'\left(\mu\right)}{N-1}\end{array}\label{eq:complete-network-eigenvalues-matrix-A}
\end{equation}
}}
\par\end{center}{\small \par}
\end{onehalfspace}

\noindent The multiplicity of $a_{0}$ and $a_{1}$ is respectively
$1$ and $N-1$, therefore in order to obtain the stochastic synchronization,
according to Section \ref{sub:The general theory}, we have to set
$a_{0}\geq0$. Let us consider the case $a_{0}=0$, namely $\Lambda S'\left(\mu\right)=\frac{1}{\tau}$.
Now, since:

\begin{onehalfspace}
\begin{center}
{\small{
\begin{equation}
S'\left(\mu\right)=\lambda\left[S\left(\mu\right)-\frac{S^{2}\left(\mu\right)}{T_{MAX}}\right]\label{eq:derivative-sigmoid-function}
\end{equation}
}}
\par\end{center}{\small \par}
\end{onehalfspace}

\noindent we obtain the algebraic equation:

\begin{onehalfspace}
\begin{center}
{\small{
\[
\Lambda\lambda\left[S\left(\mu\right)-\frac{S^{2}\left(\mu\right)}{T_{MAX}}\right]=\frac{1}{\tau}
\]
}}
\par\end{center}{\small \par}
\end{onehalfspace}

\noindent whose solutions are:

\begin{onehalfspace}
\begin{center}
{\small{
\begin{equation}
S\left(\mu_{1,2}\right)=T_{MAX}\frac{1\pm\sqrt{1-\frac{4}{\tau\Lambda\lambda T_{MAX}}}}{2}\label{eq:sigmoid-function-at-stationary-points}
\end{equation}
}}
\par\end{center}{\small \par}
\end{onehalfspace}

\noindent where $\mu_{1,2}$ are two possible stationary solutions
of the membrane potential. Moreover, from equation \ref{eq:perturbative-equation-1}
we know that:

\begin{onehalfspace}
\begin{center}
{\small{
\begin{equation}
\mu_{1,2}=\tau\left[\Lambda S\left(\mu_{1,2}\right)+\overline{I}\right]\label{eq:stationary-points}
\end{equation}
}}
\par\end{center}{\small \par}
\end{onehalfspace}

\noindent Putting together the formulae \ref{eq:sigmoid-function-at-stationary-points}
and \ref{eq:stationary-points} we obtain:

\begin{onehalfspace}
\begin{flushleft}
{\small{
\begin{equation}
\mu_{1,2}=\tau\left(\Lambda T_{MAX}\frac{1\pm\sqrt{1-\frac{4}{\tau\Lambda\lambda T_{MAX}}}}{2}+\overline{I}\right)\label{eq:stationary-points-final-formula}
\end{equation}
}}
\par\end{flushleft}{\small \par}
\end{onehalfspace}

\noindent Replace this value of $\mu_{1,2}$ in \ref{eq:stationary-points}
to obtain the final result:

\begin{onehalfspace}
\begin{center}
{\small{
\begin{equation}
T_{MAX}\frac{1\pm\sqrt{1-\frac{4}{\tau\Lambda\lambda T_{MAX}}}}{2}=S\left(\tau\left(\Lambda T_{MAX}\frac{1\pm\sqrt{1-\frac{4}{\tau\Lambda\lambda T_{MAX}}}}{2}+\overline{I}\right)\right)\label{eq:parameters-constraint-stochastic-synchronization}
\end{equation}
}}
\par\end{center}{\small \par}
\end{onehalfspace}

\begin{figure}
\begin{centering}
\includegraphics[scale=0.3]{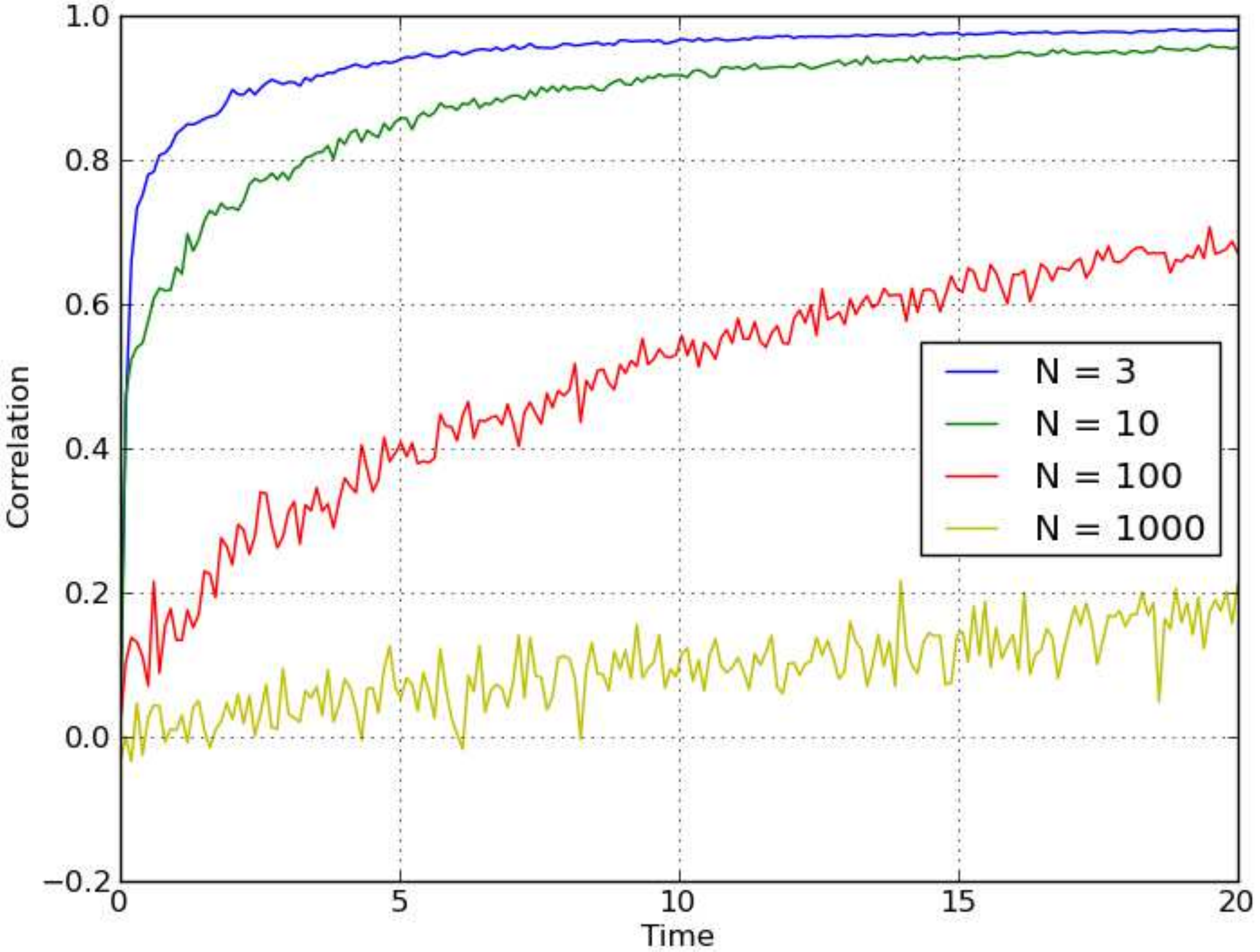}
\par\end{centering}

\caption[{\footnotesize{Stochastic synchronization}}]{{\footnotesize{\label{fig:stochastic-synchronization}Stochastic
synchronization in a fully connected network. Correlation gets closer
and closer to $1$ with a speed that depends on the number of neurons
$N$ in the system. These results have been obtained with the exact
non-linear equations \ref{eq:rate-model-exact-equations-2} and with
$1,000$ Monte Carlo simulations. The parameters are those of Table
\ref{tab:simulation-parameters-strong-2}, which are chosen in order
to satisfy the constraint \ref{eq:parameters-constraint-stochastic-synchronization}
for $\Lambda=-2\overline{I}$ and $\tau=-\frac{2}{\overline{I}}$.
The value of the external current is purposely large ($\overline{I}=-20$)
because it causes a faster convergence of the correlation to the value
$1$.}}}
\end{figure}

\noindent This non-linear algebraic equation is the constraint that
must be satisfied by all the parameters of the system in order to
have correlation equal to $1$ in the limit $t\rightarrow+\infty$.
An example of solution of this equation is $\lambda=T_{MAX}=1$, $V_{T}=0$,
$\Lambda=-2\overline{I}$ and $\tau=-\frac{2}{\overline{I}}$, $\forall\overline{I}<0$.
In this case $\mu_{1,2}=0$ and it can be used as initial condition
in order to ensure the stationarity of the system. In Figure \ref{fig:stochastic-synchronization}
we show the phenomenon of stochastic synchronization in the case of
a fully connected network, for the values of the parameters reported
in Table \ref{tab:simulation-parameters-strong-2}, which satisfy
the constraint \ref{eq:parameters-constraint-stochastic-synchronization}.
As we can see, correlation goes to $1$ more and more slowly if we
increase the number of neurons $N$ in the network. It reaches the
value $1$ asymptotically with an inverse exponential-like behavior,
with a time constant that increases with the size of the network.
For $N\rightarrow\infty$ the time constant diverges, therefore for
every finite time the system has correlation $0$. This proves that
in the thermodynamic limit there is still propagation of chaos, provided
that $C_{1}=C_{2}=C_{3}=0$. This is in perfect agreement with the
result on propagation of chaos proved in \cite{journals/siamads/TouboulHF12}\cite{baladron:inserm-00732288}\cite{Touboul11}
for independent Brownian motions, initial conditions and synaptic
weights.

\begin{table}
\begin{centering}
\begin{tabular}{|c||c||c||c|}
\hline 
Neuron & Input & Synaptic Weights & Sigmoid Function\tabularnewline
\hline 
\hline 
$\tau=0.1$ & $\overline{I}=-20$ & $\Lambda=40$ & $T_{MAX}=1$\tabularnewline
\hline 
$C_{2}=0$ & $C_{1}=0$ & $C_{3}=0$ & $\lambda=1$\tabularnewline
\hline 
 &  &  & $V_{T}=0$\tabularnewline
\hline 
\end{tabular}
\par\end{centering}

\centering{}\caption[{\footnotesize{Parameters for the simulation of the stochastic synchronization}}]{{\footnotesize{}}\label{tab:simulation-parameters-strong-2}{\footnotesize{Parameters
used for the numerical simulations of Figure \ref{fig:stochastic-synchronization}.}}}
\end{table}

\section{\noindent \label{sec:Conclusion} Conclusion}

\noindent In this work we have developed a perturbative expansion
that let us determine the dynamics and the correlation structure of
a neural network made up of a finite number of rate neurons and with
specific connectivity matrices.

\noindent The network has three independent sources of randomness
in the stochastic fluctuations of the membrane potentials (or in the
external input current), in the initial conditions and in the synaptic
weights.

\noindent Their probability distributions are supposed to be normal
with known correlation matrices, moreover we have assumed that the
system is invariant under exchange of the neural indices, without
any other assumption on the intensity of the synaptic weights.

\noindent This has allowed us to obtain analytic results, even if
the equations of the network are non-linear, which have been confirmed
numerically.

\noindent With this approach we have analyzed block circulant and
special symmetric connections obtained from the composition of circulant
and path graphs.

\noindent In the former case, when the correlations of the noise,
the initial conditions and the synaptic weights are set to zero, we
have proved that propagation of chaos of the membrane potentials increases
with the number of incoming connections per neuron in the network,
and not with the number of neurons, as previously thought.

\noindent Instead, if the correlations of the three sources of randomness
are not set to zero, propagation of chaos in general does not occur,
even if the number of connections per neuron in the network is infinite.

\noindent Moreover, for special values of the parameters of the system,
we have proved that the membrane potentials become perfectly correlated,
a phenomenon that we have called \textit{stochastic synchronization}.

\noindent This phenomenon occurs, for any finite number of neurons,
even if the correlations of the three sources of randomness are set
to zero.

\noindent In this case the correlation of the membrane potentials
starts to increase from zero and it reaches the value $1$ (namely
perfect correlation) at a time instant that increases with the number
of neurons in the network.

\noindent Therefore for an infinite number of neurons, i.e. in the
mean-field limit of the system, correlation is always equal to zero
for any finite time instant.

\noindent This is in agreement with the results proved in \cite{journals/siamads/TouboulHF12}\cite{baladron:inserm-00732288}\cite{Touboul11}.

\selectlanguage{french}%
\noindent {\footnotesize{\smallskip{}
}}{\footnotesize \par}

\selectlanguage{english}%
\noindent Moreover, we have shown how to use these perturbative expansions
in the calculation of higher order correlations, like those between
triplets, quadruplets, quintuplets etc of neurons.

\noindent So with this approach we can determine all the moments of
the joint probability density, but in general we are not able to evaluate
the probability density itself.

\noindent In principle the moments can be used to calculate the moment-generating
function, from whom we could determine the corresponding probability
density, but in practice this is not feasible.

\noindent So we only know that at the first perturbative order the
process is normal, while at the second order it is a \textit{generalized
chi-squared process}.

\noindent The generalized chi-squared process can be defined as the
sum of products between pairs of normal processes, but its probability
density is not known, since it is still an open problem.

\noindent The same issue persists at higher orders, so another approach
must be followed.

\noindent Since the probability density is usually described by the
Fokker-Planck equation, we could think to solve it perturbatively.

\noindent However for a network with finite size we cannot use the
trick of the mean-field dimensional reduction, since here we have
proved that propagation of chaos does not occur, therefore the corresponding
Fokker-Planck equation would be high-dimensional.

\noindent Another problem with this approach is the second order derivative
of the equation, that describes the diffusion process.

\noindent This term is multiplied by the intensity of the Brownian
motions $\sigma_{1}$, forcing us to use the singular perturbation
theory.

\noindent So this idea is not promising, leaving the problem open.

\selectlanguage{french}%
\noindent {\footnotesize{\smallskip{}
}}{\footnotesize \par}

\selectlanguage{english}%
\noindent Even if we have developed this perturbative expansion for
a rate model with generic evolutions of the synaptic weights, we can
apply them to other kinds of models.

\noindent For example, we can consider spiking networks, described
for example by the FitzHugh-Nagumo or the Hodgkin-Huxley model.

\noindent In this case the perturbative expansion is useful only for
the description of the sub-threshold activity, since over the threshold
the neurons are spiking and so there is no stationary solution around
which the membrane potential can be expanded.

\noindent This is a different situation compared to the rate model,
since in this latter case also a stationary solution $\mu$ describes
a spiking activity, with rate $S\left(\mu\right)$.

\selectlanguage{french}%
\noindent {\footnotesize{\smallskip{}
}}{\footnotesize \par}

\selectlanguage{english}%
\noindent To conclude, we have proved that the perturbative method
developed in this article sheds new light in the comprehension of
stochastic neural networks, with special emphasis on finite size effects
and correlation.

\noindent In particular, these results establish the relation between
the functional and anatomical connectivity of the network, a problem
that is currently intensively investigated.

\noindent Moreover, at the first perturbative order, the probability
density of the system is a multivariate normal distribution.

\noindent For such a probability density, the information quanties
like the Shannon information, the Fisher information and the transfer
entropy \cite{schreiber2000information}, can be evaluated analytically
for any finite number of neurons.

\noindent This is a big advantage since otherwise these quantities
can be evaluated only through the calculation of high-dimensional
integrals using the Monte Carlo integration.

\noindent Therefore this allows us to quantify the information processing
capabilities of the neural networks, in terms of information encoding,
storage, transmission and modification, following the same ideas already
developed in the field of automata theory \cite{alifexi_lizier_374}\cite{liz07b}\cite{liz10e}\cite{liz12a}.

\noindent We think now we are in a much better position for understanding
in detail the working principles of stochastic neural networks.

\appendix
\begin{appendices}

\section{\noindent \label{sec:Radius of convergence of the sigmoid and arctangent functions}\noun{\index{Radius of convergence of the sigmoid and arctangent functions}}Radius
of convergence of the sigmoid and arctangent functions}

\noindent In this section we compute numerically the radius of convergence
of two examples of the activation function $S\left(\cdot\right)$.
For simplicity we consider only the case with $T_{MAX}=1$ and $V_{T}=0$,
but this analysis can be extended easily to the most general case.

\subsection{\label{sub:The sigmoid function}The sigmoid function}

\noindent According to \cite{Minai92originalcontribution}, the $n$-th
order derivative of the sigmoid function:

\begin{onehalfspace}
\begin{center}
{\small{
\[
S\left(x\right)=\frac{1}{1+e^{-\lambda x}}
\]
}}
\par\end{center}{\small \par}
\end{onehalfspace}

\noindent is:

\begin{onehalfspace}
\begin{center}
{\small{
\[
S^{\left(n\right)}\left(x\right)=\lambda^{n}\sum_{k=1}^{n}\left(-1\right)^{k-1}A\left(n,k-1\right)\left[S\left(x\right)\right]^{k}\left[1-S\left(x\right)\right]^{n+1-k}
\]
}}
\par\end{center}{\small \par}
\end{onehalfspace}

\noindent where $A\left(n,k\right)$ are the so called \textit{Eulerian
numbers} \cite{1959}. Now we can rewrite this expression in the following
way:

\begin{onehalfspace}
\begin{center}
{\small{
\begin{align*}
S^{\left(n\right)}\left(x\right)= & \lambda^{n}S\left(x\right)\left[1-S\left(x\right)\right]^{n}\sum_{k=1}^{n}\left(-1\right)^{k-1}A\left(n,k-1\right)\left[S\left(x\right)\right]^{k-1}\left[1-S\left(x\right)\right]^{-\left(k-1\right)}\\
\\
= & \lambda^{n}S\left(x\right)\left[1-S\left(x\right)\right]^{n}\sum_{k=0}^{n-1}\left(-1\right)^{k}A\left(n,k\right)\left[S\left(x\right)\right]^{k}\left[1-S\left(x\right)\right]^{-k}\\
\\
= & \lambda^{n}S\left(x\right)\left[1-S\left(x\right)\right]^{n}\sum_{k=0}^{n-1}A\left(n,k\right)\left(-e^{-\lambda x}\right)^{-k}
\end{align*}
}}
\par\end{center}{\small \par}
\end{onehalfspace}

\noindent Now from \cite{Miller:arXiv0804.3611}, we know that:

\begin{onehalfspace}
\begin{center}
{\small{
\begin{equation}
\begin{array}{ccc}
Li_{-n}\left(x\right)=\frac{x^{n}}{\left(1-x\right)^{n+1}}{\displaystyle \sum_{k=0}^{n-1}}A\left(n,k\right)x^{-k}, &  & n>0,\:\left|x\right|<1\end{array}\label{eq:polylogarithm-function}
\end{equation}
}}
\par\end{center}{\small \par}
\end{onehalfspace}

\noindent where $Li_{-n}\left(\cdot\right)$ represents the so called
\textit{polylogarithm} (with negative order). Here we have omitted
the $n$-th term of the sum since $A\left(n,n\right)=0$ $\forall n>0$.
So we can write:

\begin{onehalfspace}
\begin{center}
{\small{
\begin{equation}
S^{\left(n\right)}\left(x\right)=\lambda^{n}S\left(x\right)\left[1-S\left(x\right)\right]^{n}\frac{\left(1+e^{-\lambda x}\right)^{n+1}}{\left(-e^{-\lambda x}\right)^{n}}Li_{-n}\left(-e^{-\lambda x}\right)=\left(-\lambda\right)^{n}Li_{-n}\left(-e^{-\lambda x}\right)\label{eq:sigmoid-higher-order-derivatives}
\end{equation}
}}
\par\end{center}{\small \par}
\end{onehalfspace}

\noindent This result is true only for $\left|-e^{-\lambda x}\right|<1$,
i.e. only for $x>0$. Instead, for $x<0$, we can use the relation
$S\left(-x\right)=1-S\left(x\right)$, from which we deduce that:
\begin{itemize}
\item \noindent $S^{\left(n\right)}\left(-x\right)=\left(-1\right)^{n-1}S^{\left(n\right)}\left(x\right),$
$\forall n>0$;
\item \noindent $S\left(-x\right)$ has the same radius of convergence of
$S\left(x\right)$.
\end{itemize}
\noindent So formula \ref{eq:sigmoid-higher-order-derivatives} can
be used to express $S^{\left(n\right)}\left(x\right)$ $\forall x\neq0$.
Instead for $x=0$ it gives $Li_{-n}\left(-1\right)$, that is defined
by an analytic continuation of the polylogarithm function. In this
way we can determine $S^{\left(n\right)}\left(0\right)$. Another
way is to use the following property of the Eulerian numbers:

\begin{onehalfspace}
\begin{center}
{\small{
\begin{equation}
\sum_{k=1}^{n}\left(-1\right)^{k-1}A\left(n,k-1\right)=2^{n+1}\left(2^{n+1}-1\right)\frac{B_{n+1}}{n+1}\label{eq:Bernoulli-numbers}
\end{equation}
}}
\par\end{center}{\small \par}
\end{onehalfspace}

\noindent where $B_{n}$ are the so called \textit{Bernoulli numbers}
\cite{1991}, from which we obtain:

\begin{onehalfspace}
\begin{center}
{\small{
\begin{equation}
S^{\left(n\right)}\left(0\right)=\frac{\lambda^{n}}{2^{n+1}}\sum_{k=1}^{n}\left(-1\right)^{k-1}A\left(n,k-1\right)=\lambda^{n}\left(2^{n+1}-1\right)\frac{B_{n+1}}{n+1}\label{eq:sigmoid-higher-order-derivatives-in-zero}
\end{equation}
}}
\par\end{center}{\small \par}
\end{onehalfspace}

\noindent Now we can compute the radius of convergence $R\left(x_{0}\right)$
of the Taylor series:

\begin{onehalfspace}
\begin{center}
{\small{
\[
S\left(x\right)=\sum_{n=0}^{+\infty}\frac{S^{\left(n\right)}\left(x_{0}\right)}{n!}\left(x-x_{0}\right)^{n}
\]
}}
\par\end{center}{\small \par}
\end{onehalfspace}

\noindent using the Cauchy root test:

\begin{onehalfspace}
\begin{center}
{\small{
\[
R\left(x_{0}\right)=\frac{1}{\underset{n\rightarrow+\infty}{lim\, sup}\sqrt[n]{\left|\frac{S^{\left(n\right)}\left(x_{0}\right)}{n!}\right|}}
\]
}}
\par\end{center}{\small \par}
\end{onehalfspace}

\noindent For $x_{0}=0$ we obtain:

\begin{onehalfspace}
\begin{center}
{\small{
\[
R\left(0\right)=\frac{1}{\underset{n\rightarrow+\infty}{lim\, sup}\sqrt[n]{\left|\frac{\lambda^{n}\left(2^{n+1}-1\right)\frac{B_{n+1}}{n+1}}{n!}\right|}}=\frac{\pi}{\lambda}
\]
}}
\par\end{center}{\small \par}
\end{onehalfspace}

\noindent This can be proved after the substitution $n\rightarrow2n-1$
(which is motivated by the fact that $B_{2n+1}=0$ $\forall n>0$),
using the following asymptotic expansion of the Bernoulli numbers:

\begin{onehalfspace}
\begin{center}
{\small{
\[
\begin{array}{ccc}
B_{2n}\sim\left(-1\right)^{n-1}4\sqrt{\pi n}\left(\frac{n}{\pi e}\right)^{2n}, &  & n\rightarrow+\infty\end{array}
\]
}}
\par\end{center}{\small \par}
\end{onehalfspace}

\noindent and the Stirling approximation of $\left(2n-1\right)!$.
We are not aware of any asymptotic expansion of $Li_{-n}\left(-e^{-x_{0}}\right)$
for $n\rightarrow+\infty$ and $x_{0}\neq0$, so we have to compute
the radius of convergence numerically $\forall x_{0}\neq0$.

\begin{figure}
\begin{centering}
\includegraphics[scale=0.35]{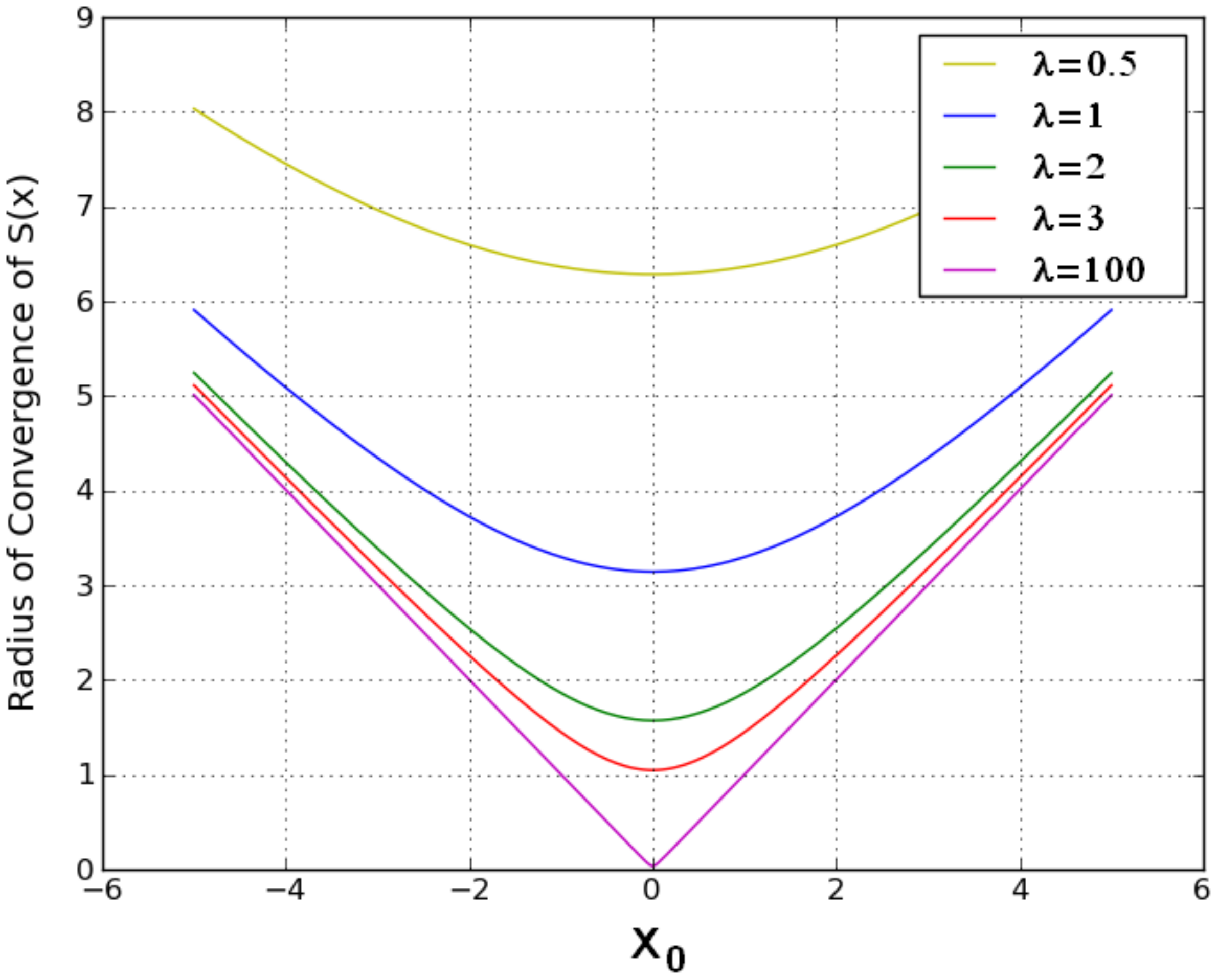}
\par\end{centering}

\caption[{\footnotesize{Radius of converge of the sigmoid function}}]{{\footnotesize{\label{fig:radius-convergence-sigmoid}Radius of convergence
$R$ of the Taylor series of the sigmoid function, in terms of the
point $x_{0}$ about which the expansion is performed. $R\left(x_{0}\right)$
has been computed numerically, for many values of the parameter $\lambda$,
which determines the slope of the sigmoid function. For large $x_{0}$
the radius of converge increases linearly since the sigmoid function
is asymptotically flat. Instead for $\lambda\rightarrow+\infty$ we
obtain $R\left(0\right)\rightarrow0$, because in that limit the sigmoid
function $S\left(x\right)$ becomes a Heaviside step function with
a discontinuity in $x=0$.}}}
\end{figure}

\noindent Figure \ref{fig:radius-convergence-sigmoid} shows the result
for different values of $\lambda$. From it we can see that the radius
of convergence of the Taylor series of $S\left(x\right)$ around the
point $x=x_{0}$ increases with $x_{0}$. This is reasonable, since
the function $S\left(x\right)$ becomes flat when $x$ is large. Moreover
for large $\lambda$ it converges to $R\left(x_{0}\right)=\left|x_{0}\right|$
and therefore it is equal to zero only for $x_{0}=0$, as it must
be. In fact, for $\lambda\rightarrow+\infty$ the function $S\left(x\right)$
converges to the Heaviside step function, which has a vertical jump
at $x=0$.

\noindent 

\subsection{\label{sub:The arctangent function}The arctangent function}

\noindent Now we calculate the radius of convergence of the arctangent
function. According to \cite{Adegoke2010}, the $n$-th order derivative
of this function is:

\begin{onehalfspace}
\begin{center}
{\small{
\[
arctan^{\left(n\right)}\left(\lambda x\right)=\lambda^{n}\frac{\left(-1\right)^{n-1}\left(n-1\right)!}{\left[1+\left(\lambda x\right)^{2}\right]^{\frac{n}{2}}}\sin\left[n\arcsin\left(\frac{1}{\sqrt{1+\left(\lambda x\right)^{2}}}\right)\right]
\]
}}
\par\end{center}{\small \par}
\end{onehalfspace}

\noindent So from the root test we obtain:

\begin{onehalfspace}
\begin{center}
{\small{
\[
R\left(x_{0}\right)=\frac{\sqrt{1+\left(\lambda x_{0}\right)^{2}}}{\lambda\underset{n\rightarrow+\infty}{lim\: sup}\frac{\sqrt[n]{\left|\sin\left[n\arcsin\left(\frac{1}{\sqrt{1+\left(\lambda x_{0}\right)^{2}}}\right)\right]\right|}}{\sqrt[n]{n}}}
\]
}}
\par\end{center}{\small \par}
\end{onehalfspace}

\noindent Now, since:

\begin{onehalfspace}
\begin{center}
{\small{
\[
\underset{n\rightarrow+\infty}{lim}\sqrt[n]{\left|\sin\left(n\arcsin\left(\frac{1}{\sqrt{1+\left(\lambda x_{0}\right)^{2}}}\right)\right)\right|}=1
\]
}}
\par\end{center}{\small \par}
\end{onehalfspace}

\noindent due to the fact that:

\begin{onehalfspace}
\begin{center}
\[
\left|\sin\left[n\arcsin\left(\frac{1}{\sqrt{1+\left(\lambda x_{0}\right)^{2}}}\right)\right]\right|\in\left[0,1\right]
\]

\par\end{center}
\end{onehalfspace}

\noindent and moreover$\underset{n\rightarrow+\infty}{lim}\sqrt[n]{n}=1$,
we obtain finally:

\begin{onehalfspace}
\begin{center}
{\small{
\[
R\left(x_{0}\right)=\frac{1}{\lambda}\sqrt{1+\left(\lambda x_{0}\right)^{2}}
\]
}}
\par\end{center}{\small \par}
\end{onehalfspace}

\noindent Therefore the radius of convergence increases with $x_{0}$,
as it must be. Moreover in the limit $\lambda\rightarrow+\infty$
it gives $R\left(x_{0}\right)=\left|x_{0}\right|$, as with the sigmoid
function.

\section{\noindent \label{sec:Higher order correlations for a fully connected neural network}\noun{\index{Higher order correlations for a fully connected neural network}}Higher
order correlations for a fully connected neural network}

\noindent Here we show how it is possible to use the perturbative
expansion to calculate the higher order correlations between the neurons.
For simplicity, we consider only the simplest case, namely a fully
connected network, even if this analysis could be extended to more
complicated connectivity matrices. Moreover we want to avoid long
expressions for the joint cumulants, therefore we consider only the
expansion of the membrane potential at the first perturbative order.
In principle this calculation can be performed at any perturbative
order, but starting from the second order (namely from the third order
terms in the covariance) the functions $Z\left(t\right)$ and $\overrightarrow{H}\left(t\right)$
in general introduce inhomogeneities in the covariance structure of
the network, therefore the higher order correlations should be calculated
using combinatorial techniques applied to the Isserlis' theorem. In
this section we avoid the issue and we focus only on the first order
perturbations. In this case the probability density of every $V_{i}\left(t\right)$
is normal and this is true also for the quantities $V_{i}\left(t\right)-\overline{V}_{i}\left(t\right)$,
which have all zero mean and the same variance, that we call $Var\left(V\left(t\right)\right)$.
Since they have zero mean we can use the Isserlis' theorem, that for
a fully connected network gives simply:

\begin{onehalfspace}
\begin{center}
{\small{
\begin{equation}
\mathbb{E}\left[\prod_{j=0}^{n-1}\left(V_{i_{j}}\left(t\right)-\overline{V}_{i_{j}}\left(t\right)\right)\right]=\begin{cases}
0, & \begin{array}{cc}
 & n\mathrm{\; odd}\end{array}\\
\\
\frac{n!}{2^{\frac{n}{2}}\left(\frac{n}{2}\right)!}\left[Cov\left(V_{i}\left(t\right),V_{j}\left(t\right)\right)\right]^{\frac{n}{2}}, & \begin{array}{cc}
 & n\mathrm{\; even}\end{array}
\end{cases}\label{eq:Isserlis-theorem}
\end{equation}
}}
\par\end{center}{\small \par}
\end{onehalfspace}

\noindent because in this case all the pairs of neurons are equivalent,
since they are all-to-all connected (instead, if the network is not
fully connected, the connected pairs give a different contribution
with the Isserlis' theorem compared to the disconnected pairs). Moreover,
the central absolute moments of a normal distribution are:

\begin{onehalfspace}
\begin{center}
{\small{
\begin{equation}
\mathbb{E}\left[\left|V_{i_{j}}\left(t\right)-\overline{V}_{i_{j}}\left(t\right)\right|^{n}\right]=\frac{2^{\frac{n}{2}}\Gamma\left(\frac{n+1}{2}\right)}{\sqrt{\pi}}\left[Var\left(V\left(t\right)\right)\right]^{\frac{n}{2}}\label{eq:normal-central-absolute-moments}
\end{equation}
}}
\par\end{center}{\small \par}
\end{onehalfspace}

\noindent and if $n$ is even we have $\Gamma\left(\frac{n+1}{2}\right)=\frac{n!}{2^{n}\left(\frac{n}{2}\right)!}\sqrt{\pi}$.
Therefore putting everything together we obtain:

\begin{onehalfspace}
\begin{center}
{\small{
\begin{equation}
Corr_{n}\left(V_{i_{0}}\left(t\right),V_{i_{1}}\left(t\right),...,V_{i_{n-1}}\left(t\right)\right)=\begin{cases}
0, & \begin{array}{cc}
 & n\mathrm{\; odd}\end{array}\\
\\
\left[Corr_{2}\left(V_{i}\left(t\right),V_{j}\left(t\right)\right)\right]^{\frac{n}{2}}, & \begin{array}{cc}
 & n\mathrm{\; even}\end{array}
\end{cases}\label{eq:higher-order-correlation-2}
\end{equation}
}}
\par\end{center}{\small \par}
\end{onehalfspace}

\noindent From this result it is interesting to observe that if there
is a perfect stochastic synchronization between pairs of neurons,
then it is ``propagated'' to all the higher order correlations with
even order, namely $Corr_{2}\left(V_{i}\left(t\right),V_{j}\left(t\right)\right)=1$
implies $Corr_{n}\left(V_{i_{0}}\left(t\right),V_{i_{1}}\left(t\right),...,V_{i_{n-1}}\left(t\right)\right)=1$,
$\forall n$ even. It is also curious to observe that all the odd
order correlations are always equal to zero, even if the size of the
network is finite.

\section{\noindent \label{sec:Proof that formula 4.6 gives real functions}\noun{\index{Proof that formula 4.6 gives real functions}}Proof
that formula 4.6 gives real functions}

\begin{flushleft}
In this appendix we want to prove that the quantities $\Phi_{ij}\left(t\right)$
and $\left[\Phi\left(t\right)\Phi^{T}\left(t\right)\right]_{ij}$,
given by formula \ref{eq:Phi-matrix-block-circulant-case}, are real
functions. The proof can be divided into four cases, namely when $R$
and $S$ are both even, both odd, $R$ even and $S$ odd, or vice
versa. Here we analyze only the first case, while the others can be
proved in a similar way.
\par\end{flushleft}

\noindent So, if $R$ and $S$ are both even, the function $\Phi_{ij}\left(t\right)$,
according to \ref{eq:Phi-matrix-block-circulant-case}, can be equivalently
rewritten as:

\begin{onehalfspace}
\begin{center}
{\small{
\begin{equation}
\Phi_{ij}\left(t\right)={\displaystyle \sum_{x=0}^{R-1}}{\displaystyle \sum_{y=0}^{S-1}}e^{\left[-\frac{1}{\tau}+e_{xS+y}S'\left(\mu\right)\right]t}f_{i,j,xS+y}\label{eq:Phi-equivalent-form-1}
\end{equation}
}}
\par\end{center}{\small \par}
\end{onehalfspace}

\noindent where now the subscripts are separated by commas, in order
to avoid confusion. Defining:

\begin{onehalfspace}
\begin{center}
{\small{
\[
g_{x,y}^{ij}=e^{\left[-\frac{1}{\tau}+e_{xS+y}S'\left(\mu\right)\right]t}f_{i,j,xS+y}
\]
}}
\par\end{center}{\small \par}
\end{onehalfspace}

\noindent formula \ref{eq:Phi-equivalent-form-1} can be rewritten
in the following symmetric way, with respect to $R$ and $S$:

\begin{onehalfspace}
\begin{center}
{\small{
\begin{align}
\Phi_{ij}\left(t\right)= & g_{0,0}^{ij}+g_{\frac{R}{2},0}^{ij}+g_{0,\frac{S}{2}}^{ij}+g_{\frac{R}{2},\frac{S}{2}}^{ij}\nonumber \\
\nonumber \\
 & +{\displaystyle \sum_{x=1}^{\frac{R}{2}-1}}\left[g_{x,0}^{ij}+g_{R-x,0}^{ij}\right]+{\displaystyle \sum_{y=1}^{\frac{S}{2}-1}}\left[g_{0,y}^{ij}+g_{0,S-y}^{ij}\right]+{\displaystyle \sum_{y=1}^{\frac{S}{2}-1}}\left[g_{\frac{R}{2},y}^{ij}+g_{\frac{R}{2},S-y}^{ij}\right]+{\displaystyle \sum_{x=1}^{\frac{R}{2}-1}}\left[g_{x,\frac{S}{2}}^{ij}+g_{R-x,\frac{S}{2}}^{ij}\right]\nonumber \\
\nonumber \\
 & +{\displaystyle \sum_{x=1}^{\frac{R}{2}-1}}{\displaystyle \sum_{y=1}^{\frac{S}{2}-1}}\left[g_{x,y}^{ij}+g_{R-x,S-y}^{ij}\right]+{\displaystyle \sum_{x=1}^{\frac{R}{2}-1}}{\displaystyle \sum_{y=1}^{\frac{S}{2}-1}}\left[g_{x,S-y}^{ij}+g_{R-x,y}^{ij}\right]\label{eq:Phi-equivalent-form-2}
\end{align}
}}
\par\end{center}{\small \par}
\end{onehalfspace}

\noindent The quantities $g_{0,0}^{ij}$, $g_{\frac{R}{2},0}^{ij}$,
$g_{0,\frac{S}{2}}^{ij}$ and $g_{\frac{R}{2},\frac{S}{2}}^{ij}$
are real numbers. Moreover, since:

\begin{onehalfspace}
\begin{center}
{\small{
\begin{align*}
f_{i,j,xS+y}= & e^{2\pi\left\{ \frac{x}{R}\left(\left\lfloor \frac{i}{S}\right\rfloor -\left\lfloor \frac{j}{S}\right\rfloor \right)+\frac{y}{S}\left(i-j\right)\right\} \iota}=f_{i,j,\left(R-x\right)S+\left(S-y\right)}^{*}\\
\\
f_{i,j,xS+\left(S-y\right)}= & e^{2\pi\left\{ \frac{x}{R}\left(\left\lfloor \frac{i}{S}\right\rfloor -\left\lfloor \frac{j}{S}\right\rfloor \right)-\frac{y}{S}\left(i-j\right)\right\} \iota}=f_{i,j,\left(R-x\right)S+y}^{*}
\end{align*}
}}
\par\end{center}{\small \par}
\end{onehalfspace}

\noindent for $0\leq y<S$ and $0\leq x<R$, and also, according to
\ref{eq:block-circulant-matrix-eigenvalues}:

\begin{onehalfspace}
\begin{center}
{\small{
\begin{align*}
e_{xS+y}= & \sum_{k=0}^{S-1}\sum_{l=0}^{R-1}e^{2\pi\left(\frac{yk}{S}+\frac{xl}{R}\right)\iota}b_{k}^{\left(l\right)}=e_{\left(R-x\right)S+\left(S-y\right)}^{*}\\
\\
e_{xS+\left(S-y\right)}= & \sum_{k=0}^{S-1}\sum_{l=0}^{R-1}e^{2\pi\left(-\frac{yk}{S}+\frac{xl}{R}\right)\iota}b_{k}^{\left(l\right)}=e_{\left(R-x\right)S+y}^{*}
\end{align*}
}}
\par\end{center}{\small \par}
\end{onehalfspace}

\noindent we conclude that:

\begin{onehalfspace}
\begin{center}
{\small{
\begin{align*}
g_{x,y}^{ij}= & \left(g_{R-x,S-y}^{ij}\right)^{*}\\
\\
g_{x,S-y}^{ij}= & \left(g_{R-x,y}^{ij}\right)^{*}
\end{align*}
}}
\par\end{center}{\small \par}
\end{onehalfspace}

\noindent For this reason, all the quantities in the square parenthesis
in formula \ref{eq:Phi-equivalent-form-2} are real, and therefore
also $\Phi_{ij}\left(t\right)$. A similar proof can be obtained for
$\left[\Phi\left(t\right)\Phi^{T}\left(t\right)\right]_{ij}$, and
in the cases when only one of $R$ and $S$, or both, are odd.

\end{appendices}

\section*{Acknowledgements}

This work was partially supported by the ERC grant \#227747 NerVi,
the FACETS-ITN Marie-Curie Initial Training Network \#237955 and the
IP project BrainScaleS \#269921.

\selectlanguage{french}%
\noindent {\footnotesize{\smallskip{}
}}{\footnotesize \par}

\selectlanguage{english}%
\bibliographystyle{unsrt}
\bibliography{Article}

\end{document}